\documentclass[12pt,a4paper]{article}
\usepackage{a4wide}
\usepackage[T1]{fontenc}
\usepackage[utf8]{inputenc}
\usepackage{lmodern}
\usepackage[USenglish]{babel}
\usepackage{microtype}
\usepackage{xcolor}
\usepackage{csquotes}
\usepackage{amsmath}
\usepackage{mathtools}
\usepackage{amssymb}
\usepackage{amsthm}

\usepackage{graphicx}
\usepackage{caption}
\usepackage{subcaption}
\usepackage{rotating}
\usepackage{float}
\usepackage[section]{placeins}

\usepackage{tikz}
\usetikzlibrary{calc,math,positioning,shapes,arrows,snakes,decorations.text}

\usepackage{pgfplots}
\pgfplotsset{compat=1.18} 
\usepgfplotslibrary{groupplots} 
\usepackage{pgfplotstable}

\usepackage{hyperref}

\newcommand{\HeatSmoothAxisHeight}{0.27}
\newcommand{\HeatSingularAxisHeight}{0.29}
\newcommand{\HeatGMRESAxisHeight}{0.36}
\newcommand{\ODELegendColumns}{3}

\newcommand{\PlaceODELegend}{%
  \node[anchor=north] at ($(group c1r1.south west)!0.5!(group c2r1.south east)+(0,-1.25cm)$)
    {\pgfplotslegendfromname{commonlegend}};
}

\DeclareMathOperator{\dom}{dom}

\DeclareMathOperator{\dist}{dist}

\newcommand{\R}{\ensuremath{\mathbb{R}}}

\newcommand{\N}{\ensuremath{\mathbb{N}}}

\newcommand{\Hmod}{\ensuremath{\mathcal{H}_T}}
\newcommand{\Hstd}{\ensuremath{\mathcal{H}}}
\newcommand{\Hcl}{\ensuremath{\mathcal{H}_{\mathrm{cl}}}}
\newcommand{\Hmu}{\ensuremath{\mathcal{H}_\mu}}

\renewcommand{\hat}[1]{\widehat{#1}}
\renewcommand{\tilde}[1]{\widetilde{#1}}

\newcommand{\norm}[1]{\ensuremath{\|#1\|}}
\newcommand{\skpr}[1]{\ensuremath{\langle#1\rangle}}

\newtheorem{theorem}{Theorem}

\newtheorem{lemma}[theorem]{Lemma}
\newtheorem{corollary}[theorem]{Corollary}

\newtheorem{definition}[theorem]{Definition}

\newtheorem{remark}[theorem]{Remark}

\hypersetup{
 pdftitle={Coercivity by Design: A Unified Approach to Parabolic Variational Formulations},
 pdfauthor={Richard Löscher and Michael Reichelt},
 hidelinks
}

\title{Coercivity by Design\\[0.3em]
A Unified Approach to Parabolic Variational Formulations}
\author{%
 \begin{tabular}{c}
  Richard L\"oscher$^{1}$, Michael Reichelt$^{2}$\\[0.8em]
  $^{1}$Graz University of Technology\\
  Institute of Applied Mathematics\\[0.5em]
  $^{2}$Helmut Schmidt University\\
  Faculty of Mechanical Engineering
 \end{tabular}%
}
\date{}

\begin{document}
\pagestyle{plain}
\maketitle

\begin{abstract}
Stable discretizations of inf–sup stable problems generally require compatible trial and test spaces, whereas coercive formulations preserve stability on every conforming Galerkin subspace. We present a unified framework based on the polar decomposition directly on the operator level that systematically constructs coercive formulations and connects existing approaches for parabolic problems. For the heat equation, it yields a test transformation acting jointly in space and time and a symmetric bilinear form coercive in the full energy norm. Since the exact transformation is computationally impractical, we approximate it by the sum of the identity and the modified Hilbert transform. We prove energy error estimates for tensor-product finite element discretizations. Numerical examples complement the theory.
\end{abstract}


\section{Introduction}
A standard approach to the numerical solution of parabolic problems
combines a spatial discretization with a time-stepping scheme
\cite{WannerHairerNorsett:}. Many classical time-stepping schemes also
admit a variational interpretation in time, for example as
Petrov--Galerkin methods \cite{AndreevSchweitzer:2014}. This perspective
connects time integration with finite element analysis and provides
a common framework for studying stability, approximation properties,
and discretization errors.

Also, the analysis of evolution problems is more commonly posed on the separation of spatial and temporal parts of the operator, e.g. by semi-group theory \cite{EngelNagelOneparameterSemigroups2000}. Space--time methods extend this perspective by formulating the
parabolic problem as a single variational operator equation over the
entire space--time domain. This allows approximation spaces to be
designed jointly in space and time and opens the way to simultaneous
space--time adaptivity. Suitable tensor-product discretizations
recover familiar time-stepping schemes
\cite{Andreev:2012,LoescherReicheltSteinbachTimeStepping:2026}.
Their space--time interpretation provides additional insight into
stability and convergence, while also offering a framework for
developing more general discretizations.

A central question in the design of such methods is whether
continuous stability is preserved under discretization. For a bounded,
coercive bilinear form, every conforming Galerkin subspace inherits
the continuous coercivity bound, and standard Galerkin theory yields
quasi-optimal error estimates. Coercivity can also facilitate the
construction of effective preconditioners, as the induced energy norm
provides a natural reference for their design and analysis.
Continuous inf--sup stability, however, need not be inherited by
arbitrary discrete trial and test spaces. A coercive formulation
therefore simplifies the construction of stable Galerkin methods.
Conforming approximation spaces can be chosen without imposing an
additional compatibility condition between discrete trial and test
spaces.

For parabolic problems, coercivity in the desired space--time norm
is not immediate from the standard variational formulation.
We pursue a systematic approach to obtaining such coercivity through
the polar decomposition. For the operators considered here, this
decomposition separates the operator into a positive self-adjoint
factor and a unitary factor. The quadratic form associated with the
positive factor defines an energy norm, while the unitary factor
determines a transformation of the test functions. The resulting
bilinear form is symmetric and coercive in this norm. Consequently,
a conforming trial space in the energy space yields a stable
Petrov--Galerkin method when paired with its image under the test
transformation. This principle is related to the construction of
optimal test spaces in discontinuous Petrov--Galerkin (DPG) methods
\cite{DemkowiczGopalakrishnan:} and to the use of transformation
operators in $T$-coercivity
\cite{CiarletTCoercivityApplication2012}.

The variational treatment of parabolic problems in fractional
Sobolev spaces has its foundations in the work of Lions and Magenes
\cite{LionsMagenesVol1}. Temporal regularity of order one half allows
a balanced treatment of trial and test functions, with the time
derivative acting between spaces of equal temporal regularity.
Early discretization error estimates in fractional Sobolev norms
were established by Baiocchi and Brezzi \cite{BrezziBaiocchi:}.

Steinbach and Zank \cite{ZankSteinbachCoerciveSpaceTimeFinite2020}
provided an explicit temporal test transformation on a finite
interval by introducing the modified Hilbert transformation.
This transformation makes the time-derivative term symmetric and
coercive and leads to conforming space--time discretizations of
the heat equation. The resulting coercivity estimate directly
controls the temporal half-order norm. This construction provides
a natural starting point for incorporating the spatial contribution
into the same coercivity estimate.

We build on this transformation-based approach by applying the
polar decomposition to the full parabolic operator. The resulting
transformation acts jointly in space and time and yields a symmetric
formulation that is coercive in the full energy norm. We establish
stability and discretization error estimates and explain the
relationship between the proposed transformation and the modified
Hilbert transformation.

The paper is organized as follows. Section~2 reviews the polar
decomposition in the operator setting used throughout the paper.
Section~3 studies the time derivative and develops the tools needed
in Section~4, where we consider a parameter-dependent ordinary
differential equation associated with the heat equation. We analyze
its polar decomposition, derive discretization error estimates,
and illustrate the results numerically. Section~5 extends the
analysis to the heat equation. Numerical experiments with
tensor-product discretizations illustrate the approximation
properties and computational performance of the proposed method,
and preconditioning strategies are discussed.
\section{The polar decomposition}
In this section, we introduce a functional-analytic framework for
constructing transformation operators that yield coercive variational
formulations. 
As a first illustrative example, that is used to get familiar with the notation, we consider the well-known Dirichlet Laplacian on a bounded domain
$\Omega\subset\mathbb{R}^d$. As an unbounded operator on $L^2(\Omega)$, it is
defined by
\begin{equation*}
    -\Delta:\dom(-\Delta)\subset L^2(\Omega)\to L^2(\Omega),
    \qquad
    \dom(-\Delta)
    :=
    \bigl\{
        u\in H_0^1(\Omega):
        \Delta u\in L^2(\Omega)
    \bigr\},
\end{equation*}
where the Laplacian is understood in the distributional sense. For ease of presentation let us assume that $\Omega$ is either convex or smoothly bounded, then 
\begin{equation*}
    \dom(-\Delta)
    =H_\Delta^1(\Omega) 
    := \{v\in H^1_0(\Omega):\, \Delta v\in L^2(\Omega)\}
    = H^2(\Omega)\cap H_0^1(\Omega).
\end{equation*}
For given $f\in L^2(\Omega)$, the strong operator formulation is then to find $u\in \dom(A)$ such that 
\begin{equation*}
    -\Delta u=f
    \qquad\text{in }L^2(\Omega).
\end{equation*}
The usual variational formulation is instead derived by multiplying by a test function and integrating by parts, and reads to find $u\in H^1_0(\Omega)$ such that 
\begin{equation*}
    a(u,v):=
    \int_\Omega \nabla u(x)\cdot\nabla v(x)\,\mathrm{d}x 
    =
    \int_\Omega f(x)v(x)\, \mathrm{d}x, 
    \quad 
    \forall v\in H^1_0(\Omega). 
\end{equation*}
Thus, the variational problem is posed in the larger space on the domain of the form
$\dom(a) = H_0^1(\Omega)$, called the \emph{energy space}. This raises the natural question on how the domain of the operator $\dom(-\Delta)$ is related to the domain of the form $\dom(a)$. 
Indeed, note that the Dirichlet Laplacian is a densely defined, strictly positive, self-adjoint operator on $L^2(\Omega)$. Its positive square root $(-\Delta)^{1/2}$ is therefore well-defined, and $\dom((-\Delta)^{1/2}) = H^1_0(\Omega)$, see \cite[Chapter~VI, Theorem~2.23]
{KatoPerturbationTheoryLinear1995}. Thus there holds the relation 
\begin{equation*}
    \dom(a) 
    = 
    H^1_0(\Omega) 
    = 
    [L^2(\Omega),H^2(\Omega)\cap H^1_0(\Omega)]_{1/2}
    =
    \dom((-\Delta)^{1/2}). 
\end{equation*} 
Hence, the natural energy space is the domain of the form, which is exactly the domain of the square root of the operator.  
The corresponding \emph{energy norm} is given by 
\begin{equation*}
    \norm{u}_a = a(u,u)^{1/2} = \norm{\nabla u}_{L^2(\Omega)}. 
\end{equation*}
The following theorem states this principle in an abstract setting. 

\begin{theorem}\label{thm:coercive for self adjoint operators}
Let $H$ be a Hilbert space and let
$A\colon\dom(A)\subset H\to H$ be a densely defined, positive,
self-adjoint operator. Let $a(\cdot,\cdot)$ denote the closed symmetric form
associated with $A$, i.e., $a(u,v) = (Au,v)_H$ for all $u,v\in \dom(A)$. Then the \emph{energy space} is denoted by 
\begin{equation*}
    X 
    :=\dom(a)
    =\dom(A^{1/2}) 
    = [H,\dom(A)]_{1/2},
\end{equation*}
where $\dom(A)$ is equipped with its graph norm and the equality is
understood with equivalent norms.
Consequently, $(X,\norm{\cdot}_X)$ where 
\begin{equation*}
    \|u\|_X^2
    :=
    \|u\|_H^2+a(u,u),
\end{equation*}
is a Hilbert space densely and continuously embedded into $H$, and identifying $H$ with its dual gives the Gelfand triple
$
    X\hookrightarrow H\cong H'\hookrightarrow X'.
$

In particular, if $A$ is strictly positive, that is, there exists $\alpha>0$ such that 
$(Av,v)_H\geq\alpha \norm{v}_H^2$ for all $v\in \dom(A)$, then the \emph{energy norm} is equivalent to 
\begin{equation*}
    \|u\|_X=a(u,u)^{1/2}. 
\end{equation*}
\end{theorem}
\begin{proof}
    The definition of the interpolation space, $\dom(A^{1/2}) = [H,\dom(A)]_{1/2}$ can be found, e.g., in \cite[Ch. 1, Sec. 2.1]{LionsMagenesVol1} and \cite[Chapter 4, Proposition 2.2]{Taylor2011}.
    The relation $\dom(a) = \dom(A^{1/2})$ is established in \cite[Chapter VI, Thm. 2.23] {KatoPerturbationTheoryLinear1995}, together with 
    \begin{equation*}
        a(u,v) = (A^{1/2}u,A^{1/2}v),\quad \forall u,v\in \dom(a). 
    \end{equation*}
    Thus, indeed if $A \geq \alpha I$ it holds for any $v\in\dom(a)$ 
    \begin{equation*}
        a(v,v) 
        = 
        \norm{A^{1/2}v}_H^2
        \geq 
        \alpha \norm{v}_H^2,
    \end{equation*}
    giving the desired norm equivalence. Since $a$ is closed and densely defined, $\dom(a)\subset H$ equipped with $\norm{\cdot}_X$ is a Hilbert space, continuously embedded in $H$. Duality then yields the Gelfand triple.   
\end{proof}
\begin{remark}\label{rem:coercive operator quasi best approxiamtion}
    The preceding theorem implies that for any densely defined, strictly positive, and self-adjoint operator $A:\dom(A)\subset H\to H$, the associated bilinear form is coercive and bounded  
    \begin{equation*}
        \norm{v}_X^2\leq a(v,v),\quad 
        a(u,v) \leq \norm{u}_X\norm{v}_X, \quad 
        \forall u,v\in X. 
    \end{equation*}
    Thus, for any $f\in X'$ the variational formulation to find $u\in X$ such that 
    \begin{equation*}
        a(u,v) = \skpr{f,v}_{X'\times X},\quad \forall v\in X,
    \end{equation*}
    admits a unique solution by the Lemma of Lax-Milgram \cite[Lemma 2.2]{ErnGuermondTheoryPracticeFinite2004}. In particular, for any conforming finite dimensional subspace $X_h\subset X$, the discrete problem to find $u_h\in X_h$ such that 
    \begin{equation*}
        a(u_h,v_h) = \skpr{f,v_h}_{X'\times X},\quad \forall v_h\in X_h,
    \end{equation*}
    admits a unique solution, satisfying the (quasi)-best approximation error estimate \cite[Lemma 2.28]{ErnGuermondTheoryPracticeFinite2004}
    \begin{equation*}
        \norm{u-u_h}_X \leq \inf_{v_h\in X_h}\norm{u-v_h}_X. 
    \end{equation*} 
\end{remark}

Though, the picture changes, if $A$ is not self-adjoint. To apply the same concept also in this case, we will use the \emph{polar decompostion}. 

Before we state the general result for unbounded linear operators, let us consider an invertible matrix $\mathbf{B}\in \mathbb{R}^{n\times n}.$
The singular value decomposition gives orthogonal matrices $\mathbf{U},\mathbf{V}\in \mathbb{R}^{n\times n}$ and a diagonal matrix $\mathbf{\Sigma}\in \mathbb{R}^{n\times n}_{>0}$ such that 
\begin{equation*}
    \mathbf{B} = \mathbf{U}\mathbf{\Sigma} \mathbf{V}^\top.
\end{equation*}
In particular, we compute $\mathbf{B}^\top \mathbf{B} = \mathbf{V}\mathbf{\Sigma}^2\mathbf{V}^\top$, which is symmetric and positive definite. Since all the entries of $\mathbf{\Sigma}$ are positive, taking the square root gives $(\mathbf{B}^\top \mathbf{B})^{1/2} = \mathbf{V}\mathbf{\Sigma}\mathbf{V}^\top.$ Consequently, we can derive the \emph{polar decompostion} given as  
\begin{equation}\label{eq:polar decomposition matrix}
    \mathbf{B} 
    =
    \mathbf{U}(\mathbf{V}^\top \mathbf{V}) \mathbf{\Sigma}\mathbf{V}^\top
    = 
    \mathbf{U}\mathbf{V}^\top (\mathbf{B}^\top \mathbf{B})^{1/2}
    =
    \mathbf{H}|\mathbf{B}|,  
\end{equation}
where $\mathbf{H} := \mathbf{U}\mathbf{V}^\top$ is an orthogonal matrix and $|\mathbf{B}| := (\mathbf{B}^\top \mathbf{B})^{1/2}$ is symmetric and positive definite. 
The same decomposition holds true for operators. 

\begin{theorem}\label{thm:polar_decomposition}
    Let $H$ be a Hilbert space and $B:\dom(B)\subset H \to H$ be a closed and densely defined operator, satisfying $\ker(B) = \ker(B^\ast) = \{0\}.$ Then there exists a unique unitary operator $\mathcal{H}:H\to H$ such that 
    \begin{equation}\label{eq:polar decomposition B}
        B = \mathcal{H} \lvert B \rvert,
        \quad 
        |B| = (B^\ast B)^{1/2}. 
    \end{equation}
    Moreover, $|B|$ is positive and $\dom(B) = \dom(|B|)$.
    
    In particular, if the there exists $\beta>0$ such that the inf-sup condition  
    \begin{equation}\label{eq:inf-sup B}
        \beta \leq \inf_{0\neq u\in \dom(B)} \sup_{0\neq v\in H} \frac{(Bu,v)_H}{\norm{u}_H\norm{v}_H},
    \end{equation}
    is satisfied, then $|B|$ is strictly positive, i.e., 
    \begin{equation*}
        (|B|v,v)_H \geq \beta \norm{v}_H^2,\quad \forall v\in \dom(B). 
    \end{equation*}
\end{theorem}
\begin{proof}
        A proof of the existence and uniqueness of $\mathcal{H}$, as well as the fact that $\dom(B) = \dom(|B|)$ and $|B|$ is positive is given in \cite[Section VIII.32]{ReedSimonMethodsModernMathematicalPhysics1980}. Further, note that for any $u\in \dom(B)$ using the Cauchy--Schwarz inequality we have that 
        \begin{equation*}
          \sup_{0\neq v\in H} \frac{(Bu,v)_H}{\norm{v}_H} 
          \leq  \sup_{0\neq v\in H} \frac{\norm{Bu}_H\norm{v}_H}{\norm{v}_H}
          = \norm{Bu}_H. 
        \end{equation*} 
        On the other hand, choosing $\widehat v = Bu\in H$ we get 
        \begin{equation*}
            \sup_{0\neq v\in H}\frac{(Bu,v)_H}{\norm{v}_H}
            \geq 
            \frac{(Bu,\widehat{v})_H}{\norm{\widehat{v}}_H} 
            = 
            \frac{\norm{Bu}_H^2}{\norm{\widehat{v}}_H} 
            = \norm{Bu}_H, 
        \end{equation*}
        and thus we obtain the relation 
        \begin{equation*}
            \norm{Bu}_H = \sup_{0\neq v\in H} \frac{(Bu,v)_H}{\norm{v}_H}. 
        \end{equation*}
        Now if \eqref{eq:inf-sup B} holds, using that $\mathcal{H}$ is isometric and $|B|$ is self-adjoint, we get  
        \begin{align*}
            \beta^2 \norm{u}_H^2
            &\leq 
            \left(\sup_{0\neq v\in H} \frac{(Bu,v)_H}{\norm{v}_H}\right)^2
            =
            \norm{Bu}_H^2 
            = 
            \norm{\mathcal{H}|B|u}_H^2 \\
            &= 
            \norm{|B|u}_H^2.  
        \end{align*}
        Thus $|B| \geq \beta I $ and therefore $(|B|u,u)_H \geq \beta\norm{u}_H^2$.  
    \end{proof}
    Whenever \eqref{eq:inf-sup B} is satisfied, $|B|:\dom(|B|)\subset H\to H$ is densely defined, strictly positive, and self-adjoint. By Theorem \ref{thm:coercive for self adjoint operators}, we can thus define the \emph{energy space} in this setting by 
    \begin{equation*}
        X := \dom(|B|^{1/2}) = [H,\dom(B)]_{1/2},
    \end{equation*}
    equipped with the \emph{energy norm}
    \begin{equation*}
        \norm{v}_X := \norm{|B|^{1/2}v}_H = \sqrt{\skpr{|B|v,v}_{X'\times X}},\quad \forall v\in X. 
    \end{equation*}
    Note that by construction $X\hookrightarrow H$ is dense and $|B|:X\to X'$ is an isomorphism, by the Lemma of Lax--Milgram. 
    \begin{lemma}\label{lem:transformed energy space}
        Let $\Hstd:H\to H$ be unitary and $X\hookrightarrow H$ be dense. Then 
        \begin{equation*}
            Y := \Hstd X = \{\Hstd v:\, v\in X\}, 
        \end{equation*}   
        equipped with the inner product 
        \begin{equation*}
            (p,q)_Y 
            = 
            (\Hstd^\ast p,\Hstd^\ast q)_X, 
            \quad \forall p,q\in Y,
        \end{equation*}     
        is a Hilbert space, continuously embedded in $H$. Moreover, $\Hstd:X\to Y$ is unitary. 
    \end{lemma}
    \begin{proof}
        By construction,
        \begin{equation*}
            (\Hstd u,\Hstd v)_Y=(u,v)_X,
            \qquad u,v\in X,
        \end{equation*}
        and $\Hstd|_X$ is bijective from $X$ onto $Y$. Hence it is
        unitary, and the completeness of $Y$ follows from that of $X$. Moreover, since $X\hookrightarrow H$ we have for any $u\in X$ that  
        $\norm{u}_H\leq C \norm{u}_X$ and for any $q\in Y$ we compute 
        \begin{equation*}
            \norm{q}_H
            =\norm{\Hstd^\ast q}_H
            \leq C\norm{\Hstd^\ast q}_X
            =C\norm{q}_Y, 
        \end{equation*}
        and thus $Y\hookrightarrow H$ continuously. Density follows from 
        \begin{equation*}
            \overline{Y}^{\norm{\cdot}_H} 
            = 
            \overline{\Hstd|_X X}^{\norm{\cdot}_H} 
            = 
            \Hstd \overline{X}^{\norm{\cdot}_H}
            =
            \Hstd H 
            = H.  
        \end{equation*} 
    \end{proof}

    Let us now turn to the operator equation, to find $u\in \dom(B)$ such that 
    \begin{equation}\label{eq:operator equation non self adjoint}
        B u = f\quad \text{in }H,
    \end{equation}
    where $B:\dom(B)\subset H\to H$ is densely defined, closed and satisfies the inf-sup condition \eqref{eq:inf-sup B}.
    Multiplying with $\Hstd^\ast$ gives 
    \begin{equation*}
        |B| u = \Hstd^\ast B u = \Hstd^\ast f \quad \text{in }H. 
    \end{equation*}
    With the associated form, we obtain the variational formulation to find $u\in \dom(|b|) = \dom(|B|^{1/2}) = X$ such that 
    \begin{equation}\label{eq:VF symmetrized for B}
        |b|(u,v) 
        :=
        \skpr{|B|u,v}_{X'\times X} 
        = 
        \skpr{\Hstd^\ast f,v}_{X'\times X}
        = 
        \skpr{f,\Hstd v}_{Y'\times Y},
        \quad 
        \forall v\in X, 
    \end{equation} 
    which is well-posed for any $f\in Y'$. Moreover, $|b|(\cdot,\cdot):X\times X\to \mathbb{R}$ is symmetric, coercive and bounded, and by the Lemma of Lax--Milgram \eqref{eq:VF symmetrized for B} admits a unique solution, satisfying 
    \begin{equation*}
        \norm{u}_X \leq \norm{f}_{Y'}. 
    \end{equation*}   
    Hence, using the polar decomposition, we can define the natural energy space also for non self-adjoint operators and reformulate the operator equation \eqref{eq:operator equation non self adjoint} into a symmetric, and coercive variational formulation.

\section{The Time Derivative and its Polar Decomposition}
In this chapter, we apply the abstract framework to the first-derivative operator in one dimension. We derive its polar decomposition and show that the resulting form domains and unitary factors are classical, well-studied objects. We first consider the derivative on an unbounded interval and then extend the analysis to a bounded interval. This reveals a fundamental connection between the classical Hilbert transformation and the modified Hilbert transform, see also \cite{FerrariSomePropertiesModified2024,SteinbachMissoni2023}. 

\subsection{The derivative on an unbounded interval}

Let us consider the derivative as an unbounded operator
\begin{equation*}
  B \, : \, \dom(B)\subset L^2(\R)  \to L^2(\R),\quad  
  B v = \partial_t v, 
\end{equation*}
with domain
$\dom(B)= H^1(\R).$
We can proof the following mapping properties.
\begin{lemma}\label{lem:unbounded derivative}
  The operator $B$ is densely defined, closed and $\ker(B) = \ker(B^\ast) = \{0\}.$ 
\end{lemma}
\begin{proof}
  It is well-known that $H^1(\R)\subset L^2(\R)$ is a dense Hilbert space and $B$ is closed, since its graph norm equals the $H^1(\R)$ norm. Moreover, $Bv=0$ implies that $v\equiv c\in\R$ and since the only constant belonging to $L^2(\R)$ is zero, we have $\ker(B)=\{0\}$. With the same reasoning the adjoint operator $B^\ast : H^1(\R)\subset L^2(\R)\to L^2(\R)$ with $B^\ast v = -\partial_t v$, satisfies $\ker(B^\ast) = \{0\}.$
\end{proof}

Due to the properties of the derivative, we can apply Theorem \ref{thm:polar_decomposition}, to get the following result. 
\begin{theorem}
  Let $B := \partial_t : H^1(\R) \subset L^2(\R) \to L^2(\R)$. The natural energy space in this setting is 
  \begin{equation*}
    X = \dom(|\partial_t|^{1/2}) = [L^2(\R),H^1(\R)]_{1/2} = H^{1/2}(\R). 
  \end{equation*}
  Moreover, 
  \begin{equation*}
    B = \Hcl (-\partial_{tt})^{1/2}, 
  \end{equation*}
  where $\Hcl:L^2(\R)\to L^2(\R)$ denotes the classical Hilbert transformation, defined as 
  \begin{equation*}
    (\Hcl v)(t) := \frac{1}{\pi} \mathrm{p.v.} \int_\R \frac{v(s)}{s-t}\, ds. 
  \end{equation*}
\end{theorem}
\begin{remark}
  Note that our sign convention for the classical Hilbert transform is
  $\Hcl=-\mathcal{H}_{\mathrm{std}}$, where the standard Hilbert transform is given by
  \begin{equation*}
    (\mathcal{H}_{\mathrm{std}}v)(t)
    = \frac{1}{\pi}\,\mathrm{p.v.}\!\int_{\R} \frac{v(s)}{t-s}\,ds.
  \end{equation*}
\end{remark}
\begin{proof}
  The defintion of the energy space is an immediate consequence of Lemma \ref{lem:unbounded derivative} and Theorem \ref{thm:polar_decomposition}. To compute the unitary factor, note that 
  $B^\ast B = -\partial_{tt}:\dom(B^\ast B) \subset L^2(\R) \to L^2(\R),$ with domain $\dom(B^\ast B) = H^2(\R)$. Denoting by $\mathcal{F}v$ the Fourier transform of a function, we have  
  \begin{equation*}
    (\mathcal{F}B v)(\omega) 
    = 
    \frac{1}{\sqrt{2\pi}} \int_\R e^{-i\omega t} \partial_t v(t)\, dt
    =
    i\omega (\mathcal{F}v)(\omega), 
   \end{equation*} 
   and 
   \begin{equation*}
      (\mathcal{F}B^\ast B v)(\omega)
      =
      \frac{1}{\sqrt{2\pi}} \int_\R e^{-i\omega t} (-\partial_{tt}v(t))\, dt 
      = 
      \omega^2 (\mathcal{F}v)(\omega). 
   \end{equation*}
   In particular, the Fourier transform diagonalizes $B^\ast B$, see e.g. \cite[Section 4.3]{EvansPartialDifferentialEquations2010}, as in operator notation we can write  
   \begin{equation*}
      (\mathcal{F}B^\ast B \mathcal{F}^{-1}\widehat{v})(\omega) = \omega^2 \widehat{v}(\omega). 
   \end{equation*}
   Moreover, by Plancherel's theorem the Fourier transform is a unitary operator, and using functional calculus we have that 
   \begin{equation*}
      |\omega| \widehat{v}(\omega) 
      =
      ((\mathcal{F} B^\ast B \mathcal{F}^{-1})^{1/2}\widehat{v})(\omega) 
      =
      (\mathcal{F} (B^\ast B)^{1/2} \mathcal{F}^{-1}\widehat{v})(\omega).
   \end{equation*}
   Thus, we have that  
   \begin{equation*}
     (\mathcal{F} B \mathcal{F}^{-1} \widehat{v})(\omega)
     =
     i \omega \widehat{v}(\omega)
     =
     i\operatorname{sign}(\omega)|\omega| \widehat{v}(\omega),
   \end{equation*} 
   where $\operatorname{sign}(\cdot)$ denotes the sign function. 
   Equivalently 
   \begin{equation*}
      (Bv)(t) 
      =
      (\mathcal{F}^{-1} i\operatorname{sign}(\omega) |\omega| \mathcal{F}v)(t)
      =
      (\mathcal{F}^{-1} i\operatorname{sign}(\omega) \mathcal{F})(\mathcal{F}^{-1} |\omega| \mathcal{F}v)(t).
   \end{equation*}
   Now, using the relations 
   \begin{equation*}
      (\mathcal{F}^{-1} |\omega| \mathcal{F}v)(t) 
      = (|B|v)(t) 
      = ((B^\ast B)^{1/2} v)(t) = (-\partial_{tt})^{1/2}v(t),
   \end{equation*}
   and 
   \cite[Equation (1.14.46)]{NIST:DLMF} 
   \begin{equation*}
      (\mathcal{F}^{-1} i\operatorname{sign}(\omega)\mathcal{F} v)(t) = (\Hcl v)(t), 
   \end{equation*}
   we obtain the desired result. 
\end{proof}

\subsection{The derivative on a bounded interval}
For some fixed $T>0$, we consider the derivative as an unbounded operator
\begin{align}
  B\, : \, \dom(B)\subset L^2(0,T) &\to L^2(0,T), 
  \quad  B v = \partial_t v, \label{eq:primal_time_derivative}
\end{align}
with domain
\begin{equation*}
  \dom(B) = H^1_{0,}(0,T) := \{ v \in H^1(0,T) \mid v(0)=0 \}.
\end{equation*}
Since $H^1(0,T)\subset C([0,T])$, the domain is well-defined. 
\begin{lemma}
  The operator $B$ is densely defined, closed and $\ker(B) = \ker(B^\ast) = \{0\}.$ Moreover, it satisfies the inf-sup condition \eqref{eq:inf-sup B}, i.e.,  
  \begin{equation*}
    \frac{\pi}{2T} \leq \inf_{0\neq u\in H^1_{0,}(0,T)} \sup_{0\neq v\in L^2(0,T)} \frac{(Bu,v)_{L^2(0,T)}}{\norm{u}_{L^2(0,T)}\norm{v}_{L^2(0,T)}}. 
  \end{equation*}
\end{lemma}
\begin{proof}
  It is well-known that $H^1_{0,}(0,T)\subset L^2(0,T)$ is dense.  
  Due to the Poincar\'{e} inequality, 
  \begin{equation}\label{eq:Poincare 1D}
    \norm{v}_{L^2(0,T)} 
    \leq 
    \frac{2T}{\pi}\norm{\partial_t v}_{L^2(0,T)},\quad \forall v\in H^1_{0,}(0,T),
  \end{equation}
  we immediately have that $B$ is closed and $\ker(B) = \{0\}.$
  Its adjoint admits the domain $\dom(B^\ast) = H^1_{,0}(0,T):= \{v\in H^1(0,T):\, v(T)=0\}$ and the action $B^\ast v = -\partial_t v.$ Since \eqref{eq:Poincare 1D} also holds true for functions in $H^1_{,0}(0,T)$, we also have $\ker(B^\ast) = \{0\}$. 
  Now for any arbitrary $0\neq u\in H^1_{0,}(0,T)$, choosing $\hat v = \partial_t u\in L^2(0,T)$ together with \eqref{eq:Poincare 1D} gives 
  \begin{equation*}
    \sup_{0\neq v\in L^2(0,T)} \frac{(Bu,v)_{L^2(0,T)}}{\norm{v}_{L^2(0,T)}}
    \geq 
    \frac{\int_0^T \partial_t u(t) \hat v(t)\, dt }{\norm{\hat v}_{L^2(0,T)}} 
    =
    \norm{\partial_t u}_{L^2(0,T)}
    \geq 
    \frac{\pi}{2T}\norm{u}_{L^2(0,T)}.
  \end{equation*}   
\end{proof}

All the assumptions of Theorem \ref{thm:polar_decomposition} are met. Before we give the main result, let us recall the definition of the \emph{modified Hilbert transformation}, as given in \cite{ZankSteinbachCoerciveSpaceTimeFinite2020}. Therefore, we consider the eigenvalue problems 
\begin{equation*}
    -\partial_{tt} s_k(t) = \lambda_k s_k(t),\, t\in (0,T), \quad s_k(0) = \partial_t s_k(T) =0,  
  \end{equation*} 
  and 
  \begin{equation*}
    -\partial_{tt} c_k(t) = \lambda_k c_k(t),\, t\in (0,T), \quad c_k(T) = \partial_t c_k(0) =0,  
  \end{equation*} 
  for $k\in \N_0.$
  Their solutions 
  \begin{equation*}
    s_k(t) = \sin\left(\left(\frac{\pi}{2}+k\pi\right)\frac{t}{T}\right),\quad 
    c_k(t) = \cos\left(\left(\frac{\pi}{2}+k\pi\right)\frac{t}{T}\right),
  \end{equation*}
  built an orthogonal basis of $H^s_{0,}(0,T)$ and $H^s_{,0}(0,T)$, for $s\in[0,1]$, respectively, satisfying for all $k,\ell\in \mathbb{N}_0,$
  \begin{equation*}
    \int_0^T s_k(t)s_\ell(t)\, dt = \int_0^T c_k(t)c_\ell(t)\, dt = \frac{T}{2}\delta_{k,\ell} 
   \quad \delta_{k,\ell}:= \begin{cases}
      1, & k=\ell,\\
      0, & k\neq \ell,
    \end{cases}
  \end{equation*}
  and their eigenvalues are given by 
  \begin{equation*}
    \lambda_k = \left(\frac{\pi}{2}+k\pi\right)^2\frac{1}{T^2}\, ,\quad k\in \mathbb{N}_0. 
  \end{equation*}
  Hence, for any function $v\in H^s_{0,}(0,T)$ there holds the representation 
  \begin{equation*}
    v(t) 
    =
    \sum_{k=0}^\infty v_k s_k(t),\quad 
    v_k = \frac{2}{T} \int_0^T v(t)s_k(t)\, dt. 
  \end{equation*}
  \begin{definition}
    The modified Hilbert transformation $\Hmod:H^s_{0,}(0,T) \to H^s_{,0}(0,T)$ is defined as 
    \begin{equation*}
      (\Hmod v)(t)
      =
      \Hmod\left(\sum_{k=0}^\infty v_k s_k(t)\right)
      := \sum_{k=0}^\infty v_k c_k(t). 
    \end{equation*} 
  \end{definition}
  We stress that the modified Hilbert transformation is unitary \cite{LoescherSteinbachZank:2024}, i.e, 
  \begin{equation*}
    \norm{\Hmod v}_{H^s_{,0}(0,T)} = \norm{v}_{H^s_{0,}(0,T)},\quad \forall v\in H^s_{0,}(0,T),\, s\in [0,1].
  \end{equation*}
  \begin{remark}
    Using that $\cos(x) = (-1)^k \sin(k\pi +\tfrac{\pi}{2} - x)$ we compute 
    \begin{align*}
      c_k(t) 
      = \cos\left(\left(k\pi +\frac{\pi}{2}\right)\frac{t}{T}\right)
      &=
      (-1)^k\sin\left(k\pi +\frac{\pi}{2} -\left(k\pi + \frac{\pi}{2}\right)\frac{t}{T}\right)\\
      &=
      (-1)^k\sin\left(\left(k\pi + \frac{\pi}{2}\right)\frac{T-t}{T}\right)\\
      &=
      (-1)^k s_k(T-t).  
    \end{align*}
    This establishes a relation between the time reversal operator and the modified Hilbert transforamtion, i.e., 
    \begin{equation*}
      (\Hmod s_k)(t) = c_k(t) = (-1)^k s_k(T-t). 
    \end{equation*} 
  \end{remark}
\begin{theorem}
  Let $B:= \partial_t:H^1_{0,}(0,T)\subset L^2(0,T)\to L^2(0,T)$. The natural energy space in this case is given as 
  \begin{equation*}
    X = \dom(|B|^{1/2}) = [L^2(0,T),H^1_{0,}(0,T)]_{1/2} = H^{1/2}_{0,}(0,T). 
  \end{equation*}
  Moreover, 
  \begin{equation*}
    B = \Hmod (-\partial_{tt})^{1/2},
  \end{equation*}
  where $\Hmod:L^2(0,T)\to L^2(0,T)$ is the modified Hilbert transformation. 
\end{theorem}
\begin{proof}
  The energy space is immediate from Theorem \ref{thm:polar_decomposition}. 
  It remains to compare the action of $B$ and $|B|$ on the basis functions. The application of $B$ gives 
  \begin{equation*}
    (B s_k)(t) = \partial_t s_k(t) = \sqrt{\lambda_k} c_k(t),
  \end{equation*}
  while 
  \begin{equation*}
    (|B|s_k)(t) = (-\partial_{tt})^{1/2} s_k(t) = \sqrt{\lambda_k} s_k(t). 
  \end{equation*}
  Thus, the unitary factor has to map 
  \begin{equation*}
    (\Hmod s_k)(t) = c_k(t),
  \end{equation*}
  which is exactly the definition of the modified Hilbert transformation. 
\end{proof}

\section{The Ordinary Differential Equation}
In this section we consider the ordinary differential equation (ODE)
\[
  \partial_t u(t)+\mu u(t)=f(t)
  \quad\text{for }t\in(0,T),
  \qquad u(0)=0.
\]
This equation also governs the temporal coefficients in an eigenfunction expansion of a parabolic evolution problem, with \(\mu\geq 0\) representing an eigenvalue of the spatial operator. In particular, stability estimates will transfer to the parabolic case, studied later on.

\subsection{The Polar Decomposition of the ordinary differential equation}
In the previous section, we have analyzed the first-derivative operator. 
For the treatment of the heat equation later on, including a zeroth order term will be required. 
For some fixed $T>0$ and $\mu\geq 0$ we consider 
\begin{align}
  B_\mu \, : \, \dom(B_\mu) \subset L^2(0,T) \to L^2(0,T),
  \quad 
  B_\mu v = \partial_t v + \mu v, \label{eq:primal_ode_operator}
\end{align}
with domain
\begin{align}
  \dom(B_\mu) &= H^1_{0,}(0,T).
\end{align}
\begin{lemma}
  The operator $B_\mu$ is densely defined, closed and $\ker(B_\mu) = \ker(B_\mu^\ast) = \{0\}.$ Moreover, it satisfies the inf-sup condition \eqref{eq:inf-sup B}, i.e.,  
  \begin{equation*}
    \frac{\pi}{2T} \leq \inf_{0\neq u\in H^1_{0,}(0,T)} \sup_{0\neq v\in L^2(0,T)} \frac{(B_\mu u,v)_{L^2(0,T)}}{\norm{u}_{L^2(0,T)}\norm{v}_{L^2(0,T)}}. 
  \end{equation*}
\end{lemma}
\begin{proof}
  We already mentioned that $H^1_{0,}(0,T)\subset L^2(0,T)$ is dense. 
  Noting that 
  \begin{align*}
    \norm{B_\mu v}_{L^2(0,T)}^2 
    &=
    \norm{\partial_t v+\mu v}_{L^2(0,T)}^2\\
    &= 
    \norm{\partial_t v}_{L^2(0,T)}^2  + 2\mu \int_0^T \partial_t v(t)v(t)\, dt  +\mu^2\norm{v}_{L^2(0,T)}^2 \\
    & =
    \norm{\partial_t v}_{L^2(0,T)}^2 + \mu |v(T)|^2 + \mu^2 \norm{v}_{L^2(0,T)}^2\\
    & \geq 
    \norm{\partial_t v}_{L^2(0,T)}^2,
  \end{align*}
  and using the Poincar\'{e} inequality \eqref{eq:Poincare 1D} gives 
  \begin{equation*}
    \norm{v}_{L^2(0,T)}
    \leq 
    \frac{2T}{\pi} \norm{B_\mu v}_{L^2(0,T)},
  \end{equation*}
  which shows that $B_\mu$ is closed and $\ker(B_\mu) = \{0\}$. Analogously, we can derive $\ker(B_\mu^\ast) = \{0\}$. 
  Moreover, since for any $u\in \dom(B_\mu)$ 
  \begin{equation*}
    \norm{B_\mu u}_{L^2(0,T)}
    =
    \sup_{0\neq v\in L^2(0,T)} \frac{(B_\mu u,v)_{L^2(0,T)}}{\norm{v}_{L^2(0,T)}},
  \end{equation*}   
  this also shows the inf-sup condition and finishes the proof. 
\end{proof}

\begin{figure}[ht]
      \centering
      \begin{tikzpicture}
        \begin{axis}[%
            axis lines=middle,
            axis on top,
            legend cell align=left,
            xlabel=$\omega$,
            domain=0:2*pi,
            xmin=0,
            xmax=7,
            ymin=-5,
            ymax=5,
            trig format plots=rad, 
            xtick={-2*pi,-3*pi/2, -pi, -pi/2,pi/2,pi,3*pi/2,2*pi},
            xticklabels={$-2\pi$, $-\frac{3\pi}{2}$, $-\pi$, $-\frac{\pi}{2}$, $\frac{\pi}{2}$,$\pi$,$\frac{3\pi}{2}$,$2\pi$},
            every axis y label/.style={rotate=0, black, at={(0.5,1.05)},},
            every axis x label/.style={rotate=0, black, at={(1.05,0.5)},},,
            font=\footnotesize,     
        ]
        \pgfplotsinvokeforeach{-5,-3,...,3}{
        \pgfmathsetmacro{\xmin}{ifthenelse(#1==-5,-2*pi,#1*pi/2+0.01)}
        \pgfmathsetmacro{\xmax}{ifthenelse(#1==3,2*pi,#1*pi/2+pi-0.01)}
        \addplot[samples=51,smooth,domain=\xmin:\xmax,forget plot]{tan(x)};
        \draw[densely dotted] (#1*pi/2,\pgfkeysvalueof{/pgfplots/ymin})
        -- (#1*pi/2,\pgfkeysvalueof{/pgfplots/ymax});
        }
        \addplot[domain=-1:8, samples=2, dashdotdotted, blue]{-x/1};
        \addlegendentry{$\mu = 1$}
        \addplot[domain=-1:8, samples=2, densely dashdotdotted, red]{-x/20};
        \addlegendentry{$\mu = 20$}
        \end{axis}
      \end{tikzpicture}
      \caption{Graphical solution of the equation $\omega \cos(\omega T) + \mu \sin(\omega T) = 0$ for $T=1$ and two different values of $\mu$. The intersections of the blue and red line with the graph of the tangent function give the solutions $\omega_{\mu,k}$ for the respective value of $\mu$.} \label{fig:graphical_solution_tan}
    \end{figure}
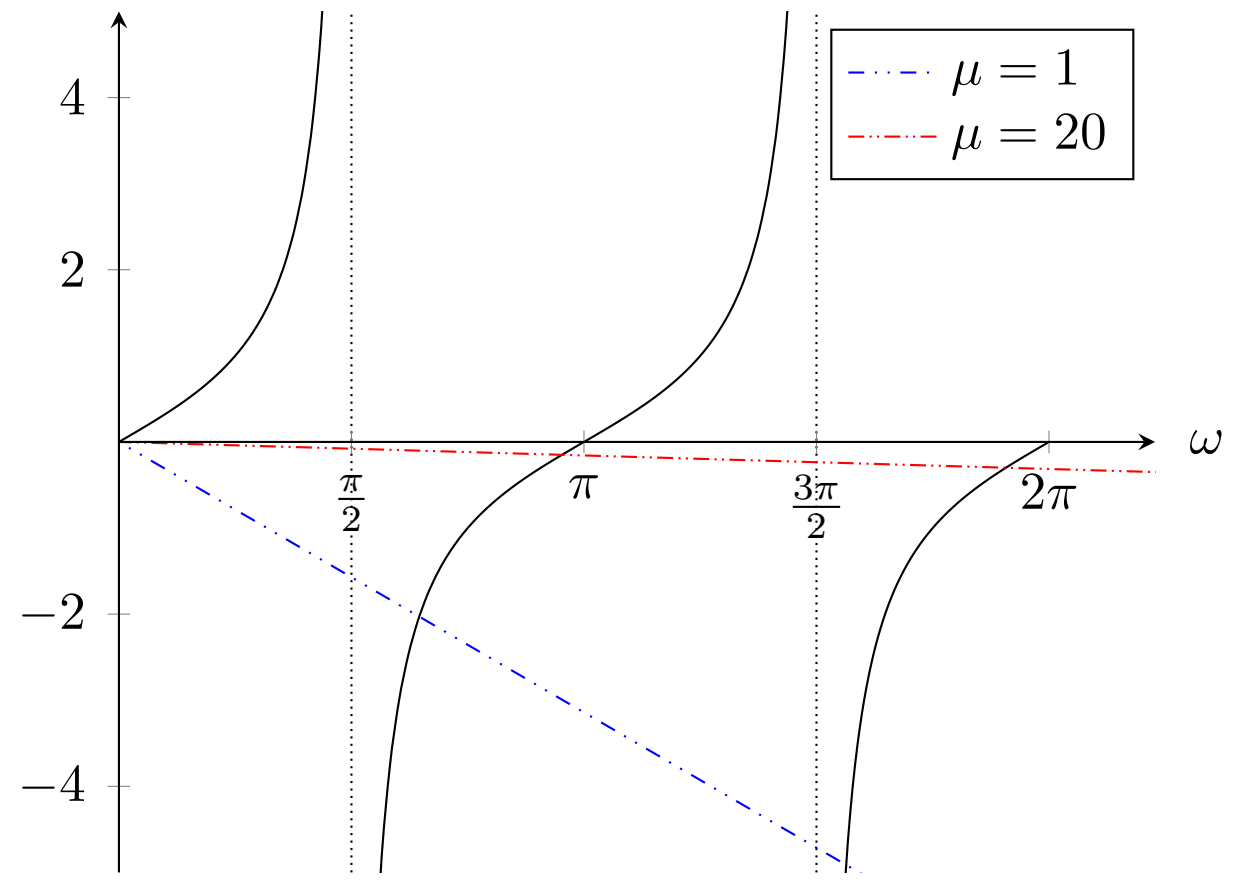

Hence, all the assumptions of Theorem \ref{thm:polar_decomposition} are met, from which we we can immediately conclude the energy space. Before we state the main result, and to gain a better understanding of the unitary operator in the polar decomposition, we will consider suitable bases of eigenfunctions. Therefore, note that  
\begin{align*}
  B_\mu^\ast B_\mu:\dom(B_\mu^\ast B_\mu)\subset L^2(0,T)\to L^2(0,T),
  \quad B_\mu^\ast B_\mu v = -\partial_{tt} v + \mu^2 v
\end{align*}
with domain $\dom(B_\mu^\ast B_\mu) = \{v\in H^2(0,T): v(0) = \partial_t v(T) + \mu v(T)=0\}$ is self-adjoint, positive and has a compactly embedded form domain and thus its spectrum is discrete, and the eigenfunctions build an orthogonal basis in $L^2(0,T)$, see e.g. \cite{ArendtUrbanPartialDifferentialEquations2023}. Indeed, one can easily check that the eigenfunctions are the solutions of 
\begin{equation}\label{eq:EVP for s_mu}
\begin{aligned}
  -\partial_{tt} s_{\mu,k}(t) + \mu^2 s_{\mu,k}(t) &= \lambda_{\mu,k} s_{\mu,k}(t)\quad \text{for }
  t\in (0,T),\\ 
  s_{\mu,k}(0) &= \partial_t s_{\mu,k}(T) + \mu s_{\mu,k}(T)=0,
\end{aligned}   
\end{equation}
which are given as 
\begin{equation*}
  s_{\mu,k}(t) = \sin(\omega_{\mu,k} t),\quad \lambda_{\mu,k} = \mu^2 + \omega_{\mu,k}^2, \quad k\in \mathbb{N}_0,   
\end{equation*}
where $\omega_{\mu,k}$ is the unique root of, see Figure \ref{fig:graphical_solution_tan}, 
\begin{equation}\label{eq:omega as root}
  \omega_{\mu,k} \cos(\omega_{\mu,k} T)  + \mu \sin(\omega_{\mu,k} T) =0, \quad 
  \omega_{\mu,k} \in \left[\left(k\pi +\frac{\pi}{2}\right) \frac{1}{T},\frac{(k+1)\pi}{T}\right).  
\end{equation}
The same holds true for the operator   
\begin{align*}
  B_\mu B_\mu^\ast:\dom(B_\mu B_\mu^\ast)\subset L^2(0,T)\to L^2(0,T),
  \quad B_\mu B_\mu^\ast v = -\partial_{tt}v + \mu^2 v
\end{align*}
with domain $\dom(B_\mu B_\mu^\ast) = \{v\in H^2(0,T): v(T) = -\partial_t v(0) + \mu v(0) =0\}$. 
The eigenfunctions in this case are the solutions of 
\begin{align*}
  -\partial_{tt} c_{\mu,k}(t) + \mu^2 c_{\mu,k}(t) &= \lambda_{\mu,k} c_{\mu,k}(t)\quad \text{for }
  t\in (0,T),\\ 
  c_{\mu,k}(T) &= -\partial_t c_{\mu,k}(0) + \mu c_{\mu,k}(0)=0,
\end{align*} 
which are given as 
\begin{equation*}
  c_{\mu,k}(t) = \frac{\omega_{\mu,k}\cos(\omega_{\mu,k} t) + \mu\sin(\omega_{\mu,k}t)}{\sqrt{\omega_{\mu,k}^2 + \mu^2}},\quad \lambda_{\mu,k} = \mu^2 + \omega_{\mu,k}^2.  
\end{equation*}

\begin{remark}
  Using trigonometric relations and \eqref{eq:omega as root}, we compute 
  \begin{align*}
    s_{\mu,k}(T-t) 
    &= \sin(\omega_{\mu,k} (T-t))
    = \sin(\omega_{\mu,k}T)\cos(\omega_{\mu,k}t) - \cos(\omega_{\mu,k}T) \sin(\omega_{\mu,k}t)\\
    &= \sin(\omega_{\mu,k}T)\left(\cos(\omega_{\mu,k}t) + \frac{\mu}{\omega_{\mu,k}}\sin(\omega_{\mu,k}t)\right). 
  \end{align*}
  Moreover, using that $\sin(\omega_{\mu,k}T)^2+\cos(\omega_{\mu,k}T)^2=1$ and \eqref{eq:omega as root}, we compute that 
  \begin{equation*}
    \sin(\omega_{\mu,k}T) = (-1)^k \frac{\omega_{\mu,k}}{\sqrt{\omega_{\mu,k}^2 + \mu^2}},
  \end{equation*}
  to conclude the relation 
  \begin{equation}\label{eq:time reversal vs space time hilbert}
    c_{\mu,k}(t) = (-1)^k s_{\mu,k}(T-t). 
  \end{equation}
\end{remark}

For any $k,\ell\in\mathbb{N}_0$, we can compute the orthogonality  
\begin{equation}\label{eq:orthogonality in B_mu basis}
  \int_0^T s_{\mu,k}(t)s_{\mu,\ell}(t)\, dt 
  = 
  \int_0^T c_{\mu,k}(t)c_{\mu,\ell}(t)\, dt 
  =
  \left(\frac{T}{2} + \frac{\mu}{2\lambda_{\mu,k}}\right)\delta_{k,\ell},
\end{equation}
and for any $v\in L^2(0,T)$ there hold the representations 
\begin{equation}\label{eq:representation v in B_mu s basis}
  v(t) = \sum_{k=0}^\infty v_{\mu,k} s_{\mu,k}(t),\quad 
  v_{\mu,k} = \left(\frac{T}{2} + \frac{\mu}{2\lambda_{\mu,k}}\right)^{-1}\int_0^T v(t)s_{\mu,k}(t)\, dt, 
\end{equation}
and 
\begin{equation}\label{eq:representation v in B_mu c basis}
  v(t) = \sum_{k=0}^\infty \overline{v}_{\mu,k} c_{\mu,k}(t),\quad 
  \overline{v}_{\mu,k} = \left(\frac{T}{2} + \frac{\mu}{2\lambda_{\mu,k}}\right)^{-1}\int_0^T v(t)c_{\mu,k}(t)\, dt. 
\end{equation}

As in the case for the first-derivative, this motivates the following definition. 
\begin{definition}\label{def:ParabolicHilbertTransform}
  We define the \emph{parabolic Hilbert transformation} $\Hmu:L^2(0,T)\to L^2(0,T)$ as 
  \begin{equation*}
    (\Hmu v)(t) 
    = \Hmu \left(\sum_{k=0}^\infty  v_{\mu,k} s_{\mu,k}(t)\right) 
    = \sum_{k=0}^\infty v_{\mu,k} c_{\mu,k}(t). 
  \end{equation*}
\end{definition}

\begin{lemma}\label{lem:parabolic-hilbert-properties}
  Let $T>0$ and $\mu\geq0$. The parabolic Hilbert transformation
  $\mathcal H_\mu$ has the following properties.
  \begin{enumerate}
    \item[(i)]
      For every $s\in[0,1]$, the operator
      $\mathcal H_\mu:H^s_{0,}(0,T)\longrightarrow H^s_{,0}(0,T)$
      is unitary when the spaces are equipped with the spectral
      norms
      \begin{equation*}
        \norm{v}_{H^s_{0,}(0,t),\mu}^2 := \sum_{k=0}^\infty \left(\frac{T}{2} + \frac{\mu}{2\lambda_{\mu,k}}\right)\lambda_{\mu,k}^s v_{\mu,k}^2, 
      \end{equation*}
      and 
      \begin{equation*}
        \norm{w}_{H^s_{,0}(0,t),\mu}^2 := \sum_{k=0}^\infty \left(\frac{T}{2} + \frac{\mu}{2\lambda_{\mu,k}}\right)\lambda_{\mu,k}^s \overline{w}_{\mu,k}^2,
      \end{equation*}
      respectively. For fixed $\mu\geq 0$, the norm equivalences hold
      \begin{align*}
        &2^{s-1}\big(\norm{v}_{H^s_{0,}(0,T)}^2 + \mu^{2s}\norm{v}_{L^2(0,T)}^2\big)\\
        &\hskip 5mm \leq \norm{v}_{H^{s}_{0,}(0,T),\mu}^2 \\
        &\hskip 10mm \leq 2\big(\norm{v}_{H^s_{0,}(0,T)}^2 + \mu^{2s} \norm{v}_{L^2(0,T)}^2\big),\\
        &2^{s-1}\big(\norm{w}_{H^s_{,0}(0,T)}^2 + \mu^{2s}\norm{w}_{L^2(0,T)}^2\big)\\
        &\hskip 5mm \leq \norm{w}_{H^{s}_{,0}(0,T),\mu}^2 \\
        &\hskip 10mm \leq 2\big(\norm{w}_{H^s_{,0}(0,T)}^2 + \mu^{2s} \norm{w}_{L^2(0,T)}^2\big).
      \end{align*} 
    \item[(ii)]
      It holds that 
      \[
        \|\mathcal H_\mu-\mathcal H_T\|
        _{L^2(0,T)\to L^2(0,T)}\leq \frac{2}{\sqrt{3}} T\mu.
      \]
      In particular, $\mathcal{H}_{0} = \Hmod$. 

    \item[(iii)]
      As $\mu\to\infty$, the operators $\mathcal H_\mu$
      converge strongly to the identity on $L^2(0,T)$, i.e.,
      \[
        \lim_{\mu\to\infty}
        \|\mathcal H_\mu v-v\|_{L^2(0,T)}=0,
        \qquad\forall v\in L^2(0,T).
      \]
  \end{enumerate}
\end{lemma}

\begin{proof}
  \emph{(i)}
  Since $\mathcal H_\mu$ preserves the expansion coefficients, we have 
  \[
    \|\mathcal H_\mu v\|_{H^{s}_{,0}(0,T),\mu}^2
    =
    \sum_{k=0}^{\infty}
    \left(\frac{T}{2} + \frac{\mu}{2\lambda_{\mu,k}}\right)\lambda_{\mu,k}^{s}v_{\mu,k}^2
    =
    \norm{v}_{H^{s}_{0,}(0,T),\mu}^2.
  \]
  
  To show the norm equivalence, firstly, note that by the expansion \eqref{eq:representation v in B_mu s basis} and the orthogonality \eqref{eq:orthogonality in B_mu basis} we have that  
  \begin{align*}
    \norm{v}_{L^2(0,T)}^2
    & = 
    \int_0^T \left(\sum_{k=0}^\infty v_{\mu,k} s_{\mu,k}(t)\right)^2 \, dt
    = 
    \sum_{k,\ell=0}^\infty v_{\mu,k}v_{\mu,\ell}\int_0^T s_{\mu,k}(t)s_{\mu,\ell}(t)\, dt \\
    & =
    \sum_{k=0}^\infty \left(\frac{T}{2} + \frac{\mu}{2\lambda_{\mu,k}}\right) v_{\mu,k}^2
    = \norm{v}_{H^0_{0,},\mu}^2. 
  \end{align*}
  Secondly, for $s=1$, we have that 
  \begin{align*}
    \norm{v}_{H^1_{0,}(0,T),\mu}^2 
    &=
    \sum_{k=0}^\infty \left(\frac{T}{2} + \frac{\mu}{2\lambda_{\mu,k}}\right)\lambda_{\mu,k}v_{\mu,k}^2\\
    & = 
    \sum_{k,\ell=0}^\infty v_{\mu,k}v_{\mu,\ell}\int_0^T \lambda_{\mu,k}s_{\mu,k}(t)s_{\mu,\ell}(t)\, dt.
  \end{align*}
  Using that $s_{\mu,k}$ is an eigenfunction of \eqref{eq:EVP for s_mu} and integrating by parts then gives 
  \begin{align*}
    &\sum_{k,\ell=0}^\infty v_{\mu,k}v_{\mu,\ell}\int_0^T \big(-\partial_{tt}s_{\mu,k}(t) + \mu^2 s_{\mu,k}(t)\big)s_{\mu,\ell}(t)\, dt \\
    & = 
    \sum_{k,\ell=0}^\infty v_{\mu,k}v_{\mu,\ell}\left[ \int_0^T \partial_t s_{\mu,k}(t)\partial_ts_{\mu,\ell}(t)\, dt 
    + \mu^2 \int_0^T s_{\mu,k}(t)s_{\mu,\ell}(t)\, dt 
    +\mu s_{\mu,k}(T)s_{\mu,\ell}(T) \right]\\  
    & = 
    \norm{\partial_t v}_{L^2(0,T)}^2 + \mu^2 \norm{v}_{L^2(0,T)}^2 + \mu v(T)^2.  
  \end{align*}
  Moreover, since 
  \begin{align*}
    0&\leq \mu|v(T)|^2 
    = \int_0^T \mu\frac{d}{dt} |v(t)|^2 \, dt 
    = 2\int_0^T \partial_t v(t)\mu v(t)\, dt \\
    &\leq \norm{\partial_t v}_{L^2(0,T)}^2 +\mu^2 \norm{v}_{L^2(0,T)}^2,
  \end{align*}
  there holds the norm equivalence 
  \begin{equation*}
    \norm{\partial_t v}_{L^2(0,T)}^2 + \mu^2\norm{v}_{L^2(0,T)}^2 
    \leq 
    \norm{v}_{H^1_{0,}(0,T),\mu}^2
    \leq 
    2(\norm{\partial_t v}_{L^2(0,T)}^2 +\mu^2 \norm{v}_{L^2(0,T)}^2),
  \end{equation*}
  for all $v\in H^1_{0,}(0,T).$
  Analogously, using the expansion \eqref{eq:representation v in B_mu c basis}, we get 
  \begin{equation*}
    \norm{w}_{H^1_{,0}(0,T),\mu}^2 
    =
    \norm{\partial_t w}_{L^2(0,T)}^2 + \mu^2 \norm{w}_{L^2(0,T)}^2 + \mu |w(0)|^2. 
  \end{equation*}
  and the norm equivalence
  \begin{equation*}
    \norm{\partial_t w}_{L^2(0,T)}^2 + \mu^2\norm{w}_{L^2(0,T)}^2 
    \leq 
    \norm{w}_{H^1_{,0}(0,T),\mu}^2
    \leq 
    2(\norm{\partial_t w}_{L^2(0,T)}^2 +\mu^2 \norm{w}_{L^2(0,T)}^2),
  \end{equation*}
  for all $w\in H^1_{,0}(0,T)$. Now, let us consider the representation
  \begin{equation*}
    v(t) = \sum_{k=0}^\infty v_k s_k(t),\quad v_k = \frac{2}{T}\int_0^T v(t)s_k(t)\, dt. 
  \end{equation*}
  Then 
  \begin{equation*}
    \norm{v}_{L^2(0,T)}^2 
    =
    \frac{T}{2}\sum_{k=0}^\infty v_k^2, 
  \end{equation*}
  and
  \begin{equation*}
    \norm{\partial_t v}_{L^2(0,T)}^2 +\mu^2 \norm{v}_{L^2(0,T)}^2 
    =
    \frac{T}{2}\sum_{k=0}^\infty (\lambda_k +\mu^2) v_k^2. 
  \end{equation*}
  Interpolating the norms, using the established norm equivalences, gives   
  \begin{equation*}
    \frac{T}{2} \sum_{k=0}^\infty (\lambda_k+\mu^2)^s v_k^2 
    \leq 
    \norm{v}_{H^s_{0,}(0,T),\mu}^2 
    \leq  T\sum_{k=0}^\infty (\lambda_k +\mu^2)^sv_k^2.
  \end{equation*}
  Further, using that $2^{s-1}(\lambda_k^s+\mu^{2s})\leq (\lambda_k +\mu^2)^s \leq \lambda_k^s + \mu^{2s}$
  as well as 
  \begin{equation*}
    \norm{v}_{H^{s}_{0,}(0,T)}^2 = \frac{T}{2}\sum_{k=0}^\infty \lambda_k^s v_k^2, 
  \end{equation*}
  we obtain 
  \begin{align*}
    &2^{s-1} \left(\norm{v}_{H^s_{0,}(0,T)}^2 + \mu^{2s}\norm{v}_{L^2(0,T)}^2\right)\\
    &\hskip 5 mm \leq 
    \norm{v}_{H^s_{0,}(0,T),\mu}^2 \\
    &\hskip 10mm \leq
    2\left(\norm{v}_{H^s_{0,}(0,T)}^2 + \mu^{2s}\norm{v}_{L^2(0,T)}^2\right). 
  \end{align*}
  The norm equivalence in $H_{,0}^1(0,T)$ can be derived analogously, using the expansion 
  \begin{equation*}
    w(t) = \sum_{k=0}^\infty \overline{w}_k c_k(t),\quad 
    \overline{w}_k = \frac{2}{T} \int_0^T w(t)c_k(t)\, dt. 
  \end{equation*}
   
  \emph{(ii).} 
  For any fixed $\mu\geq 0$ let us consider 
  \begin{equation*}
    0\leq \delta_{\mu,k} := T(\omega_{\mu,k} - \omega_{0,k}) 
    = T\omega_{\mu,k} - \left(\frac{\pi}{2} + k\pi\right). 
  \end{equation*}
  Using trigonometric identities it holds that  
  \begin{equation*}
    \sin(\omega_{\mu,k}T)
    =
    \sin\left(\left(\frac{\pi}{2} + k\pi\right) + \delta_{\mu,k}\right)
    =
    (-1)^k \cos(\delta_{\mu,k}),
  \end{equation*}
  as well as 
  \begin{equation*}
    \cos(\omega_{\mu,k}T) 
    =
    \cos\left(\left(\frac{\pi}{2} + k\pi\right) + \delta_{\mu,k}\right)
    =
    (-1)^{k+1}\sin(\delta_{\mu,k}),
  \end{equation*}
  and thus the equation \eqref{eq:omega as root} for $\omega_{\mu,k}$ becomes 
  \begin{equation*}
    -\omega_{\mu,k}\sin(\delta_{\mu,k}) + \mu \cos(\delta_{\mu,k}) =0. 
  \end{equation*}
  Reordering gives 
  \[
    \delta_{\mu,k} = \arctan\left(\frac{\mu}{\omega_{\mu,k}}\right) 
    \leq
    \frac{\mu}{\omega_{\mu,k}}
    \leq 
    \frac{\mu}{\omega_{0,k}}.  
  \]
  Consequently,
  \[
    0\leq\omega_{\mu,k}-\omega_{0,k}
    \leq\frac{\mu}{T\omega_{0,k}}.
  \]
  and using that 
  \[
    |\sin(\omega_{\mu,k}t)-\sin(\omega_{0,k} t)|
    \leq|\omega_{\mu,k}-\omega_{0,k}|\,t = \delta_{\mu,k} \frac{t}{T}
  \]
  we obtain
  \begin{align*}
    \|s_{\mu,k}-s_{0,k}\|_{L^2(0,T)}
    & \leq \sqrt{\frac{\delta_{\mu,k}^2 T}{3}} 
    \leq \sqrt{\frac{T}{3}}\frac{\mu}{\omega_{0,k}}.
  \end{align*}
  Due to \eqref{eq:time reversal vs space time hilbert} we also have 
  \begin{equation*}
    \norm{c_{\mu,k} - c_{0,k}}_{L^2(0,T)} 
    \leq 
    \sqrt{\frac{T}{3}}\frac{\mu}{\omega_{0,k}},
  \end{equation*}
  and, using that $\Hmu:L^2(0,T)\to L^2(0,T)$ is unitary, we further compute  
  \begin{align*}
    \norm{(\Hmu-\Hmod)s_{0,k}}_{L^2(0,T)} 
    & = \norm{\Hmu(s_{0,k} - s_{\mu,k}) + c_{\mu,k} - c_{0,k}}_{L^2(0,T)} \\
    & \leq \norm{s_{\mu,k} - s_{0,k}}_{L^2(0,T)} + \norm{c_{\mu,k} -c_{0,k}}_{L^2(0,T)}\\
    & \leq 2 \sqrt{\frac{T}{3}}\frac{\mu}{\omega_{0,k}}. 
  \end{align*}
  Now for any $v\in L^2(0,T)$ we have that  
  \begin{align*}
    \norm{(\Hmu-\Hmod)v}_{L^2(0,T)}^2
    & = 
    \left\lVert(\Hmu-\Hmod)\left(\sum_{k=0}^\infty v_k s_{0,k}\right)\right\rVert^2_{L^2(0,T)}\\
    & \leq \left(\sum_{k=0}^\infty |v_k|\norm{(\Hmu-\Hmod)s_{0,k}}_{L^2(0,T)}\right)^2\\
    & \leq \left(\sum_{k=0}^\infty v_k^2\right)\left(\sum_{k=0}^\infty \norm{(\Hmu-\Hmod)s_{0,k}}_{L^2(0,T)}^2\right)\\
    & \leq \left(\sum_{k=0}^\infty v_k^2\right)\left(\frac{4}{3}T\mu^2\sum_{k=0}^\infty \frac{1}{\omega_{0,k}^2}\right).
  \end{align*}
  Using that 
  \begin{equation*}
    \sum_{k=0}^\infty v_k^2 = \frac{2}{T}\norm{v}_{L^2(0,T)}^2, 
  \end{equation*}
  as well as 
  \begin{equation*}
    \sum_{k=0}^\infty \frac{1}{\omega_{0,k}^2}
    =
    \frac{T^2}{\pi^2}\sum_{k=0}^\infty \frac{1}{\left(k+\frac{1}{2}\right)^2}
    =
    \frac{T^2}{2}, 
  \end{equation*}
  finishes the proof of \emph{(ii).}
  
  \emph{(iii).}
  For every fixed $k\in\mathbb N_0$ \eqref{eq:omega as root} implies 
  \[
    \omega_{\mu,k}\longrightarrow
    \omega_{\infty,k}:=\frac{(k+1)\pi}{T}
    \qquad\text{as }\mu\to\infty.
  \]
  The functions 
  \begin{equation*}
    s_{\infty,k}(t) := \sin(\omega_{\infty,k}t) = \sin\left((k+1)\pi \frac{t}{T}\right), \quad k\in \mathbb{N}_0,
  \end{equation*}
  built an orthogonal basis in $L^2(0,T)$. Moreover, we have that 
  \begin{equation*}
    s_{\mu,k}(t) \longrightarrow s_{\infty,k}(t),\quad \text{as }\mu\to \infty,
  \end{equation*}
  and thus 
  \begin{equation*}
    c_{\mu,k}(t) = (-1)^{k}s_{\mu,k}(T-t) \longrightarrow (-1)^k s_{\infty,k}(T-t) = s_{\infty,k}(t),\quad \text{as }\mu\to \infty,
  \end{equation*}
  Again using that $\Hmu:L^2(0,T)\to L^2(0,T)$ is unitary, we get 
  \begin{align*}
    \|(\mathcal H_\mu-I)s_{\infty,k}\|_{L^2(0,T)}
    &\leq
    \|\mathcal H_\mu(s_{\infty,k}-s_{\mu,k})\|_{L^2(0,T)}
    +\|c_{\mu,k}-s_{\infty,k}\|_{L^2(0,T)}
    \\
    &\leq 
    \|s_{\infty,k}-s_{\mu,k}\|_2+\|c_{\mu,k}-s_{\infty,k}\|_{L^2(0,T)}
    \longrightarrow0,
  \end{align*}
  for $\mu\to\infty,$ which finishes the proof. 
\end{proof}
\begin{remark}\label{rem:parabolic-hilbert-no-norm-limit}
  The strong convergence
  $\mathcal H_\mu\to I$ as $\mu\to\infty$
  does not hold in operator norm on $L^2(0,T)$.
  Indeed, for any $\mu \geq 0$, we compute 
  \begin{equation*}
    \int_0^T s_{\mu,k}(t)c_{\mu,k}(t)\, dt 
    =
    \frac{1+\mu T}{2\sqrt{\lambda_{\mu,k}}},
  \end{equation*}
  and therefore by the orthogonality \eqref{eq:orthogonality in B_mu basis} 
  \begin{align*}
    \norm{(\Hmu - I)s_{\mu,k}}_{L^2(0,T)}^2 
    &=
    \norm{\Hmu s_{\mu,k}}_{L^2(0,T)}^2 - 2(\Hmu s_{\mu,k},s_{\mu,k})_{L^2(0,T)} + \norm{s_{\mu,k}}_{L^2(0,T)}^2\\
    & = 2\norm{s_{\mu,k}}_{L^2(0,T)}^2 - 2(c_{\mu,k},s_{\mu,k})_{L^2(0,T)}\\
    & = 2\left(\frac{T}{2} + \frac{\mu}{2\lambda_{\mu,k}}\right) - \frac{1+\mu T}{\sqrt{\lambda_{\mu,k}}}.
  \end{align*}
  Recalling that $\lambda_{\mu,k} = \omega_{\mu,k}^2+\mu^2$ and $\omega_{\mu,k}\to \infty$ as $k\to \infty$
  this gives 
  \begin{equation*}
    \frac{\norm{(\Hmu-I)s_{\mu,k}}_{L^2(0,T)}^2}{\norm{s_{\mu,k}}_{L^2(0,T)}^2} 
    = 2 - \frac{1+\mu T}{\sqrt{\lambda_{\mu,k}}}\left(\frac{T}{2} + \frac{\mu}{2\lambda_{\mu,k}}\right)^{-1}
    \to 2\quad \text{as }k \to \infty. 
  \end{equation*}
  Consequently,
  \[
    \|\mathcal H_\mu-I\|_{L^2(0,T)\to L^2(0,T)}
    \geq\sqrt2
    \qquad\text{for every finite }\mu\geq0.
  \]
  Since any operator-norm limit would also be the
  strong limit, no operator-norm limit exists as
  $\mu\to\infty$.
\end{remark}

\begin{theorem}
  Let $B_\mu := \partial_t +\mu: H^1_{0,}(0,T)\subset L^2(0,T)\to L^2(0,T)$. The natural energy space in this case is given as 
  \begin{equation*}
    X_\mu = \dom(|B_\mu|^{1/2}) = [L^2(0,T),H^1_{0,}(0,T)]_{1/2} = H^{1/2}_{0,}(0,T), 
  \end{equation*}
  equipped with the norm 
  \begin{equation*}
    \norm{v}_{X_\mu} = \big(\norm{v}_{H^{1/2}_{0,}(0,T)}^2 + \mu \norm{v}_{L^2(0,T)}^2\big)^{1/2}.
  \end{equation*}
  Moreover, 
  \begin{equation*}
    B_\mu = \Hmu (-\partial_{tt} + \mu^2)^{1/2}, 
  \end{equation*}
  where $\Hmu:L^2(0,T)\to L^2(0,T)$ is the parabolic Hilbert transformation. 
\end{theorem}
\begin{proof}
  The energy space is again immediate from Theorem \ref{thm:polar_decomposition}. 
  The energy norm is given as 
  \begin{equation*}
    \norm{|B_\mu|^{1/2} v}_{L^2(0,T)}^2 = \sum_{k=0}^\infty \left(\frac{T}{2} + \frac{\mu}{2\lambda_{\mu,k}}\right)\lambda_{\mu,k}^{1/2} v_{\mu,k}^2 = \norm{v}_{H^{1/2}_{0,},\mu}^2, 
  \end{equation*}
  which we have shown to be equivalent to $\norm{v}_{X_\mu}$ in Lemma \ref{lem:parabolic-hilbert-properties}.
  Finally, to obtain the polar decomposition, we compare the action on the basis functions. We compute that 
  \begin{equation*}
    |B_\mu|s_{\mu,k}(t)
    = (B_\mu ^\ast B_\mu)^{1/2} s_{\mu,k}(t) 
    = \sqrt{\lambda_{\mu,k}} s_{\mu,k}(t).
  \end{equation*}
  while 
  \begin{equation*}
    (B_\mu s_{\mu,k})(t) = \omega_{\mu,k} \cos(\omega_{\mu,k} t) + \mu \sin(\omega_{\mu,k}t) = \sqrt{\lambda_{\mu,k}} c_{\mu,k}(t). 
  \end{equation*}
  Thus, the unitary factor has to satisfy 
  \begin{equation*}
    (\Hmu s_{\mu,k})(t) = c_{\mu,k}(t),
  \end{equation*}
  which is exactly the parabolic Hilbert transformation. 
\end{proof}

The polar decomposition immediately leads to a coercive variational formulation of the ODE. 
\begin{theorem}\label{thm:ODE coercive formulation}
  Let $\mu\geq 0$ and $f\in [H^{1/2}_{,0}(0,T)]'$.  
  Then the variational formulation to find $u\in H^{1/2}_{0,}(0,T)$ such that 
  \begin{equation}\label{eq:VF for ODE coercive}
    |b_\mu|(u,v) = \skpr{f,\Hmu v}_{[H^{1/2}_{,0}(0,T)]'\times H^{1/2}_{,0}(0,T)}, \quad \forall v\in H^{1/2}_{0,}(0,T),
  \end{equation}
  where $|b_\mu|(\cdot,\cdot):H^{1/2}_{0,}(0,T)\times H^{1/2}_{0,}(0,T)\to \R$
  \begin{equation*}
    |b_\mu|(u,v) := \skpr{\partial_t u + \mu u,\Hmu v}_{[H^{1/2}_{,0}(0,T)]'\times H^{1/2}_{,0}(0,T)},
  \end{equation*} 
  admits a unique solution, satisfying 
  \begin{equation}
    \norm{u}_{X_\mu} = \big(\norm{u}_{H^{1/2}_{0,}(0,T)}^2+\mu\norm{u}_{L^2(0,T)}^2\big)^{1/2}
    \leq 2 \norm{f}_{[H^{1/2}_{,0}(0,T)]'}. 
  \end{equation}
\end{theorem} 
\begin{proof}
  Since $\Hmu^\ast B_\mu = |B_\mu|$ the bilinear form in \eqref{eq:VF for ODE coercive} is exactly 
  \begin{align*}
    |b_\mu|(u,v) 
    &= \skpr{|B_\mu| u,v}_{[H^{1/2}_{0,}(0,T)]'\times H^{1/2}_{0,}(0,T)}
    = (|B_\mu|^{1/2} u,|B_\mu|^{1/2}v)_{L^2(0,T)}.
  \end{align*}
  Thus, by Lemma \ref{lem:parabolic-hilbert-properties}, we get the ellipticity and boundedness, as  
  \begin{equation}
    |b_\mu|(u,u) 
    =
    \norm{|B_\mu|^{1/2} u}_{L^2(0,T)}^2 
    \geq 
    \frac{1}{\sqrt{2}}\norm{u}_{X_\mu}^2 \label{eq:OdeEllipticity}
  \end{equation}
  and 
  \begin{align}
    |b_\mu|(u,v)
    &=
    (|B_\mu|^{1/2}u,|B_\mu|^{1/2}v)_{L^2(0,T)}
    \leq
    \norm{|B_\mu|^{1/2}u}_{L^2(0,T)}\norm{|B_\mu|^{1/2}v}_{L^2(0,T)}\\
    &\leq 
    2\norm{u}_{X_\mu}\norm{v}_{X_\mu}. \label{eq:OdeBoundedness}
  \end{align}
  Moreover, since for any $w\in H^{1/2}_{,0}(0,T)$ 
  \begin{align*}
    \norm{w}_{H^{1/2}_{,0}(0,T)}^2
    &=
    \frac{T}{2}\sum_{k=0}^\infty \lambda_k^{1/2} \overline{w}_k^2
    \leq 
    \frac{T}{2} \sum_{k=0}^\infty (\lambda_k +\mu^2)^{1/2} \overline{w}_k^2
    \leq \norm{w}_{H^{1/2}_{,0}(0,T),\mu}, 
  \end{align*}
  we have, using the properties of Lemma \ref{lem:parabolic-hilbert-properties}, that 
  \begin{equation*}
    \norm{\Hmu v}_{H^{1/2}_{,0}(0,T)}
    \leq 
    \norm{\Hmu v}_{H^{1/2}_{,0}(0,T),\mu}
    \leq 
    \norm{v}_{H^{1/2}_{0,}(0,T),\mu}
    \leq 
    \sqrt{2}\norm{v}_{X_\mu}
  \end{equation*} 
  and thus 
  \begin{align*}
    \skpr{f,\Hmu v}_{[H^{1/2}_{,0}(0,T)]'\times H^{1/2}_{,0}(0,T)}
    &\leq 
    \norm{f}_{[H^{1/2}_{,0}(0,T)]'}\norm{\Hmu v}_{H^{1/2}_{,0}(0,T)}\\
    &\leq  
    \sqrt{2}\norm{f}_{[H^{1/2}_{,0}(0,T)]'}\norm{v}_{X_\mu}. 
  \end{align*}
  The Lemma of Lax--Milgram guarantees a unique solution and establishes the estimate. 
\end{proof}

\subsection{Discretization, Preconditioning and Numerical Examples}\label{sec:ode_numerical_examples}
Due to the establishes coercivity of the bilinear form $|b_\mu|(\cdot,\cdot)$ for any conforming finite dimensional subspace $X_{h} \subset X_\mu$, the discrete variational formulation to find $u_h\in X_{h}$ such that 
\begin{equation}\label{eq:discrete variational form ode}
  |b_\mu|(u_h,v_h)
  =
  \skpr{f,\Hmu v_h}_{[H^{1/2}_{,0}(0,T)]'\times H^{1/2}_{,0}(0,T)}, \quad \forall v_h\in X_{h},
\end{equation}
admits a unique solution. By standard means using Galerkin orthogonality, one can derive the quasi-best approximation error estimate 
\begin{equation}\label{eq:quasi beast approximation Hmu}
  \norm{u-u_h}_{X_\mu} \leq 2\sqrt{2}\inf_{v_h\in X_{h}} \norm{u-v_h}_{X_\mu}. 
\end{equation}
Anyhow, realizing the numerical implementation of \eqref{eq:discrete variational form ode} is only possible, if $\Hmu$ can be computed. Indeed, for fixed $\mu\geq 0$, one might think of doing so. However, since our aim is a stable discretization of the heat equation later on, where $\mu$ resembles an eigenvalue of the spatial Laplacian, we need a transformation operator, that acts independently of $\mu$. 
In view of the properties of $\Hmu$, established in Lemma \ref{lem:parabolic-hilbert-properties}, we choose the sum of its limiting operators $\widehat{\mathcal{H}} := \Hmod + I$ as an approximation. This was already considered in \cite{LarssonSchwab:2015} for an unbounded domain and in \cite{ReicheltSteinbachHybridSpaceTime2025} for a bounded domain. 
We will show that this choice leads to an unbounded bilinear form in general. However, we are able to give a quasi-best approximation error estimate, up to a logarithmic term. 

Let $0=t_0<t_1<\ldots<t_N=T$ be a quasi-uniform decomposition of $(0,T)$, with mesh sizes $h_i = t_{i}-t_{i-1}$, $i=1,\ldots,N$. For ease of presentation, we restrict our analysis to the case of piecewise linear, globally continuous trial functions $X_{h} = S_{h}^1(0,T)\cap H^1_{0,}(0,T)$ and the discrete variational formulation to find $u_h\in X_{h}$ such that 
\begin{equation}\label{eq:discrete variational form ode perturbed}
  |\widehat{b}_\mu|(u_h,v_h) = \skpr{f,\widehat{\mathcal{H}} v_h}_{[H^{1/2}_{,0}(0,T)]'\times H^{1/2}_{,0}(0,T)},\quad \forall v_h\in X_{h},
\end{equation} 
with the perturbed bilinear form 
\begin{equation*}
  |\widehat{b}_\mu|(u_h,u_h) = (\partial_t u_h+\mu u_h, (\Hmod + I)u_h)_{L^2(0,T)}. 
\end{equation*}
Indeed, the bilinear form is coercive.  
\begin{lemma}\label{lem:coercivity of hat bmu}
  It holds 
  \begin{equation*}
    |\widehat{b}_\mu|(u_h,u_h) 
    \geq 
    \norm{u_h}_{X_\mu}^2,\quad 
    \forall u_h\in X_{h}. 
  \end{equation*}
\end{lemma}
\begin{proof}
  We have that  
\begin{align*}
  |\widehat{b}_\mu|(u_h,u_h) 
  &= (\partial_t u_h,\Hmod u_h)_{L^2(0,T)} 
  + (\partial_t u_h,u_h)_{L^2(0,T)}\\
  &\quad + 
  \mu (u_h,\Hmod u_h)_{L^2(0,T)} 
  +
  \mu (u_h,u_h)_{L^2(0,T)}.
\end{align*}
Due to 
\begin{equation*}
  (\partial_t u_h,u_h)_{L^2(0,T)}
  = \frac{1}{2}\int_0^T \frac{d}{dt}[u_h(t)]^2\, dt 
  =\frac{1}{2}|u_h(T)|^2 \geq 0,
\end{equation*}
and the property of the modified Hilbert transformation, see \cite[Lemma 2.2]{LoescherSteinbachZank:2024}, 
\begin{equation*}
  (u_h,\Hmod u_h)_{L^2(0,T)} \geq 0, 
\end{equation*}
we obtain 
\begin{equation*}
  |\widehat{b}_\mu|(u_h,u_h)
  \geq 
  \norm{u_h}_{H^{1/2}_{0,}(0,T)}^2 + \mu \norm{u_h}_{L^2(0,T)}^2. 
\end{equation*}
\end{proof}

However, the bilinear form is not bounded in the corresponding norm. 
\begin{lemma}
  There does not exist $C>0$ such that
  \begin{equation*}
    |\widehat{b}_\mu|(u,v) \leq C \norm{u}_{H^{1/2}_{0,}(0,T)}\norm{v}_{H^{1/2}_{0,}(0,T)},\quad \forall u,v\in H^{1/2}_{0,}(0,T). 
  \end{equation*}  
\end{lemma} 
\begin{proof}
  Since $\partial_t :H^{1/2}_{0,}(0,T)\to [H^{1/2}_{,0}(0,T)]'$ is an isomorphism, see \cite[Eq. (2.7)]{ZankSteinbachCoerciveSpaceTimeFinite2020} the problematic term is $\skpr{\partial_t u,v}_{[H^{1/2}_{,0}(0,T)]'\times H^{1/2}_{,0}(0,T)}$, when testing with $v\in H^{1/2}_{0,}(0,T)$. To see this, let 
  $v=u\in H^1_{0,}(0,T)$ and consider  
  \begin{equation*}
    (\partial_t u,u)_{L^2(0,T)} 
    =\frac{1}{2}\int_0^T \frac{d}{dt} |u(t)|^2\, dt 
    = 
    \frac{1}{2}|u(T)|^2.  
  \end{equation*}
  Thus boundedness would require 
  \begin{equation}\label{eq:auxiliary bound in H12}
    |u(T)|^2 \leq C \norm{u}_{H^{1/2}_{0,}(0,T)}^2,\quad \forall u\in H^{1/2}_{0,}(0,T).  
  \end{equation} 
  For some fixed $N\in \N$, let us consider 
  \begin{equation*}
    u_N(t) = \sqrt{\frac{2}{\pi L_N}}\sum_{k=0}^N \frac{(-1)^k}{k+\frac{1}{2}} s_k(t) \in H^1_{0,}(0,T), \quad L_N := \sum_{k=0}^N \frac{1}{k+\frac{1}{2}}.
  \end{equation*}
  Then we compute 
  \begin{equation*}
    \norm{u_N}_{H^{1/2}_{0,}(0,T)}^2 
    =
    \frac{T}{2}\frac{2}{\pi L_N}\sum_{k=0}^N \frac{\left(\frac{\pi}{2}+k\pi\right)\frac{1}{T}}{\left(k+\frac{1}{2}\right)^2}
    =
    \frac{1}{L_N} \sum_{k=0}^N \frac{1}{k+\frac{1}{2}} = 1.  
  \end{equation*}
  However, since $s_k(T) = (-1)^k$, we have 
  \begin{align*}
    |u_N(T)| = \sqrt{\frac{2}{\pi L_N}}\sum_{k=0}^N \frac{1}{k+\frac{1}{2}} = \sqrt{\frac{2L_N}{\pi}} \to \infty, \quad \text{as }N\to\infty,
  \end{align*}
  violating \eqref{eq:auxiliary bound in H12}. 
\end{proof}

Anyhow, on finite dimensional spaces, we are able to give a bound of the bilinear form depending on the mesh size $h$, following an estimate established in \cite{McLeanSteinbach1999}.  
\begin{lemma}\label{lem:auxiliary lemma boundedness discrete bilinear form}
  Let $v_h\in X_{h}$ and $w_h(t) = v_h(t)-v_h(T)$. For a quasi-uniform mesh, there holds 
  \begin{equation*}
    \norm{w_h}_{H^{1/2}_{,0}(0,T)} \leq C \left(1+|\log\left(h\right)|\right)\norm{v_h}_{H^{1/2}_{0,}(0,T)}, 
  \end{equation*}
  with a constant $C>0$ independent of $\mu$. 
\end{lemma}
\begin{proof}
   We use the norm representations, see \cite{LionsMagenesVol1}, 
   \begin{align*}
    \norm{v_h}_{H^{1/2}_{0,}(0,T)}^2 
    &\simeq \norm{v_h}_{H^{1/2}(0,T)}^2 + \int_0^T \frac{v_h(t)^2}{t}\,dt,\\
    \norm{w_h}_{H^{1/2}_{,0}(0,T)}^2 
    &\simeq \norm{w_h}_{H^{1/2}(0,T)}^2 + \int_0^T \frac{w_h(t)^2}{T-t}\,dt, 
   \end{align*}
   and the Sobolev-Slobodeckij norm for $u\in H^{1/2}(0,T)$ 
   \begin{equation*}
    \norm{u}_{H^{1/2}(0,T)}^2 
    = \norm{u}_{L^2(0,T)}^2 + |u|_{H^{1/2}(0,T)}^2, 
    \;
    |u|_{H^{1/2}(0,T)}^2 = \int_0^T \int_0^T \frac{|u(t)-u(s)|^2}{|t-s|^2}\, ds \, dt.   
   \end{equation*}
   Since $w_h(T) = 0$, on the last element $w_h(t) = (T-t)w_h(t_{N-1})/h_N$, and thus 
   \begin{equation*}
    \int_{T-h_N}^{T} \frac{w_h(t)^2}{T-t}
    =
    \int_{T-h_N}^{T} \frac{(T-t)w_h(t_{N-1})^2}{h_N^2}\, dt
    =
    \frac{w_h(t_{N-1})^2}{2}.  
   \end{equation*}
   Consequently, as by quasi-uniformity $h_N\simeq h$, 
   \begin{align*}
    \int_0^T \frac{w_h(t)^2}{T-t}\, dt 
    & =  \int_{0}^{T-h_N} \frac{w_h(t)^2}{T-t}\, dt +\frac{w_h(t_{N-1})^2}{2}\\
    & \leq \norm{w_h}_{L^\infty(0,T)}^2\left(\int_0^{T-h_N} \frac{1}{T-t}\, dt + \frac{1}{2}\right)\\
    & = \left(\log\frac{T}{h_N}+ \frac{1}{2}\right)\norm{w_h}_{L^\infty(0,T)}^2\\
    & \leq \left(|\log(T)| + |\log(h_N)|+1\right)\norm{w_h}_{L^\infty(0,T)}^2 \\
    & \leq C(1+|\log(h)|)\norm{w_h}_{L^\infty(0,T)}^2. 
   \end{align*}
   Noting that 
   \begin{equation*}
    \norm{w_h}_{L^2(0,T)}^2 \leq T\norm{w_h}_{L^\infty(0,T)}^2
    = 
    T\norm{v_h-v_h(T)}_{L^\infty(0,T)}^2 
    \leq 
    4T \norm{v_h}_{L^\infty(0,T)}^2, 
   \end{equation*}
   and for the semi norm, it holds that  
   \begin{equation*}
    |w_h|_{H^{1/2}(0,T)} = |v_h|_{H^{1/2}(0,T)}, 
   \end{equation*}
   we get 
   \begin{align*}
    \norm{w_h}_{H^{1/2}_{,0}(0,T)}^2 
    & = 
    \norm{w_h}_{H^{1/2}(0,T)}^2 + \int_0^T \frac{w_h(t)^2}{T-t}\, dt \\
    & = \norm{w_h}_{L^2(0,T)}^2 + |w_h|_{H^{1/2}(0,T)}^2 + \int_0^T \frac{w_h(t)^2}{T-t}\, dt \\
    & \leq 4\big(T + C(1+|\log(h)|)\big)\norm{v_h}_{L^\infty(0,T)}^2 + |v_h|_{H^{1/2}(0,T)}^2 \\
    & \leq C(1+|\log(h)|)\norm{v_h}_{L^\infty(0,T)}^2 + \norm{v_h}_{H^{1/2}(0,T)}^2.  
   \end{align*}
   Further, using the inverse inequality \cite[Lemma 1]{Dryja1988}  
   \begin{equation}\label{eq:inverse estimate H12}
    \norm{u_h}_{L^\infty(0,T)}^2 
    \leq C\left(1+|\log(h)|\right) \norm{u_h}_{H^{1/2}(0,T)}^2,
    \quad \forall u_h\in S_h^1(0,T),
   \end{equation}
   we finally obtain  
   \begin{align*}
      \norm{w_h}_{H^{1/2}_{,0}(0,T)}^2 
      & \leq 
      C \left(1+|\log(h)|\right)^2 \norm{v_h}_{H^{1/2}(0,T)}^2\\
      &\leq C \left(1+|\log(h)|\right)^2 \norm{v_h}_{H^{1/2}_{0,}(0,T)}^2.    
   \end{align*}
\end{proof}

\begin{lemma}\label{lem:discrete boundedness bilinear form}
  There holds 
  \begin{equation*}
    |\widehat{b}_\mu|(u_h,v_h) 
    \leq 
    C(1+|\log(h)|)\norm{u_h}_{X_\mu}\norm{v_h}_{X_\mu},\quad 
    \forall u_h,v_h\in X_{h}. 
  \end{equation*}
\end{lemma}
\begin{proof}
  We have 
  \begin{align*}
  |\widehat{b}_\mu|(u_h,v_h) 
  &= (\partial_t u_h,\Hmod v_h)_{L^2(0,T)} 
  + (\partial_t u_h,v_h)_{L^2(0,T)}\\
  &\quad + 
  \mu (u_h,\Hmod v_h)_{L^2(0,T)} 
  +
  \mu (u_h,v_h)_{L^2(0,T)}.
  \end{align*}
  Using that 
  \begin{equation*}
    (\partial_t u_h,\Hmod v_h)_{L^2(0,T)}\leq \norm{u_h}_{H^{1/2}_{0,}(0,T)}\norm{v_h}_{H^{1/2}_{0,}(0,T)},
  \end{equation*}
  and 
  \begin{equation*}
    (u_h,\Hmod v_h)_{L^2(0,T)} \leq \norm{u_h}_{L^2(0,T)}\norm{v_h}_{L^2(0,T)}, 
  \end{equation*}
  it remains to bound $(\partial_t u_h,v_h)_{L^2(0,T)}$. Defining $w_h(t) = v_h(t)-v_h(T)$, we compute 
  \begin{align*}
    (\partial_t u_h,v_h)_{L^2(0,T)}
    & = (\partial_t u_h,w_h)_{L^2(0,T)} + v_h(T)\int_0^T \partial_t u_h(t)\, dt\\ 
    &= (\partial_t u_h,w_h)_{L^2(0,T)} + u_h(T)v_h(T)\\
    & \leq \norm{\partial_t u_h}_{[H^{1/2}_{,0}(0,T)]'}\norm{w_h}_{H^{1/2}_{,0}(0,T)} + u_h(T)v_h(T)\\
    & = \norm{u_h}_{H^{1/2}_{0,}(0,T)}\norm{w_h}_{H^{1/2}_{,0}(0,T)} + u_h(T)v_h(T).  
  \end{align*}
  From Lemma \ref{lem:auxiliary lemma boundedness discrete bilinear form}, we get 
  \begin{equation*}
    \norm{u_h}_{H^{1/2}_{0,}(0,T)}\norm{w_h}_{H^{1/2}_{,0}(0,T)}
    \leq 
    C(1+|\log(h)|) \norm{u_h}_{H^{1/2}_{0,}(0,T)}\norm{v_h}_{H^{1/2}_{0,}(0,T)}. 
  \end{equation*}
  Moreover, using the inverse inequality \eqref{eq:inverse estimate H12}, we can estimate 
  \begin{align*}
    u_h(T)v_h(T) 
    &\leq 
    \norm{u_h}_{L^\infty(0,T)}\norm{v_h}_{L^\infty(0,T)}\\
    &\leq 
    C(1+|\log(h)|)\norm{u_h}_{H^{1/2}(0,T)}\norm{v_h}_{H^{1/2}(0,T)}\\
    &\leq 
    C(1+|\log(h)|)\norm{u_h}_{H^{1/2}_{0,}(0,T)}\norm{v_h}_{H^{1/2}_{0,}(0,T)},
  \end{align*}
  which finishes the proof. 
\end{proof}

The main result is the following. 
\begin{theorem}\label{thm:DVF perturbed ODE main result}
  For any $f\in [H^{1-s}_{,0}(0,T)]'$, $s>1/2$ the discrete variational formulation \eqref{eq:discrete variational form ode perturbed} admits a unique solution $u_h\in X_{h}$. Moreover, if the unique solution of 
  \begin{equation}\label{eq:regular solution u in ode form}
    \partial_t u + \mu u = f \quad \text{in }[H^{1-s}_{,0}(0,T)]',  
  \end{equation} 
  satisfies  $u\in H^{1/2}_{0,}(0,T)\cap H^s(0,T)$, $s\in (1/2,2]$, there holds the error estimate 
  \begin{equation*}
    \norm{u-u_h}_{X_{\mu}} \leq C h^{s-\frac{1}{2}}\sqrt{1+\mu h}\norm{u}_{H^s(0,T)}. 
  \end{equation*}
\end{theorem}
\begin{proof}
Since $(\mathcal H_T+I)X_{h}\subset H^1(0,T)\subset H^{1-s}(0,T) = H^{1-s}_{,0}(0,T),$ for any $s>1/2$, see \cite[Theorem 3.40]{McLeanStronglyEllipticSystems2000}, 
the functional
\[
  F_h(v_h):=\langle f,(\mathcal H_T+I)v_h\rangle_{[H^{1-s}_{,0}(0,T)]'\times H^{1-s}_{,0}(0,T)}
\]
is well defined and continuous. Moreover, the bilinear form $|\widehat{b}_\mu|:X_{h}\times X_{h}\to \R$ is coercive (Lemma \ref{lem:coercivity of hat bmu}) and bounded for any fixed $h>0$ (Lemma \ref{lem:discrete boundedness bilinear form}). Consequently, by the Lemma of Lax--Milgram there exists a
unique solution $u_h\in X_{h}$.

To prove the error estimate let $I_h^1:C([0,T])\to X_{h}$ be the nodal
piecewise linear interpolant. Recall, $H^s_{0,}(0,T)\hookrightarrow C([0,T])$ for $s>1/2$. 
Thus we can split the error
\begin{equation*}
  \norm{u-u_h}_{X_\mu} 
  \leq 
  \norm{u-I_h^1 u}_{X_\mu} + \norm{I_h^1 u - u_h}_{X_\mu}. 
\end{equation*}
Again, since $(\Hmod +I)(X_{h}) \subset H^{1-s}_{,0}(0,T)$, testing \eqref{eq:regular solution u in ode form} with $(\Hmod + I)v_h$ is well-defined, and we obtain the Galerkin orthogonality
\begin{equation*}
  |\widehat{b}_\mu|(u-u_h,v_h) = 0,\quad \forall v_h\in X_{h}. 
\end{equation*}
Together with the coercivity of Lemma \ref{lem:coercivity of hat bmu}, we obtain 
\begin{equation}\label{eq:auxiliary coercivity}
  \norm{I_h^1 u-u_h}_{X_\mu}^2 
  \leq 
  |\widehat{b}_\mu|(I_h^1 u-u_h, I_h^1u-u_h) 
  =
  |\widehat{b}_\mu|(I_h^1u-u,I_h^1u-u_h).  
\end{equation}
It remains to bound the bilinear form. Therefore note that $\eta(t) := u(t)-I_h^1u(t)$ satisfies $\eta(0)=\eta(T)=0$, implying 
\begin{equation*}
  \skpr{\partial_t \eta,1}_{[H^{1-s}_{,0}(0,T)]'\times H^{1-s}_{,0}(0,T)}
  = 
  \eta(T)-\eta(0) =0. 
\end{equation*} 
Consequently, using integration by parts and the inverse inequality 
\begin{equation}\label{eq:inverse inequality in time h1 to h12}
  \norm{\partial_t v_h}_{L^2(0,T)} \leq C h^{-1/2}\norm{v_h}_{H^{1/2}_{0,}(0,T)},  
\end{equation}
see \cite[(4.17)]{FaustmannMelenkPraetorius2021}, we compute  
\begin{align*}
  &\skpr{\partial_t \eta, v_h}_{[H^{1-s}_{,0}(0,T)]'\times H^{1-s}_{,0}(0,T)}
  = \skpr{\partial_t \eta, w_h}_{[H^{1-s}_{,0}(0,T)]'\times H^{1-s}_{,0}(0,T)} \\
  &= -(\eta, \partial_t w_h)_{L^2(0,T)}  = -(\eta, \partial_t v_h)_{L^2(0,T)} \\
  &\leq \norm{\eta}_{L^2(0,T)}\norm{\partial_t v_h}_{L^2(0,T)}\\
  &\leq ch^{-1/2}\norm{\eta}_{L^2(0,T)}\norm{v_h}_{H^{1/2}_{0,}(0,T)}. 
\end{align*}
Furthermore, we have 
\begin{align*}
  &|\widehat{b}_\mu|(u-I_h^1 u,v_h)
   = |\widehat{b}_\mu|(\eta,v_h)\\
  &\quad  = \skpr{\partial_t \eta, \Hmod v_h}_{[H^{1-s}_{,0}(0,T)]'\times H^{1-s}_{,0}(0,T)} + \skpr{\partial_t \eta, v_h}_{[H^{1-s}_{,0}(0,T)]'\times H^{1-s}_{,0}(0,T)} \\
  & \qquad + \mu(\eta,\Hmod v_h)_{L^2(0,T)} + \mu (\eta, v_h)_{L^2(0,T)}\\
  & \quad \leq \norm{\eta}_{H^{1/2}_{0,}(0,T)}\norm{v_h}_{H^{1/2}_{0,}(0,T)} + ch^{-1/2}\norm{\eta}_{L^2(0,T)}\norm{v_h}_{H^{1/2}_{0,}(0,T)}\\
  &\qquad +2\mu\norm{\eta}_{L^2(0,T)}\norm{v_h}_{L^2(0,T)}\\ 
  &\quad  \leq C\big(\norm{\eta}_{H^{1/2}_{0,}(0,T)} + (h^{-1/2}+\sqrt{\mu})\norm{\eta}_{L^2(0,T)}\big)\norm{v_h}_{X_\mu}.  
\end{align*}
Thus, choosing $v_h = I_h^1 u-u_h$, by \eqref{eq:auxiliary coercivity}, we get  
\begin{equation*}
  \norm{I_h^1 u-u_h}_{X_\mu} 
  \leq 
  C\big(\norm{u-I_h^1 u}_{H^{1/2}_{0,}(0,T)} + (h^{-1/2} + \sqrt{\mu})\norm{u-I_h^1 u}_{L^2(0,T)}\big).   
\end{equation*}
By standard interpolation error estimates 
\begin{align*}
  \|u-I_h^1 u\|_{L^2(0,T)}
  &\leq C h^s\|u\|_{H^s_{0,}(0,T)},\\
  \|u-I_h^1 u\|_{H^{1/2}_{0,}(0,T)}
  &\leq C h^{s-1/2}\|u\|_{H^s_{0,}(0,T)}.
\end{align*}
the result follows. 
\end{proof}

\begin{remark}
  In order to establish the error estimate in Theorem \ref{thm:DVF perturbed ODE main result}, we need to assume that the solution satisfies $u\in H^{1/2}_{0,}(0,T)\cap H^s(0,T)$ for $s> 1/2$, preventing quasi-optimality. Indeed we cannot have a quasi-best approximation error estimate, as using the transformation $\widehat{\mathcal{H}}$ the test space is non-conforming, see \cite[Remark 4.9]{VeeserZanotti2018}. 
\end{remark}

Note, that $X_{h} \subset H^{s}_{0,}(0,T)$ for $s < 3/2$. To translate the variational problem into a linear system, we choose the Lagrange basis functions $\lbrace \varphi_i \rbrace_{i=1}^{N}$ for $X_h$ to obtain 
\begin{align}
  \mathbf{B}_\mu \underline{u} = \underline{f}, \label{eq:ode_linear_system}
\end{align}
where the system matrix $\mathbf{B}_\mu \in \R^{N \times N}$ and the right-hand side $\underline{f} \in \R^{N}$ are given by
\begin{align*}
  (\mathbf{B}_\mu)_{i,j} &= \lvert \hat{b}_\mu \rvert(\varphi_j, \varphi_i) 
  = (\partial_t \varphi_j + \mu \varphi_j, \mathcal{H}_T \varphi_i + \varphi_i )_{L^2(0,T)}, \nonumber \\
  &= (\partial_t \varphi_j, \mathcal{H}_T \varphi_i )_{L^2(0,T)} + \mu (\varphi_j, \varphi_i )_{L^2(0,T)}\\
  &\quad  + (\partial_t \varphi_j, \varphi_i )_{L^2(0,T)} + \mu ( \varphi_j, \mathcal{H}_T \varphi_i )_{L^2(0,T)}, \nonumber\\
  \underline{f}_i &= l_h(\varphi_i) = \langle f, \mathcal{H}_T \varphi_i + \varphi_i \rangle_{[H^{1-s}_{,0}(0,T)]'\times H^{1-s}_{,0}(0,T)}.
\end{align*}
Note, that the terms not involving the modified Hilbert transform are standard finite element matrices. The dense matrices and vectors involving the modified Hilbert transform can be efficiently computed using series expansions based on the Legendre chi function \cite{ZankExactRealizationModified2021} and numerical quadrature of weakly singular integral representations \cite{Zank2023} for arbitrary meshes. More recent FFT techniques for uniform meshes even outperform these approaches and are applicable for arbitrary polynomial degrees. Open source software for this approach can be found in \cite{Reichelt_hmod_2025}. 

\subsubsection{Preconditioning}
To solve the linear system \eqref{eq:ode_linear_system}, we aim to use a preconditioned iterative solver. Note, that the system matrix $\mathbf{B}_\mu$ is non-symmetric but positive definite. From Lemma \ref{lem:discrete boundedness bilinear form} we know that the bilinear form $\lvert \hat{b}_\mu \rvert$ is elliptic and bounded on the discrete space $X_{h}$ up to a logarithmic factor with respect to the $X_\mu$ norm.
Moreover, the norm itself is induced by the operator
\begin{align}
  D_\mu \, : \, X_\mu &\rightarrow X_\mu', \nonumber \\
  \langle D_\mu u, v \rangle_{X_\mu'\times X_\mu} &:= ( u , v )_{X_\mu} = \langle \partial_t u, \mathcal{H}_T v \rangle_{[H_{,0}^{1/2}(0,T)]'\times H^{1/2}_{,0}(0,T)} + \mu ( u , v )_{L^2(0,T)}.
\end{align}
Its discrete counterpart is the matrix $\mathbf{D}_\mu \in \mathbb{R}^{N\times N}$ given as
\begin{align}
  (\mathbf{D}_\mu)_{i,j} &= (\partial_t \varphi_j,\mathcal{H}_T \varphi_i)_{L^2(0,T)} + \mu (\varphi_j,\varphi_i)_{L^2(0,T)},
\end{align}
which can be used as a preconditioner for $\mathbf{B}_\mu$. By Lemma \ref{lem:coercivity of hat bmu} and Lemma \ref{lem:auxiliary lemma boundedness discrete bilinear form}, we immediately get the following result. 
\begin{lemma}\label{lem:spectral_equivalence_ode}
  The matrices $\mathbf{B}_\mu$ and $\mathbf{D}_\mu$ satisfy the spectral equivalence
  \begin{align}
    \underline{v}^T \mathbf{D}_\mu \, \underline{v} 
    \leq \underline{v}^T \mathbf{B}_\mu \, \underline{v} 
    \leq C(1+|\log(h)|) \, \underline{v}^T \mathbf{D}_\mu \, \underline{v}, \quad \forall \underline{v} \in \R^{N},
  \end{align}
  where the constant $C > 0$ is independent of $h$ and $\mu$.
\end{lemma}
We will now proceed and first show, how a symmetric preconditioning will affect the number of GMRES iterations to solve the preconditioned system. Afterward, we will also show numerical results for the left preconditioned system. The results for the symmetric preconditioning are summarized in the following theorem.

\begin{theorem}\label{thm:SymmetricalPreconditioning}
  Consider the preconditioned system
  \begin{align}
    \mathbf{\tilde{B}} \, \underline{w} = \underline{g},
  \end{align}
  where
  \begin{align}
    \mathbf{\tilde{B}} = \mathbf{D}_\mu^{-1/2} \, \mathbf{B}_\mu \, \mathbf{D}_\mu^{-1/2}, \quad \underline{w} = \mathbf{D}_\mu^{1/2} \, \underline{u}, \quad \underline{g} = \mathbf{D}_\mu^{-1/2} \, \underline{f}.
  \end{align}
  Then the $k$-th GMRES residual $\underline{r}_k = \underline{g} - \mathbf{\tilde{B}}\underline{w}_k$ satisfies
  \begin{align}
      \frac{\lVert \underline{r}_k \rVert_2}{\lVert \underline{r}_0 \rVert_2} \leq \left(1 - \frac{1}{C^2(1+\lvert \log(h) \rvert)^2}\right)^{k/2},
  \end{align}
  where the constant $C > 0$ is from the boundedness of Lemma \ref{lem:discrete boundedness bilinear form} and $\lVert \cdot \rVert_2$ is the standard Euclidean norm.
  \begin{proof}
    First, by setting $\underline{v} = \mathbf{D}^{1/2}_\mu \, \underline{y}$ we obtain from the spectral equivalence
    \begin{align}
      (\underline{v}, \underline{v} )_2 
      &\leq (\mathbf{D}_\mu \underline{y}, \underline{y} )_2 
      \leq (\mathbf{B}_\mu \underline{y}, \underline{y} )_2 
      = ( \mathbf{\tilde{B}} \underline{v}, \underline{v} )_2 
      \leq \xi(h) (\mathbf{D}_\mu \underline{y}, \underline{y})_2 
      = \xi(h) (\underline{v}, \underline{v} )_2
    \end{align}
    for arbitrary $\underline{y} \in \R^{N}$, where $( \cdot, \cdot )_2$ is the standard Euclidean inner product and $\xi(h) = C(1+\lvert \log(h) \rvert)$. This means that the field of values of $\mathbf{\tilde{B}}$ is contained in $[1, \xi(h)]$, i.e.
    \begin{align}
      W(\mathbf{\tilde{B}}) 
      = \left\lbrace \frac{( \mathbf{\tilde{B}} \underline{v}, \underline{v} )_2}{\lVert \underline{v} \rVert_2^2} : \underline{v} \in \R^{N} \backslash \lbrace 0 \rbrace \right\rbrace \subseteq [1, \xi(h)].
    \end{align}
    Further note, that
    \begin{align}
      \lVert \mathbf{\tilde{B}} \rVert_2 &= \sup_{\underline{v}, \underline{x} \in \R^{N} \backslash \lbrace 0 \rbrace} \frac{(  \mathbf{\tilde{B}} \underline{v}, \underline{x})_2}{\lVert \underline{v} \rVert_2 \lVert \underline{x} \rVert_2} = \sup_{\underline{v}, \underline{x} \in \R^{N} \backslash \lbrace 0 \rbrace} \frac{(  \mathbf{\tilde{B}} \mathbf{D}_\mu^{1/2}\underline{v}, \mathbf{D}_\mu^{1/2} \underline{x})_2}{\lVert \mathbf{D}_\mu^{1/2} \underline{v} \rVert_2 \lVert \mathbf{D}_\mu^{1/2} \underline{x} \rVert_2} \nonumber \\
      &\leq \sup_{\underline{v}, \underline{x} \in \R^{N} \backslash \lbrace 0 \rbrace} \frac{(  \mathbf{B}_\mu \underline{v}, \underline{x})_2}{\lVert \mathbf{D}_\mu^{1/2} \underline{v} \rVert_2 \lVert \mathbf{D}_\mu^{1/2} \underline{x} \rVert_2} \leq  \xi(h) \sup_{\underline{v}, \underline{x} \in \R^{N} \backslash \lbrace 0 \rbrace} \frac{\lVert \mathbf{D}_\mu^{1/2} \underline{v} \rVert_2 \lVert \mathbf{D}_\mu^{1/2} \underline{x} \rVert_2}{\lVert \mathbf{D}_\mu^{1/2} \underline{v} \rVert_2 \lVert \mathbf{D}_\mu^{1/2} \underline{x} \rVert_2} \nonumber \\
      &\leq \xi(h)
    \end{align}
    where we have used the boundedness of Lemma \ref{lem:discrete boundedness bilinear form}. Note, that this is the same constant as in the spectral equivalence of Lemma \ref{lem:spectral_equivalence_ode}. We can now use the Elman-Starke-Beckermann style convergence result for GMRES from \cite{BeckermannGoreinovTyrtyshnikovRemarks2005} to obtain
    \begin{align}
      \frac{\lVert \underline{r}_k \rVert_2}{\lVert \underline{r}_0 \rVert_2} \leq \left(1 - \frac{\nu^2}{\lVert \tilde{\mathbf{B}_\mu} \rVert_2^2}\right)^{k/2} \leq \left(1 - \frac{1}{\xi^2(h)}\right)^{k/2},
    \end{align}
    where $\underline{r}_k$ is the residual after $k$ GMRES iterations and 
    \begin{align}
      \nu = \dist\left(0, W(\mathbf{\tilde{B}})\right) \geq 1.
    \end{align}
  \end{proof}
\end{theorem}
This result shows that the number of preconditioned GMRES iterations to reduce the residual by a fixed factor is only dependent on the logarithm of the mesh size $h$ and independent of the parameter $\mu$. The numerical examples below report the iteration counts under mesh refinement for both symmetric norm preconditioning and BPX.

\begin{figure}[H]
  \centering
  \begin{tikzpicture}

    \begin{groupplot}[
        group style={
            group size=2 by 1,
            horizontal sep=2cm,
        },
        xlabel={$\mu$},
        ylabel={GMRES iterations},
        xmode=log, 
        ymax = 220,
        ymin = 0,
        legend style={
            font=\small,
            draw=none,
        },
        legend style={
            draw=black,        
            fill=white,        
            fill opacity=1,  
            text opacity=1,    
            font=\small,
        },
        legend columns=\ODELegendColumns,
        legend cell align=left,
        width=\textwidth*0.488,
    ]

    \nextgroupplot[
        legend to name=commonlegend,
    ]
        \addplot[no marks, black, very thick] table [x=mu, y=sPrec, col sep=comma, mark=none]
            {simulations/ode/gmres_iterations_symmetric_pol1.csv};
        \addlegendentry{$\widetilde{\mathcal{H}}  = I$}

        \addplot[no marks, dashed, red, very thick] table [x=mu, y=hybPrec, col sep=comma, mark=none]
            {simulations/ode/gmres_iterations_symmetric_pol1.csv};
        \addlegendentry{$\widetilde{\mathcal{H}} = \mathcal{H}_T + I$}

        \addplot[no marks, dotted, blue, very thick] table [x=mu, y=hilbPrec, col sep=comma, mark=none]
            {simulations/ode/gmres_iterations_symmetric_pol1.csv};
        \addlegendentry{$\widetilde{\mathcal{H}} = \mathcal{H}_T$}

    \nextgroupplot[
        legend to name=commonlegendright,
    ]
        \addplot[no marks, black, very thick] table [x=mu, y=standard, col sep=comma, mark=none]
            {simulations/ode/gmres_iterations_symmetric_pol1.csv};
        \addlegendentry{$\widetilde{\mathcal{H}} = I$}

        \addplot[no marks, dashed, red, very thick] table [x=mu, y=hyb, col sep=comma, mark=none]
            {simulations/ode/gmres_iterations_symmetric_pol1.csv};
        \addlegendentry{$\widetilde{\mathcal{H}} = \mathcal{H}_T + I$}

        \addplot[no marks, dotted, blue, very thick] table [x=mu, y=hilb, col sep=comma, mark=none]
            {simulations/ode/gmres_iterations_symmetric_pol1.csv};
        \addlegendentry{$\widetilde{\mathcal{H}} = \mathcal{H}_T$}
    \end{groupplot}

  \PlaceODELegend

\end{tikzpicture}
  \caption{Number of GMRES iterations for the symmetrically preconditioned system (left) and unpreconditioned system (right) for different values of $\mu$ and different choices of $\widetilde{\mathcal{H}}$, where $X_{h} = S_{h}^1(0,T)\cap H^1_{0,}(0,T)$ with $N = 200$ intervals.}
  \label{fig:ode_gmres_iterations_symm}
\end{figure}
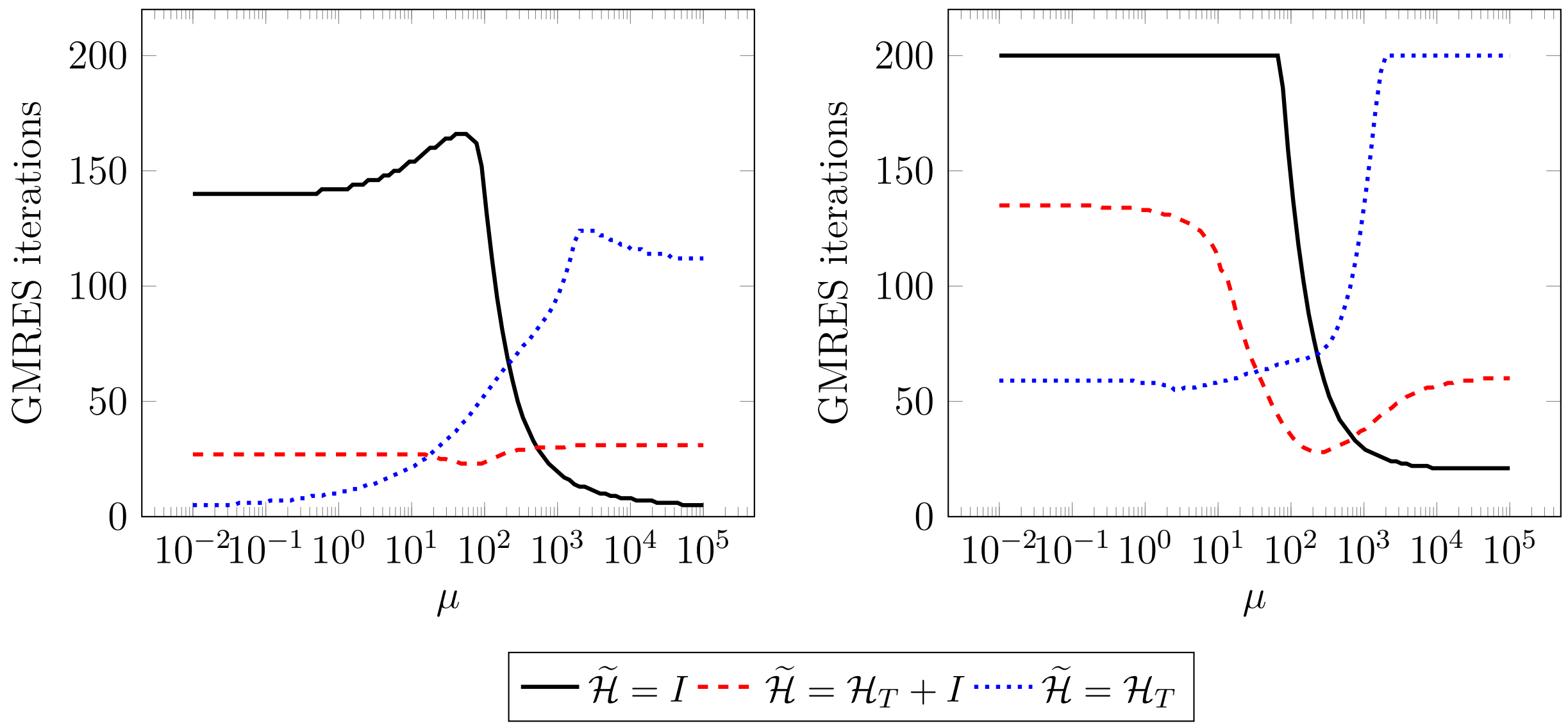

To illustrate, that we indeed need to consider both summands in the transformation operator $\widehat{\mathcal{H}} = \mathcal{H}_T + I$ to be able to precondition the system matrix $\mathbf{B}_\mu$ robustly w.r.t. the parameter $\mu$, we show numerical results for different choices of the transformation. For this purpose we define the system matrix
\begin{align}
  (\mathbf{B}_\mu)_{i,j} = (\partial_t \varphi_j + \mu \varphi_j, \widetilde{\mathcal{H}} \varphi_i )_{L^2(0,T)},
\end{align}
where $\widetilde{\mathcal{H}}$ is chosen as either $\widetilde{\mathcal{H}} = I$, corresponding to the formulation of \cite{Steinbach2015},  $\widetilde{\mathcal{H}} = \mathcal{H}_T$, corresponding to \cite{ZankSteinbachCoerciveSpaceTimeFinite2020}, or $\widetilde{\mathcal{H}} = \widehat{\mathcal{H}}$. Figure \ref{fig:ode_gmres_iterations_symm} shows the number of GMRES iterations for the symmetrically preconditioned system and the unpreconditioned system for different values of $\mu$ and the different choices of $\widetilde{\mathcal{H}}$. The initial guess for the GMRES is the zero vector, the right-hand side is set to one and the algorithm stops if the initial residual is reduced by a factor $10^{-12}$. We observe, that for the first choice, that is stable w.r.t. the $\mu$ part of the norm, the number of GMRES iterations decreases with growing $\mu$. For the second choice, that is stable w.r.t. the $H^{1/2}$ part of the norm, the number of GMRES iterations grows with growing $\mu$. Only for the third choice, that is stable w.r.t. both parts of the norm, the number of GMRES iterations remains constant for all values of $\mu$.The robustness of the preconditioner w.r.t. mesh refinement will be illustrated when considering numerical examples later.

\paragraph{BPX Preconditioner}
The symmetric preconditioner in Theorem \ref{thm:SymmetricalPreconditioning} requires the inverse square root of the dense norm matrix $\mathbf{D}_\mu$. A multilevel representation of the norm yields a BPX preconditioner for the system matrix $\mathbf{B}_\mu$.
Consider a sequence of uniformly refined meshes of the interval $(0, T)$ with mesh sizes $h_{j}$ for $j = 0, 1, \ldots $ and the corresponding nested discrete spaces
\begin{align}
  X_{h_0} \subset X_{h_1} \subset \ldots \subset X_{\mu}.
\end{align}
The operators
\begin{align}
  B^{1/2} = \sum_{j=0}^\infty h_{j}^{-1} \left( Q_{h_j}^1 - Q_{h_{j-1}}^1 \right), \quad 
  B^0 = \sum_{j=0}^\infty \left( Q_{h_j}^1 - Q_{h_{j-1}}^1 \right)
\end{align}
are equivalent to the operators $\mathcal{H}_T^{-1}\partial_t$ and $I$, respectively, where $Q_{h_j}^1$ is the $L^2$-projection onto the space $X_{h_j}$ and $Q_{h_{-1}}^1 = 0$.
So our norm inducing operator $D = \mathcal{H}_T^{-1}\partial_t + \mu I$ is spectrally equivalent to the BPX operator
\begin{align}
  B_\mu = B^{1/2} + \mu B^0 = \sum_{j=0}^\infty (h_{j}^{-1} +\mu) \left( Q_{h_j}^1 - Q_{h_{j-1}}^1 \right).
\end{align}
Denoting by $Q_j=Q_{h_j}^1$ and $w_j=(h_j^{-1}+\mu)^{-1}$ the differences $Q_j-Q_{j-1}$ project onto mutually orthogonal detail spaces, and the inverse is given as 
\begin{align}
  B_\mu^{-1} = \sum_{j=0}^\infty w_j (Q_j-Q_{j-1}).
\end{align}
On $X_{h_L}$ the sum ends at $j=L$. Collecting the coefficients of each $Q_j$ gives $w_j-w_{j+1}$ for $j<L$ and $w_L$ for $Q_L=I$.
Let $\mathbf{R}_j$ denote prolongation from $X_{h_j}$ to $X_{h_L}$ and $\mathbf{M}_j$ the level-$j$ mass matrix. The coefficient matrix of $Q_j$ is $\mathbf{R}_j\mathbf{M}_j^{-1}\mathbf{R}_j^\top\mathbf{M}_L$. Thus the discrete multilevel inverse acting on load vectors is
\begin{align}
  \mathbf{S}_{\mu,L}
  = w_L\mathbf{M}_L^{-1}
    + \sum_{j=0}^{L-1}(w_j-w_{j+1})\mathbf{R}_j\mathbf{M}_j^{-1}\mathbf{R}_j^\top .
  \label{eq:ode_bpx_inverse}
\end{align}
Since $B^0=I$, the mass term $\mu I$ enters both norm operators with the same weight. The multilevel equivalence for the temporal term therefore gives $\mathbf{S}_{\mu,L}\simeq\mathbf{D}_\mu^{-1}$ uniformly in $L$ and $\mu\geq 0$.

\subsubsection{Numerical Examples} \label{sec:ode:numerical_examples}
We now turn to numerical examples to illustrate the approximation properties as well as the efficiency of the preconditioner. We consider $T=1$
and different values of $\mu > 0$, as well as solutions of different regularity. Namely, 
\begin{align}
  u_r(t) = \sin(\pi t) \in C^\infty(\overline{I}),
\end{align}
to obtain a regular solution, and 
\begin{align}
  u_s(t) = \sqrt{t} \in H^{s}_{0,}(0,T),\quad s<1, 
\end{align}
to obtain a solution with reduced regularity, obeying a singularity at $t=0$. 
We discretize the problem using the piecewise polynomial spaces $X_{h} = S_{h}^p(0,T)\cap H^1_{0,}(0,T)$ of degrees $p = 1,2,3,5$ on a sequence of uniformly refined meshes with $N$ intervals using the variational formulation \eqref{eq:discrete variational form ode perturbed}.

\begin{figure}[h!]
  \centering
  \def\csvfolder{simulations/ode/big_convergence_study_ode_hybrid_smooth/}
  \def\muVal{10}
  \def\yminVal{0.5*1e-13}
  \begin{tikzpicture}
  \small
  \begin{groupplot}[
    group style={
      group size=2 by 2,       
      horizontal sep=1.8cm,    
      vertical sep=1.8cm,      
    },
    width=0.385\linewidth,      
    height=0.33\linewidth,
    scale only axis,
    grid=both,
  ]

    \nextgroupplot[
      xlabel={$N$},
      ylabel={$\|u - u_h\|_{X_\mu}$},
      ylabel near ticks,
      ylabel shift = -7pt,
      ymax = 0.5*1e1,
      ymin = \yminVal,
      ymode=log, xmode=log,
      legend pos=north east,
    ]
    \addplot+[
      thick,
      mark=*,
    ]
    table[
      col sep=comma,
      x expr={abs(\thisrow{mu}-\muVal)<1e-5 ? \thisrow{nt} : nan},
      y expr={abs(\thisrow{mu}-\muVal)<1e-5 ? \thisrow{Xerr} : nan},
    ]{\csvfolder/results_p1.csv};
    \addlegendentry{$p_t=1$}

    \addplot+[
      thick,
      mark=x,
    ]
    table[
      col sep=comma,
      x expr={abs(\thisrow{mu}-\muVal)<1e-5 ? \thisrow{nt} : nan},
      y expr={abs(\thisrow{mu}-\muVal)<1e-5 ? \thisrow{Xerr} : nan},
    ]{\csvfolder/results_p2.csv};
    \addlegendentry{$p_t=2$}

    \addplot+[
      thick,
      mark=o,
    ]
    table[
      col sep=comma,
      x expr={abs(\thisrow{mu}-\muVal)<1e-5 ? \thisrow{nt} : nan},
      y expr={abs(\thisrow{mu}-\muVal)<1e-5 ? \thisrow{Xerr} : nan},
    ]{\csvfolder/results_p3.csv};
    \addlegendentry{$p_t=3$}

    \addplot+[
      thick,
      mark=star,
    ]
    table[
      col sep=comma,
      x expr={abs(\thisrow{mu}-\muVal)<1e-5 ? \thisrow{nt} : nan},
      y expr={abs(\thisrow{mu}-\muVal)<1e-5 ? \thisrow{Xerr} : nan},
    ]{\csvfolder/results_p5.csv};
    \addlegendentry{$p_t=5$}

    \nextgroupplot[
      xlabel={$N$},
      ylabel={$\|u - u_h\|_{L^2}$},
      ylabel near ticks,
      ylabel shift = -7pt,
      ymax = 0.5*1e1,
      ymin = \yminVal,
      ymode=log, xmode=log,
      legend pos=north east,
    ]
    \addplot+[
      thick,
      mark=*,
    ]
    table[
      col sep=comma,
      x expr={abs(\thisrow{mu}-\muVal)<1e-5 ? \thisrow{nt} : nan},
      y expr={abs(\thisrow{mu}-\muVal)<1e-5 ? \thisrow{L2err} : nan},
    ]{\csvfolder/results_p1.csv};
    \addlegendentry{$p_t=1$}

    \addplot+[
      thick,
      mark=x,
    ]
    table[
      col sep=comma,
      x expr={abs(\thisrow{mu}-\muVal)<1e-5 ? \thisrow{nt} : nan},
      y expr={abs(\thisrow{mu}-\muVal)<1e-5 ? \thisrow{L2err} : nan},
    ]{\csvfolder/results_p2.csv};
    \addlegendentry{$p_t=2$}

    \addplot+[
      thick,
      mark=o,
    ]
    table[
      col sep=comma,
      x expr={abs(\thisrow{mu}-\muVal)<1e-5 ? \thisrow{nt} : nan},
      y expr={abs(\thisrow{mu}-\muVal)<1e-5 ? \thisrow{L2err} : nan},
    ]{\csvfolder/results_p3.csv};
    \addlegendentry{$p_t=3$}

    \addplot+[
      thick,
      mark=star,
    ]
    table[
      col sep=comma,
      x expr={abs(\thisrow{mu}-\muVal)<1e-5 ? \thisrow{nt} : nan},
      y expr={abs(\thisrow{mu}-\muVal)<1e-5 ? \thisrow{L2err} : nan},
    ]{\csvfolder/results_p5.csv};
    \addlegendentry{$p_t=5$}

    \nextgroupplot[
      xlabel={$N$},
      ylabel={$G_{\text{sym}}$},
      ylabel near ticks,
      ylabel shift = -7pt,
      xmode=log,
      legend pos=south east,
    ]
    \addplot+[
      thick,
      mark=*,
    ]
    table[
      col sep=comma,
      x expr={abs(\thisrow{mu}-\muVal)<1e-5 ? \thisrow{nt} : nan},
      y expr={abs(\thisrow{mu}-\muVal)<1e-5 ? \thisrow{Gsymm} : nan},
    ]{\csvfolder/results_p1.csv};
    \addlegendentry{$p_t=1$}

    \addplot+[
      thick,
      mark=x,
    ]
    table[
      col sep=comma,
      x expr={abs(\thisrow{mu}-\muVal)<1e-5 ? \thisrow{nt} : nan},
      y expr={abs(\thisrow{mu}-\muVal)<1e-5 ? \thisrow{Gsymm} : nan},
    ]{\csvfolder/results_p2.csv};
    \addlegendentry{$p_t=2$}

    \addplot+[
      thick,
      mark=o,
    ]
    table[
      col sep=comma,
      x expr={abs(\thisrow{mu}-\muVal)<1e-5 ? \thisrow{nt} : nan},
      y expr={abs(\thisrow{mu}-\muVal)<1e-5 ? \thisrow{Gsymm} : nan},
    ]{\csvfolder/results_p3.csv};
    \addlegendentry{$p_t=3$}

    \addplot+[
      thick,
      mark=star,
    ]
    table[
      col sep=comma,
      x expr={abs(\thisrow{mu}-\muVal)<1e-5 ? \thisrow{nt} : nan},
      y expr={abs(\thisrow{mu}-\muVal)<1e-5 ? \thisrow{Gsymm} : nan},
    ]{\csvfolder/results_p5.csv};
    \addlegendentry{$p_t=5$}

    \nextgroupplot[
      xlabel={$N$},
      ylabel={$G_{\text{bpx}}$},
      ylabel near ticks,
      ylabel shift = -7pt,
      xmode=log,
      legend pos=south east,
    ]
    \addplot+[
      thick,
      mark=*,
    ]
    table[
      col sep=comma,
      x expr={abs(\thisrow{mu}-\muVal)<1e-5 ? \thisrow{nt} : nan},
      y expr={abs(\thisrow{mu}-\muVal)<1e-5 ? \thisrow{Gbpx} : nan},
    ]{\csvfolder/results_p1.csv};
    \addlegendentry{$p_t=1$}

    \addplot+[
      thick,
      mark=x,
    ]
    table[
      col sep=comma,
      x expr={abs(\thisrow{mu}-\muVal)<1e-5 ? \thisrow{nt} : nan},
      y expr={abs(\thisrow{mu}-\muVal)<1e-5 ? \thisrow{Gbpx} : nan},
    ]{\csvfolder/results_p2.csv};
    \addlegendentry{$p_t=2$}

    \addplot+[
      thick,
      mark=o,
    ]
    table[
      col sep=comma,
      x expr={abs(\thisrow{mu}-\muVal)<1e-5 ? \thisrow{nt} : nan},
      y expr={abs(\thisrow{mu}-\muVal)<1e-5 ? \thisrow{Gbpx} : nan},
    ]{\csvfolder/results_p3.csv};
    \addlegendentry{$p_t=3$}

    \addplot+[
      thick,
      mark=star,
    ]
    table[
      col sep=comma,
      x expr={abs(\thisrow{mu}-\muVal)<1e-5 ? \thisrow{nt} : nan},
      y expr={abs(\thisrow{mu}-\muVal)<1e-5 ? \thisrow{Gbpx} : nan},
    ]{\csvfolder/results_p5.csv};
    \addlegendentry{$p_t=5$}

  \end{groupplot}

\end{tikzpicture}
  \caption{Convergence data for the regular solution $u_r$ and $\mu=\muVal$, where $G_{sym}$ and $G_{bpx}$ denote the GMRES iteration counts with symmetric norm preconditioning and left BPX preconditioning, respectively.}
  \label{fig:ode_convergence_plots:ur:mu10}
\end{figure}
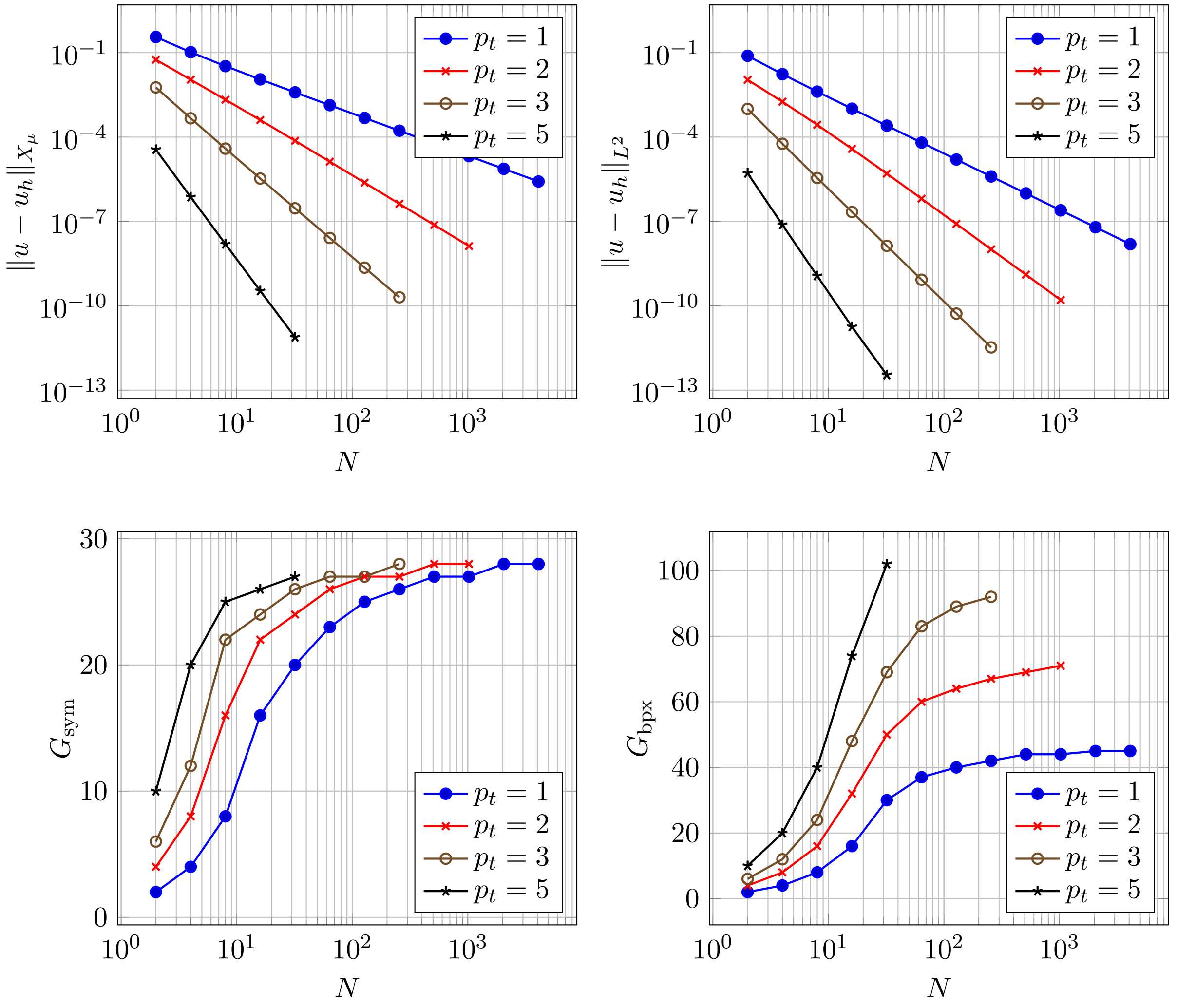

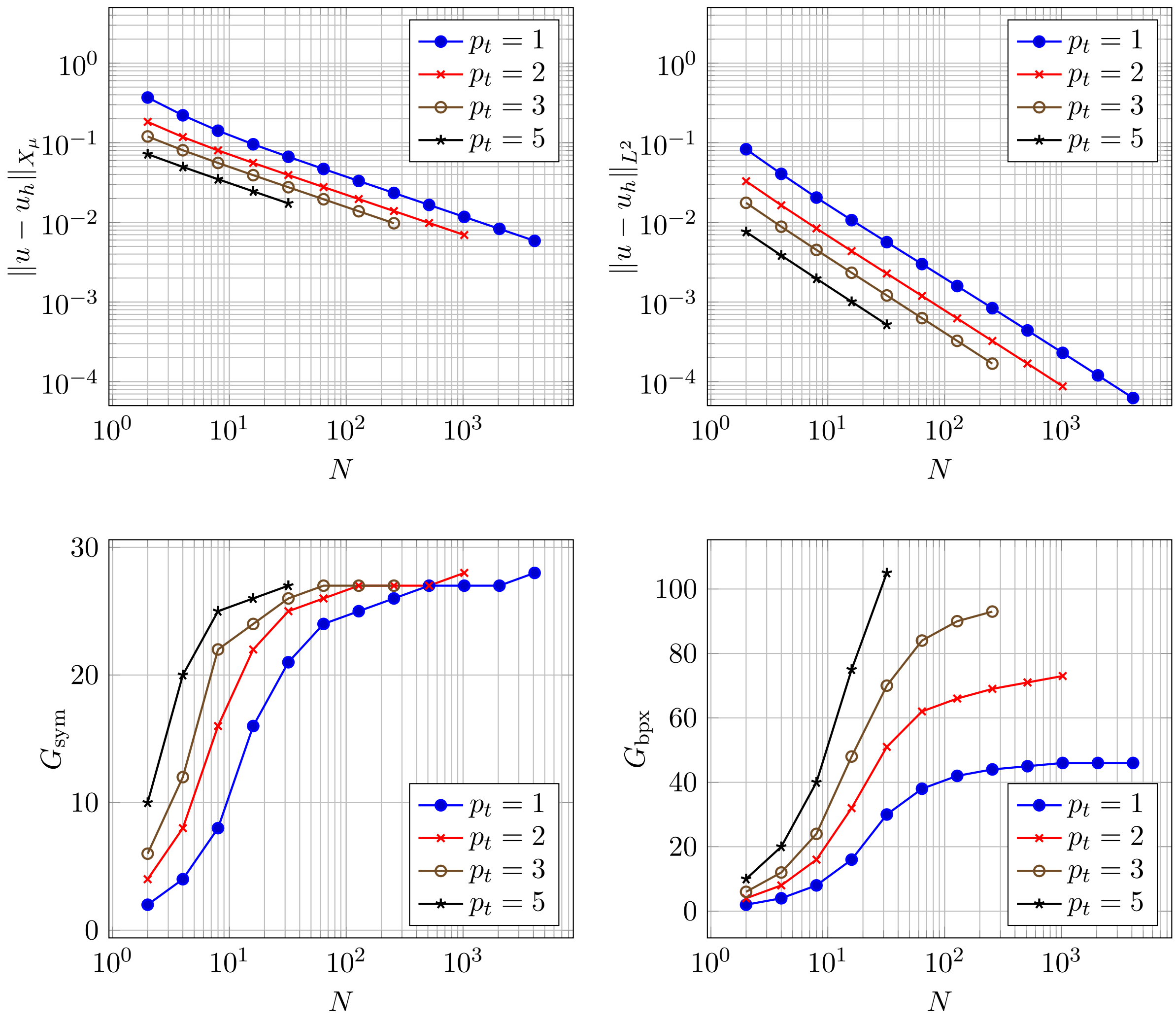
\begin{figure}[h!]
  \centering
  \def\csvfolder{simulations/ode/big_convergence_study_ode_hybrid_singular/}
  \def\muVal{10}
  \def\yminVal{5e-5}
  \begin{tikzpicture}
  \small
  \begin{groupplot}[
    group style={
      group size=2 by 2,       
      horizontal sep=1.8cm,    
      vertical sep=1.8cm,      
    },
    width=0.385\linewidth,      
    height=0.33\linewidth,
    scale only axis,
    grid=both,
  ]

    \nextgroupplot[
      xlabel={$N$},
      ylabel={$\|u - u_h\|_{X_\mu}$},
      ylabel near ticks,
      ylabel shift = -7pt,
      ymax = 0.5*1e1,
      ymin = \yminVal,
      ymode=log, xmode=log,
      legend pos=north east,
    ]
    \addplot+[
      thick,
      mark=*,
    ]
    table[
      col sep=comma,
      x expr={abs(\thisrow{mu}-\muVal)<1e-5 ? \thisrow{nt} : nan},
      y expr={abs(\thisrow{mu}-\muVal)<1e-5 ? \thisrow{Xerr} : nan},
    ]{\csvfolder/results_p1.csv};
    \addlegendentry{$p_t=1$}

    \addplot+[
      thick,
      mark=x,
    ]
    table[
      col sep=comma,
      x expr={abs(\thisrow{mu}-\muVal)<1e-5 ? \thisrow{nt} : nan},
      y expr={abs(\thisrow{mu}-\muVal)<1e-5 ? \thisrow{Xerr} : nan},
    ]{\csvfolder/results_p2.csv};
    \addlegendentry{$p_t=2$}

    \addplot+[
      thick,
      mark=o,
    ]
    table[
      col sep=comma,
      x expr={abs(\thisrow{mu}-\muVal)<1e-5 ? \thisrow{nt} : nan},
      y expr={abs(\thisrow{mu}-\muVal)<1e-5 ? \thisrow{Xerr} : nan},
    ]{\csvfolder/results_p3.csv};
    \addlegendentry{$p_t=3$}

    \addplot+[
      thick,
      mark=star,
    ]
    table[
      col sep=comma,
      x expr={abs(\thisrow{mu}-\muVal)<1e-5 ? \thisrow{nt} : nan},
      y expr={abs(\thisrow{mu}-\muVal)<1e-5 ? \thisrow{Xerr} : nan},
    ]{\csvfolder/results_p5.csv};
    \addlegendentry{$p_t=5$}

    \nextgroupplot[
      xlabel={$N$},
      ylabel={$\|u - u_h\|_{L^2}$},
      ylabel near ticks,
      ylabel shift = -7pt,
      ymax = 0.5*1e1,
      ymin = \yminVal,
      ymode=log, xmode=log,
      legend pos=north east,
    ]
    \addplot+[
      thick,
      mark=*,
    ]
    table[
      col sep=comma,
      x expr={abs(\thisrow{mu}-\muVal)<1e-5 ? \thisrow{nt} : nan},
      y expr={abs(\thisrow{mu}-\muVal)<1e-5 ? \thisrow{L2err} : nan},
    ]{\csvfolder/results_p1.csv};
    \addlegendentry{$p_t=1$}

    \addplot+[
      thick,
      mark=x,
    ]
    table[
      col sep=comma,
      x expr={abs(\thisrow{mu}-\muVal)<1e-5 ? \thisrow{nt} : nan},
      y expr={abs(\thisrow{mu}-\muVal)<1e-5 ? \thisrow{L2err} : nan},
    ]{\csvfolder/results_p2.csv};
    \addlegendentry{$p_t=2$}

    \addplot+[
      thick,
      mark=o,
    ]
    table[
      col sep=comma,
      x expr={abs(\thisrow{mu}-\muVal)<1e-5 ? \thisrow{nt} : nan},
      y expr={abs(\thisrow{mu}-\muVal)<1e-5 ? \thisrow{L2err} : nan},
    ]{\csvfolder/results_p3.csv};
    \addlegendentry{$p_t=3$}

    \addplot+[
      thick,
      mark=star,
    ]
    table[
      col sep=comma,
      x expr={abs(\thisrow{mu}-\muVal)<1e-5 ? \thisrow{nt} : nan},
      y expr={abs(\thisrow{mu}-\muVal)<1e-5 ? \thisrow{L2err} : nan},
    ]{\csvfolder/results_p5.csv};
    \addlegendentry{$p_t=5$}

    \nextgroupplot[
      xlabel={$N$},
      ylabel={$G_{\text{sym}}$},
      ylabel near ticks,
      ylabel shift = -7pt,
      xmode=log,
      legend pos=south east,
    ]
    \addplot+[
      thick,
      mark=*,
    ]
    table[
      col sep=comma,
      x expr={abs(\thisrow{mu}-\muVal)<1e-5 ? \thisrow{nt} : nan},
      y expr={abs(\thisrow{mu}-\muVal)<1e-5 ? \thisrow{Gsymm} : nan},
    ]{\csvfolder/results_p1.csv};
    \addlegendentry{$p_t=1$}

    \addplot+[
      thick,
      mark=x,
    ]
    table[
      col sep=comma,
      x expr={abs(\thisrow{mu}-\muVal)<1e-5 ? \thisrow{nt} : nan},
      y expr={abs(\thisrow{mu}-\muVal)<1e-5 ? \thisrow{Gsymm} : nan},
    ]{\csvfolder/results_p2.csv};
    \addlegendentry{$p_t=2$}

    \addplot+[
      thick,
      mark=o,
    ]
    table[
      col sep=comma,
      x expr={abs(\thisrow{mu}-\muVal)<1e-5 ? \thisrow{nt} : nan},
      y expr={abs(\thisrow{mu}-\muVal)<1e-5 ? \thisrow{Gsymm} : nan},
    ]{\csvfolder/results_p3.csv};
    \addlegendentry{$p_t=3$}

    \addplot+[
      thick,
      mark=star,
    ]
    table[
      col sep=comma,
      x expr={abs(\thisrow{mu}-\muVal)<1e-5 ? \thisrow{nt} : nan},
      y expr={abs(\thisrow{mu}-\muVal)<1e-5 ? \thisrow{Gsymm} : nan},
    ]{\csvfolder/results_p5.csv};
    \addlegendentry{$p_t=5$}

    \nextgroupplot[
      xlabel={$N$},
      ylabel={$G_{\text{bpx}}$},
      ylabel near ticks,
      ylabel shift = -7pt,
      xmode=log,
      legend pos=south east,
    ]
    \addplot+[
      thick,
      mark=*,
    ]
    table[
      col sep=comma,
      x expr={abs(\thisrow{mu}-\muVal)<1e-5 ? \thisrow{nt} : nan},
      y expr={abs(\thisrow{mu}-\muVal)<1e-5 ? \thisrow{Gbpx} : nan},
    ]{\csvfolder/results_p1.csv};
    \addlegendentry{$p_t=1$}

    \addplot+[
      thick,
      mark=x,
    ]
    table[
      col sep=comma,
      x expr={abs(\thisrow{mu}-\muVal)<1e-5 ? \thisrow{nt} : nan},
      y expr={abs(\thisrow{mu}-\muVal)<1e-5 ? \thisrow{Gbpx} : nan},
    ]{\csvfolder/results_p2.csv};
    \addlegendentry{$p_t=2$}

    \addplot+[
      thick,
      mark=o,
    ]
    table[
      col sep=comma,
      x expr={abs(\thisrow{mu}-\muVal)<1e-5 ? \thisrow{nt} : nan},
      y expr={abs(\thisrow{mu}-\muVal)<1e-5 ? \thisrow{Gbpx} : nan},
    ]{\csvfolder/results_p3.csv};
    \addlegendentry{$p_t=3$}

    \addplot+[
      thick,
      mark=star,
    ]
    table[
      col sep=comma,
      x expr={abs(\thisrow{mu}-\muVal)<1e-5 ? \thisrow{nt} : nan},
      y expr={abs(\thisrow{mu}-\muVal)<1e-5 ? \thisrow{Gbpx} : nan},
    ]{\csvfolder/results_p5.csv};
    \addlegendentry{$p_t=5$}

  \end{groupplot}

\end{tikzpicture}
  \caption{Convergence data for the singular solution $u_s$ and $\mu=\muVal$, where $G_{sym}$ and $G_{bpx}$ denote the GMRES iteration counts with symmetric norm preconditioning and left BPX preconditioning, respectively.}
  \label{fig:ode_convergence_plots:us:mu10}
\end{figure}

As solver, we apply unrestarted GMRES from a zero initial guess, using either left BPX preconditioning via \eqref{eq:ode_bpx_inverse} or symmetric preconditioning with the norm inducing matrix $\mathbf{D}_\mu$. We require a relative Euclidean residual of at most $10^{-12}$, measured in the original system for BPX and in the symmetrically transformed system for the comparison. We measure the errors in the $L^2(0,T)$- and $X_\mu$-norms using the recovered solution from the symmetric case. We show the results for the smooth solution $u_r$ in Figure \ref{fig:ode_convergence_plots:ur:mu10} and for the singular solution $u_s$ in Figure \ref{fig:ode_convergence_plots:us:mu10}, both for $\mu = 10$.
For the smooth solution and $\mu=10$, the measured energy orders on the finest two grids are $1.50$, $2.50$, $3.50$, and $5.47$ for $p=1,2,3,5$, respectively. For $u_s(t)=\sqrt{t}$ they are close to $0.50$ for all four degrees.

\begin{remark}
  The iteration counts for the BPX preconditioner are quite high, as one would expect for this type of preconditioner. However, the BPX preconditioner is very easy to implement and applies in a matrix free fashion. Furthermore, more sophisticated preconditioners such as multigrid \cite{TrottenbergOosterleeSchullerMultigrid2001} or domain decomposition \cite{TarekDomainDecompositionMethods2008} are possible. Early results for the multigrid preconditioner show its feasibility and efficiency, but a detailed study is beyond the scope of this work and will be subject of future research.
\end{remark}


\section{Coercive Formulation of the Heat Equation and its Discretization}
In this section we consider the heat equation in a bounded (Lipschitz) domain $\Omega \subset \R^d$, $d=1,2,3$, with homogeneous Dirichlet boundary conditions on a finite time interval $(0,T)$, i.e.
\begin{align}\label{eq:heat equation}
    \partial_t u(t,x) -\Delta_x u(t,x) &= f(t,x), \quad (t,x) \in Q = (0,T) \times \Omega, \nonumber \\
    u(t,x) &= 0, \quad \quad  \quad \, (t,x) \in \Sigma = (0,T) \times \partial\Omega, \\
    u(0,x) &= 0, \quad \quad \quad \, \, x \in \Omega.\nonumber
\end{align}
as a model problem for a linear parabolic PDE.
In order to apply the operator theoretic framework developed in Section 4, we will transform the PDE into a sequence of ODEs using the spectral decomposition of the spatial Laplacian. Therefore consider the eigenfunctions of the spatial Dirichlet-Laplacian, given as the solutions of 
\begin{align}
    -\Delta_x \phi_i = \mu_i \phi_i \quad\text{on } \Omega, \quad \phi_i = 0 \quad \text{on }\partial\Omega,\quad i\in \mathbb{N}. 
\end{align}
Additionally requiring $\norm{\phi_i}_{L^2(\Omega)}=1$ they build an orthonormal basis of $L^2(\Omega)$, an orthogonal basis of $H^1_0(\Omega)$ and the eigenvalues $\{\mu_i\}_{i\in \mathbb{N}}$ are positive and monotonically increasing. Hence, any function $u\in L^2(Q)$ can be represented as 
\begin{align}
    u(t,x) = \sum_{i=0}^{\infty} u_i(t) \phi_i(x), \label{eq:uInEigenfunctions}
\end{align}
with time dependent coefficients $u_i\in L^2(0,T).$

\begin{lemma}
    The operator 
    \begin{equation*}
        \mathcal{B}:\dom(\mathcal{B})\subset L^2(Q)\to L^2(Q),\quad \mathcal{B}u = \partial_t u-\Delta_x u   
    \end{equation*}
    with domain 
    \begin{equation*}
        \dom(\mathcal{B}) = H^1_{0,}(0,T;L^2(\Omega))\cap L^2(0,T;H^1_{\Delta_x}(\Omega)),
    \end{equation*}
    is densely defined, closed and $\ker(\mathcal{B}) = \ker(\mathcal{B}^\ast) = \{0\}.$ Moreover, it satisfies the inf-sup condition \eqref{eq:inf-sup B}, i.e., there exists $\beta>0$ such that 
    \begin{equation*}
        \beta \leq \inf_{0\neq u\in \dom(\mathcal{B})} \sup_{0\neq v\in L^2(Q)} \frac{(\mathcal{B}u,v)_{L^2(Q)}}{\norm{u}_{L^2(Q)}\norm{v}_{L^2(Q)}}. 
    \end{equation*}
\end{lemma}
\begin{proof}
    It is easy to see that the operator is densely defined. 
    For any $u\in \dom(\mathcal{B})$ we have that 
    \begin{align*}
        \norm{\mathcal{B}u}_{L^2(Q)}^2
        & =
        \norm{\partial_t u - \Delta_x u}_{L^2(Q)}^2\\
        & = \norm{\partial_t u}_{L^2(Q)}^2 + 2(\partial_t u,- \Delta_x u)_{L^2(Q)} + \norm{\Delta_x u}_{L^2(Q)}^2. 
    \end{align*}
    Due to the Poincar\'{e} inequality \eqref{eq:Poincare 1D} we have 
    \begin{equation*}
        \norm{\partial_t u}_{L^2(Q)}^2
        \geq \left(1+\frac{4\pi^2}{T^2}\right)^{-1}\norm{u}_{H^1_{0,}(0,T;L^2(\Omega))}^2  
    \end{equation*}
    holds true. Moreover, $\norm{\Delta_x u}_{L^2(Q)}$ defines a norm on $L^2(0,T;H^1_{\Delta_x}(\Omega))$, see \cite[Eq. (3.1)]{FeuerleLoescherSteinbachUrban2026}. For the remaining term, consider the representation \eqref{eq:uInEigenfunctions} to compute 
    \begin{align*}
        (\partial_t u,-\Delta_x u)_{L^2(Q)}
        & = 
        \left(\sum_{i=1}^\infty\partial_t u_i\phi_i,\sum_{j=1}^\infty u_j(-\Delta_x\phi_j)\right)_{L^2(Q)}
        \\
        & = \sum_{i,j=1}^\infty \mu_j\big(\partial_t u_i,u_j\big)_{L^2(0,T)}
        \big(\phi_i,\phi_j\big)_{L^2(\Omega)}
        \\
        &=
        \sum_{i=1}^\infty \mu_i(\partial_t u_i,u_i)_{L^2(0,T)}
        =
        \sum_{i=1}^\infty \mu_i \frac{|u_i(T)|^2}{2} 
        \geq 0.
    \end{align*} 
    So altogether 
    \begin{equation*}
        \norm{\mathcal{B} u}_{L^2(Q)}^2 
        \geq 
        \min\left\{1,\left(1+\frac{4T^2}{\pi^2}\right)^{-1}\right\} \big(\norm{u}_{H^1_{0,}(0,T;L^2(\Omega))}^2 + \norm{u}_{L^2(0,T;H^1_\Delta(\Omega))}^2\big),
    \end{equation*}
    which shows closedness and $\ker(\mathcal{B}) =\{0\}$. 
    By self-duality of $L^2(Q)$ we have that 
    \begin{equation*}
        \norm{\mathcal{B} u}_{L^2(Q)} 
        = 
        \sup_{0\neq v\in L^2(Q)} \frac{(\mathcal{B}u,v)_{L^2(Q)}}{\norm{v}_{L^2(Q)}}, 
    \end{equation*}
    implying the inf-sup condition. Moreover, the adjoint is given as 
    $\mathcal{B}^\ast u = -\partial_t u - \Delta_x u$ with domain 
    \begin{equation*}
        \dom(\mathcal{B}^\ast) = H^1_{,0}(0,T;L^2(\Omega))\cap L^2(0,T;H^1_{\Delta_x}(\Omega)),
    \end{equation*}
    and an analogous computation shows that $\ker(\mathcal{B}^\ast) = \{0\}$. 
\end{proof}
So all the assumption of Theorem \ref{thm:polar_decomposition} are satisfied, which immediately gives the energy space for the heat equation 
\begin{equation*}
    \mathcal{X} 
    = \dom(|\mathcal{B}|^{1/2}) 
    = [L^2(Q), \dom(\mathcal{B})]_{1/2}
    = H^{1/2}_{0,}(0,T;L^2(\Omega)) \cap L^2(0,T;H^1_0(\Omega)).   
\end{equation*} 
\begin{remark}
    Note that the energy space $\mathcal{X}$ is exactly the Lions-Magenes space \cite[§5]{LionsMagenesVol1} derived for the heat equation, that was also used in \cite{ZankSteinbachCoerciveSpaceTimeFinite2020}. 

    As in the previous section, we can compute the unitary factor of the polar decomposition to be 
    \begin{align}
        \mathcal{H}u 
        = \mathcal{H}\left(\sum_{i=1}^\infty u_i \phi_i\right)
        := \sum_{i=1}^\infty \big(\mathcal{H}_{\mu_i} u_i\big) \phi_i,
    \end{align}
    where $\mathcal{H}_{\mu_i}$ is the parabolic Hilbert transformation from Definition \ref{def:ParabolicHilbertTransform}. 
    However, as mentioned already in the last section, this transformation is cumbersome to realize in practice. Hence, for a practical realization we will replace it by an approximation. 
\end{remark}

\subsection{Discretization, Preconditioning and Numerical Examples}\label{sec:pde_numerical_examples}
As in the case of the ODE, we will again approximate the transformation $\mathcal{H}$ by $\widehat{\mathcal{H}}:= \mathcal{H}_T + I$ to obtain a practical discretization for the heat equation. Let again $0 = t_0 < t_1 < \dots < t_{N_t} = T$ be a quasi-uniform partition with mesh size $h_t$ and $\Omega_h$ be a conforming and quasi-uniform triangulation of the spatial domain $\Omega$ with mesh size $h_x$. We consider the tensor-product trial and test space
\begin{align}
    \mathcal{X}_h = X_{h_t} \otimes V_{h_x}, \label{eq:tensorSpace}
\end{align}
where
\begin{align}
    X_{h_t} = S_{h_t}^{1}(0,T) \cap H_{0,}^{1}(0,T), \quad
    V_{h_x} = S_{h_x}^{1}(\Omega_h) \cap H_0^1(\Omega),
\end{align}
and consider the discrete variational formulation to find $u_h \in \mathcal{X}_h$ such that
\begin{align}\label{eq:perturbed DVF heat}
    \lvert \widehat{b} \rvert(u_h, v_h) = \langle f, \widehat{\mathcal{H}} v_h \rangle_Q, \quad \forall v_h \in \mathcal{X}_h, 
\end{align}
with the perturbed bilinear form
\begin{align}
    \lvert \widehat{b} \rvert(u_h, v_h) = ( \partial_t u_h, \widehat{\mathcal{H}} v_h )_{L^2(Q)} + \langle \nabla_x u_h, \widehat{\mathcal{H}} \nabla_x  v_h \rangle_{L^2(Q)}.
\end{align}
Analogous to Lemma \ref{lem:coercivity of hat bmu} and Lemma \ref{lem:discrete boundedness bilinear form}, the bilinear form $|\widehat{b}|(\cdot,\cdot)$ admits the following properties. 
\begin{lemma}
    The perturbed bilinear form $\lvert \widehat{b} \rvert(\cdot, \cdot)$ elliptic
    \begin{align}\label{eq:perturbed bilinearform heat coercivity}
        \lvert \widehat{b} \rvert(u_h, u_h) \geq \| u_h \|_{\mathcal{X}}^2, \quad \forall u_h \in \mathcal{X}_h,
    \end{align}
    and bounded up to a logarithmic factor, i.e.,
    \begin{align}\label{eq:perturbed bilinearform heat boundedness}
        \lvert \widehat{b} \rvert(u_h, v_h) \leq C \left(1 + \lvert \log(h_t) \rvert\right) \| u_h \|_{\mathcal{X}} \| v_h \|_{\mathcal{X}}, \quad \forall u_h, v_h \in \mathcal{X}_h,
    \end{align}
    for some constant $C > 0$ independent of $h_t$ and $h_x$.
\end{lemma}
\begin{proof}
    Consider the discrete eigenvalue problem 
    \begin{equation*}
        (\nabla_x \phi_{h_x,i},\nabla_x v_{h_x})_{L^2(\Omega)}
        =
        \widetilde{\mu}_{i}(\phi_{h_x,i},v_{h_x})_{L^2(\Omega)},\quad \forall v_{h_x}\in V_{h_x}. 
    \end{equation*}
    The normalized eigenfunctions $\phi_{h_x,i}$, $i=1,\ldots,\text{dim}(V_{h_x}),$ built an $L^2$-orthogonal basis of $V_{h_x}$. Thus each $u_{h}\in \mathcal{X}_h$ admits the representation 
    \begin{equation*}
        u_{h}(x,t) = \sum_{i=1}^{\text{dim}(V_{h_x})} \widetilde{u}_i(t) \phi_{h_x,i}(x),\quad \widetilde{u}_i(t) = \int_\Omega u_h(x,t)\phi_{h_x,i}(x)\, dx. 
    \end{equation*}
    Therefore, the bilinear form becomes 
    \begin{align*}
    |\widehat{b}|(u_h,v_h)
    &=
    (\partial_t u_h,\widehat{\mathcal{H}}v_h)_{L^2(Q)}
    +(\nabla_x u_h,\widehat{\mathcal{H}}\nabla_x v_h)_{L^2(Q)}
    \\
    &=
    \sum_{i,j=1}^{\text{dim}(V_{h_x})}
    \Big[
        (\partial_t\widetilde{u}_i,
         \widehat{\mathcal{H}}\widetilde{v}_j)_{L^2(0,T)}
        (\phi_{h_x,i},\phi_{h_x,j})_{L^2(\Omega)}
    \\
    &\hspace{5em}
        +(\widetilde{u}_i,
          \widehat{\mathcal{H}}\widetilde{v}_j)_{L^2(0,T)}
        (\nabla_x\phi_{h_x,i},
         \nabla_x\phi_{h_x,j})_{L^2(\Omega)}
    \Big]
    \\
    &=
    \sum_{i=1}^{\text{dim}(V_{h_x})}
    \Big[
        (\partial_t\widetilde{u}_i,
         \widehat{\mathcal{H}}\widetilde{v}_i)_{L^2(0,T)}
        +\widetilde{\mu}_{i}
        (\widetilde{u}_i,
         \widehat{\mathcal{H}}\widetilde{v}_i)_{L^2(0,T)}
    \Big]
    \\
    &=
    \sum_{i=1}^{\text{dim}(V_{h_x})}
    |\widehat{b}_{\widetilde{\mu}_{i}}|
    (\widetilde{u}_i,\widetilde{v}_i).
\end{align*}
    Thus, the results from Lemmata \ref{lem:coercivity of hat bmu} and \ref{lem:discrete boundedness bilinear form} directly transfer to $|\widehat{b}(\cdot,\cdot)|$. 
\end{proof}

The main result is the following. 
\begin{theorem}\label{thm:main result heat}
    For any $f\in [H^{1-s}_{,0}(0,T;L^2(\Omega))\cap L^2(0,T;H^1_0(\Omega))]'$, $s>1/2.$ the discrete variational formulation \eqref{eq:perturbed DVF heat} admits a unique solution $u_h\in \mathcal{X}_h$. Moreover, if the unique solution of 
    \begin{equation}\label{eq:heat equatoin distributional}
        \partial_t u - \Delta_x u = f\quad \text{in } [H^{1-s}_{,0}(0,T;L^2(\Omega))\cap L^2(0,T;H^1_0(\Omega))]'
    \end{equation} 
    satisfies $u\in \mathcal{X}\cap H^{s}(0,T;L^2(\Omega))\cap L^2(0,T;H^r(\Omega))$, for $s\in (1/2,2],$ and $r\in [1,2]$, there holds the error estimate 
    \begin{align*}
    \|u-u_h\|_{\mathcal Q}
    \leq C\left(
        h_x^\varrho+
        h_t^{s-1/2}\sqrt{1+\frac{h_t}{h_x^2}}
    \right)
    \left(
        \|u\|_{H^s(0,T;L^2(\Omega))}+
        \|u\|_{L^2(0,T;H^r(\Omega))}
    \right),
\end{align*}
where $\varrho=\min\left\{r-1,\;r\left(1-\frac{1}{2s}\right)\right\}$ and $C$ is independent of $h_t$ and $h_x$.
\end{theorem}
\begin{proof}
Since $H^{1-s}_{0,}(0,T) = H^{1-s}(0,T)$, for $s>1/2$, see \cite[Theorem 3.40]{McLeanStronglyEllipticSystems2000}, we have that  
\begin{equation*}
    Y_s = H^{1-s}_{,0}(0,T;L^2(\Omega))\cap L^2(0,T;H^1_0(\Omega)) = H^{1-s}(0,T;L^2(\Omega))\cap L^2(0,T;H^1_0(\Omega)). 
\end{equation*}
Furthermore, 
$\widehat{\mathcal{H}}(\mathcal{X}_h) \subset H^1(Q) \subset Y_s$ and therefore the functional 
\begin{equation*}
    F_h(v_h) := \langle f,\widehat{\mathcal{H}}v_h\rangle_{Y_s'\times Y_s}
\end{equation*}
is well-defined and continuous. 
For fixed $h>0$ the established discrete boundedness \eqref{eq:perturbed bilinearform heat boundedness} and coercivity \eqref{eq:perturbed bilinearform heat coercivity} imply
unique solvability of the discrete variational formulation \eqref{eq:perturbed DVF heat} by the lemma of Lax--Milgram.

To obtain the error estimate let $I_{h_t}^1: C([0,T])\to X_{h_t}$ denote the nodal interpolation operator and let $Q_{h_x}^1: L^2(\Omega) \to V_{h_x}$ denote the $L^2$-projection. Since $H^{s}_{0,}(0,T;L^2(\Omega))\hookrightarrow C([0,T];L^2(\Omega))$ the function
\begin{equation*}
    I_{h_t}^1 Q_{h_x}^1 u = Q_{h_x}^1 I_{h_t}^1 u \in \mathcal{X}_h, 
\end{equation*}
is well defined, see \cite[Lemma 2.8.4]{ZankDiss}.  
This permits to split the error 
\begin{equation*}
    \norm{u-u_h}_{\mathcal{X}}
    \leq 
    \norm{\eta}_{\mathcal{X}} + \norm{Q_{h_x}^1 I_{h_t}^1 u-u_h}_{\mathcal{X}},
\end{equation*}
where $\eta := u-Q_{h_x}^1 I_{h_t}^1 u$. 
Noting that $\eta(0)=0$ and $\eta(T)= (I-Q_{h_x}^1)u(T)$, we obtain by integration by parts 
\begin{align*}
    \skpr{\partial_t \eta,v_h}_{Y_s'\times Y_s} 
    &= (\eta(T),v_h(T))_{L^2(\Omega)} - (\eta,\partial_t v_h)_{L^2(Q)}.   
\end{align*}
Due to the orthogonality of the $L^2$-projection the term 
\begin{align*}
    (\eta(T),v_h(T))_{L^2(\Omega)} 
    &= ((I-Q_{h_x}^1)u(T),v_h(T))_{L^2(\Omega)}\\
    &= (u(T),(I-Q_{h_x}^1) v_h(T))_{L^2(\Omega)}=0,
\end{align*}
vanishes.  
Then, by the Cauchy--Schwarz inequality and the inverse inequality \eqref{eq:inverse inequality in time h1 to h12}, we get 
\begin{align*}
    \skpr{\partial_t \eta,v_h}_{Y_s'\times Y_s}
    &= -(\eta,\partial_t v_h)_{L^2(Q)} = (Q_{h_x}^1 \eta,\partial_t v_h)_{L^2(Q)}\\
    &\leq \norm{Q_{h_x}^1\eta}_{L^2(Q)}\norm{\partial_t v_h}_{L^2(Q)}\\
    &\leq C h_t^{-1/2} \norm{Q_{h_x}^1\eta}_{L^2(Q)}\norm{v_h}_{H^{1/2}_{0,}(0,T;L^2(\Omega))}.  
\end{align*}
For any $v_h\in \mathcal{X}_h$ recall that 
\begin{align*}
    |\widehat{b}|(\eta,v_h)
    & = \skpr{\partial_t \eta, \widehat{\mathcal{H}} v_h}_{Y_s'\times Y_s} + (\nabla_x \eta, \widehat{\mathcal{H}}\nabla_x v_h)_{L^2(Q)}\\
    & = \skpr{\partial_t \eta, v_h}_{Y_s'\times Y_s} + \skpr{\partial_t \eta, \Hmod v_h}_{Y_s'\times Y_s} \\
    &\quad + (\nabla_x \eta, \nabla_x v_h)_{L^2(Q)} + (\nabla_x \eta, \Hmod \nabla_x v_h)_{L^2(Q)}.  
\end{align*}
Using the established bound, as well as the properties of the modified Hilbert transformation \cite[Eq. (2.9)]{ZankSteinbachCoerciveSpaceTimeFinite2020}, we get 
\begin{equation}\label{eq:auxiliary bound bilinear form}
    \begin{aligned}
    |\widehat{b}|(\eta,v_h)
    & \leq Ch_t^{-1/2} \norm{Q_{h_x}^1\eta}_{L^2(Q)}\norm{v_h}_{H^{1/2}_{0,}(0,T;L^2(\Omega))} \\
    &\hskip 3em + \norm{\eta}_{H^{1/2}_{0,}(0,T;L^2(\Omega))}\norm{v_h}_{H^{1/2}_{0,}(0,T;L^2(\Omega))} \\
    & \hskip 5em  + 2\norm{\eta}_{L^2(0,T;H^1_0(\Omega))}\norm{v_h}_{L^2(0,T;H^1_0(\Omega))}\\
    & \leq C\big(\norm{\eta}_{\mathcal{X}} + h_t^{-1/2}\norm{Q_{h_x}^1\eta}_{L^2(Q)}\big)\norm{v_h}_{\mathcal{X}}. 
\end{aligned}    
\end{equation}
As $\widehat{\mathcal{H}}(\mathcal{X}) \subset Y_s$, testing \eqref{eq:heat equatoin distributional} with $v_h\in \mathcal{X}_h\subset Y_s$ is valid and immediately gives the Galerkin orthogonality
\[
    |\widehat{b}|(u-u_h,v_h)=0
    \qquad\forall v_h\in\mathcal Q_h.
\]
Together with the coercivity \eqref{eq:perturbed bilinearform heat coercivity} this gives 
\begin{equation}\label{eq:aux coercive and GO}
    \begin{aligned}
        \norm{Q_{h_x}^1 I_{h_t}^1 u-u_h}_{\mathcal{X}}^2 
        & \leq |\widehat{b}|(Q_{h_x}^1 I_{h_t}^1 u-u_h,Q_{h_x}^1 I_{h_t}^1 u-u_h)\\
        &= |\widehat{b}|(Q_{h_x}^1 I_{h_t}^1 u-u,Q_{h_x}^1 I_{h_t}^1 u-u_h)\\
        &= -|\widehat{b}|(\eta,Q_{h_x}^1 I_{h_t}^1 u-u_h),
    \end{aligned}    
\end{equation}
and since $Q_{h_x}^1 I_{h_t}^1 u-u_h\in \mathcal{X}_h$ the bound \eqref{eq:auxiliary bound bilinear form} implies 
\begin{align*}
    \norm{Q_{h_x}^1 I_{h_t}^1 u-u_h}_{\mathcal{X}}
     \leq C\big(\norm{\eta}_{\mathcal{X}} + h_t^{-1/2}\norm{Q_{h_x}^1\eta}_{L^2(Q)}\big). 
\end{align*}
Altogther, we obtain 
\begin{equation*}
    \norm{u-u_h}_{\mathcal{X}}
    \leq \norm{\eta}_{\mathcal{X}} + \norm{Q_{h_x}^1 I_{h_t}^1 u-u_h}_{\mathcal{X}} 
    \leq C\big(\norm{\eta}_{\mathcal{X}} + h_t^{-1/2}\norm{Q_{h_x}^1\eta}_{L^2(Q)}\big).  
\end{equation*}
It remains to bound the appearing terms. Firstly, we have that 
\begin{equation*}
    \norm{Q_{h_x}^1\eta}_{L^2(Q)} 
    = \norm{Q_{h_x}^1(u-I_{h_t}^1 u)}_{L^2(Q)}
    \leq \norm{u-I_{h_t}^1 u}_{L^2(Q)}
    \leq C h_t^s \norm{u}_{H^s(0,T;L^2(\Omega))}.  
\end{equation*}
For the second term we have that 
\begin{align*}
    \norm{\eta}_{\mathcal{X}}
    &\leq 
    \norm{u-Q_{h_x}^1 u}_{H^{1/2}_{0,}(0,T;L^2(\Omega))} +\norm{u-Q_{h_x}^1 u}_{L^2(0,T;H^1_0(\Omega))}\\
    &\quad +  \norm{Q_{h_x}^1 (I-I_{h_t}^1)u}_{H^{1/2}_{0,}(0,T;L^2(\Omega))} + \norm{Q_{h_x}^1 (I-I_{h_t}^1)u}_{L^2(0,T;H^1_0(\Omega))}.  
\end{align*}
By standard best approximation error estimates we get  
\begin{align*}
    \norm{u-Q_{h_x}^1u}_{L^2(0,T;L^2(\Omega))} &\leq ch_x^{r} \norm{u}_{L^2(0,T;H^r(\Omega))},\\
    \norm{u-Q_{h_x}^1u}_{L^2(0,T;H^1_0(\Omega))} &\leq ch_x^{r-1} \norm{u}_{L^2(0,T;H^r(\Omega))}.  
\end{align*}
For $\sigma = \min\{s,1\}>1/2$ we have the stability estimate 
\begin{equation*}
    \norm{u-Q_{h_x}^1u}_{H^\sigma_{0,}(0,T;L^2(\Omega))} 
    \leq 
    c\norm{u}_{H^s(0,T;L^2(\Omega))}. 
\end{equation*}
Interpolating, see \cite[Lemma B.1]{McLeanStronglyEllipticSystems2000}, gives for 
$\theta=1/(2\sigma)$
that 
\begin{align*}
    \norm{u-Q_{h_x}^1 u}_{H^{1/2}_{0,}(0,T;L^2(\Omega))}
    & \leq c \norm{u-Q_{h_x}u}_{L^2(0,T;L^2(\Omega))}^{1-\theta} \norm{u-Q_{h_x}u}_{H^\sigma_{0,}(0,T;L^2(\Omega))}^\theta\\
    &\leq ch_x^{r(1-\theta)}\norm{u}_{L^2(0,T;H^r(\Omega))}^{1-\theta}\norm{u}_{H^s(0,T;L^2(\Omega))}^\theta
\end{align*}
and using Young's inequality this becomes 
\begin{equation*}
    \norm{u-Q_{h_x}^1 u}_{H^{1/2}_{0,}(0,T;L^2(\Omega))}
    \leq 
    c h_x^{r(1-\theta)} \big(\norm{u}_{L^2(0,T;H^r(\Omega))} + \norm{u}_{H^s(0,T;L^2(\Omega))}\big). 
\end{equation*}
Moreover, due to the stability of $L^2$-projection, we get  
\begin{align*}
    \norm{Q_{h_x}^1 (u-I_{h_t}^1u)}_{H^{1/2}_{0,}(0,T;L^2(\Omega))}
    &\leq 
    \norm{u-I_{h_t}^1u}_{H^{1/2}_{0,}(0,T;L^2(\Omega))}\\ 
    &\leq 
    ch_t^{s-1/2}\norm{u}_{H^s(0,T;L^2(\Omega))}. 
\end{align*}
For the last term, applying the inverse inequality \cite[Cor. 1.141]{ErnGuermondTheoryPracticeFinite2004} gives 
\begin{align}\label{eq:aux inverse inequality for error estimate}
    \norm{Q_{h_x}^1(u-I_{h_t}^1 u)}_{L^2(0,T;H^1_0(\Omega))} 
    &\leq
    ch_x^{-1}\norm{Q_{h_x}^1(u-I_{h_t}^1u)}_{L^2(Q)}\notag\\
    & \leq ch_x^{-1}\norm{u-I_{h_t}^1 u}_{L^2(Q)}\\
    &\leq ch_x^{-1}h_t^{s}\norm{u}_{H^s(0,T;L^2(\Omega))}.\notag 
\end{align} 
Altogether, we then have  
\begin{align*}
    \norm{\eta}_{\mathcal{X}}
    &=
    c \big[
       (h_x^{r(1-\theta)} +h_x^{r-1})\norm{u}_{L^2(0,T;H^r(\Omega))} 
       \\
       &\hskip 3em  + (h_x^{r(1-\theta)} + h_t^{s-1/2} + h_x^{-1}h_t^s)\norm{u}_{H^s(0,T;L^2(\Omega))} 
    \big]\\
    & = 
    c \Big[
       (h_x^{r(1-\theta)} +h_x^{r-1})\norm{u}_{L^2(0,T;H^r(\Omega))} 
       \\
       &\hskip 3em  + \left(h_x^{r(1-\theta)} + h_t^{s-1/2}\Big(1+\tfrac{h_t^{1/2}}{h_x}\Big)\right)\norm{u}_{H^s(0,T;L^2(\Omega))} 
    \Big] \\
    & \leq  
    c\Big[
        h_x^\varrho + h_t^{s-1/2}\sqrt{1+\tfrac{h_t}{h_x^2}}      
    \Big]\big(
        \norm{u}_{H^s(0,T;L^2(\Omega))} + \norm{u}_{L^2(0,T;H^r(\Omega))}
    \big), 
\end{align*}
where $\varrho = \min\{r-1,r\left(1-\frac{1}{2s}\right)\}. $
\end{proof}
\begin{corollary}\label{cor:regular solution estimate}
    If the solution is sufficiently smooth, i.e., $u\in \mathcal{X}$ satisfies 
    \begin{equation*}
        u\in H^{\sigma}(0,T;H^1_0(\Omega))\cap L^2(0,T;H^2(\Omega)),\quad \sigma\in(1/2,2],  
    \end{equation*}
    we can replace the inverse inequality in \eqref{eq:aux inverse inequality for error estimate} by the $H^1$-stability of the $L^2$-projection \cite{BramblePasciakSteinbach2002,DieningStornTscherpel2021} to obtain 
    \begin{align*}
        \norm{Q_{h_x}^1 (u-I_{h_t}^1u)}_{L^2(0,T;H^1_0(\Omega))}
        & \leq c\norm{u-I_{h_t}^1u}_{L^2(0,T;H^1_0(\Omega))}\\
        & \leq ch_t^{\sigma}\norm{u}_{H^{s}(0,T;H^1_0(\Omega))}. 
    \end{align*}
    The overall estimate then reads 
    \begin{equation*}
        \norm{u-u_h}_{\mathcal{X}}
        \leq 
        C(u)(h_x + h_t^{\sigma-1/2}), 
    \end{equation*}
    where $C(u)$ only depends on the regularity of $u$. 
\end{corollary}
\begin{remark}
    The error estimate in Theorem~\ref{thm:main result heat}
is balanced by the choice $h_t\simeq h_x^{r/s}$.
In particular, for solutions with only the parabolic regularity
\[
    u\in H^1(0,T;L^2(\Omega))
    \cap L^2(0,T;H^2(\Omega)),
\]
balancing the two error contributions requires the usual
parabolic scaling $h_t\simeq h_x^2$.

For the smoother solutions considered in
Corollary~\ref{cor:regular solution estimate},
Table~\ref{tab:comparisonScaling} compares the convergence rates
guaranteed by the two estimates for different relations between
the temporal and spatial mesh sizes. In particular, if $\sigma\geq 1$ isotropic scaling is sufficient to obtain optimal orders of convergence. 

\begin{table}[H]
    \centering
    \caption{Convergence rates w.r.t. $h_x$ for
        $u\in H^s(0,T;H^1_0(\Omega))
        \cap L^2(0,T;H^2(\Omega))$, $s\in(1/2,2]$.}
    \label{tab:comparisonScaling}
    \begin{tabular}{c|c|c}
        Scaling relation
        & Theorem~\ref{thm:main result heat}
        & Corollary~\ref{cor:regular solution estimate} \\
        \hline
        $h_t \simeq h_x$
        & $\mathcal{O}(h_x^{s-1})$
        & $\mathcal{O}(h_x^{\min\{1,s-1/2\}})$ \\[2pt]
        $h_t \simeq h_x^2$
        & $\mathcal{O}(h_x^{\min\{1,2s-1\}})$
        & $\mathcal{O}(h_x^{\min\{1,2s-1\}})$ \\[2pt]
                $h_t \simeq h_x^{2/(2s-1)}$
        & $\mathcal{O}\!\left(
            h_x^{\min\{2-1/s,\,1/(2s-1)\}}
          \right)$
        & $\mathcal{O}(h_x)$ \\
    \end{tabular}
\end{table}
\end{remark}

After a suitable ordering of degrees of freedom, the arising finite element matrix of the variational problem \eqref{eq:perturbed DVF heat} can be represented via Kronecker products as
\begin{align}
    \mathbf{B}_h = \mathbf{A}_t \otimes \mathbf{M}_x + \mathbf{M}_t \otimes \mathbf{K}_x, \label{eq:heat_kronecker_system}
\end{align}
where
\begin{align}
    \left(\mathbf{A}_t\right)_{i,j} 
    &= (\partial_t \varphi_j, (I + \mathcal{H}_T) \varphi_i )_{L^2(0,T)}, 
    \quad 
    &\left(\mathbf{M}_x\right)_{i,j} 
    &= (\psi_j, \psi_i )_{L^2(\Omega)}, \\
    \left(\mathbf{M}_t\right)_{i,j} 
    &= (\varphi_j, (I + \mathcal{H}_T) \varphi_i)_{L^2(0,T)}, 
    \quad &\left(\mathbf{K}_x\right)_{i,j} 
    &= (\nabla_x \psi_j, \nabla_x \psi_i )_{L^2(\Omega)}
\end{align}
and $X_{h_t} = \text{span}\{\varphi_i\}_{i=1}^{N_t}$, $V_{h_x} = \text{span}\{\psi_i\}_{i=1}^{\dim \left( V_{h_x} \right)}$.

\begin{remark}
    For the use in iterative or factorization based solvers, it is not necessary to form the full Kronecker product matrix explicitly. Instead, one can exploit the Kronecker product structure to perform matrix-vector multiplications efficiently, see e.g. \cite{ZankLangerEfficientDirect2021}. Even though the matrices involving the modified Hilbert transform are dense, it is possible to apply them efficiently using fast transforms which are implemented in the software package \cite{Reichelt_hmod_2025}.
\end{remark}

\subsubsection{BPX Preconditioning}\label{sec:heat:BPX}
We construct a BPX preconditioner in analogy to the ODE case. In time, we assume again a hierarchy of nested finite element spaces
\begin{align}
    X_0 \subset X_1 \subset \ldots \subset X_L = X_{h_t} \subset X_{L+1} \subset \ldots,
\end{align}
where $X_L$ represents the temporal finite element space of the discretization level $L$ in which we want to solve the problem. In space, we consider a single finite element space $V_{h_x} \subset H^1_0(\Omega)$. The norm inducing operator $D \, : \, \mathcal{X} \to \mathcal{X}'$ is now defined via
\begin{align}
    \langle D u, v \rangle_{\mathcal{X}'\times \mathcal{X}} 
    = \langle \partial_t u, \Hmod v \rangle_{[H^{1/2}_{,0}(0,T;L^2(\Omega))]'\times H^{1/2}_{,0}(0,T;L^2(\Omega))} 
    + 
    (\nabla_x u, \nabla_x v)_{L^2(Q)},
\end{align}
or in operator notation
\begin{align}
    D = \Hmod'\partial_t \otimes I_x + I_t \otimes (-\Delta_x),
\end{align}
where $I_t$ and $I_x$ are the identity operators on $L^2(0,T)$ and $L^2(\Omega)$, respectively.
 Denoting the $L^2$-orthogonal projection onto $X_\ell$ by $Q_\ell \, : \, L^2(0,T) \to X_\ell$, the following spectral equivalence for the temporal parts
 \begin{align}
        \Hmod' \partial_t \simeq \sum_{\ell=0}^\infty h_{t,\ell}^{-1} (Q_\ell - Q_{\ell-1}), \quad I_t = \sum_{\ell=0}^\infty (Q_\ell - Q_{\ell-1})
 \end{align}
 holds, where $h_{t,\ell}$ is the temporal mesh size on level $\ell$ and $Q_{-1}=0$.
In the finest temporal basis $\{\varphi_i\}$ let $\mathbf{P}_{-1}=0$ and define the matrices of the projections by
\begin{align}
    (\mathbf{P}_\ell)_{i,j} = (\varphi_j, Q_\ell \varphi_i )_{L^2(0,T)}.
\end{align}
With $\mathbf{A}_\ell=h_{t,\ell}^{-1}\mathbf{M}_x+\mathbf{K}_x$, the multilevel matrix
\begin{align}
    \mathbf{D}_h
    = \sum_{\ell=0}^L (\mathbf{P}_\ell - \mathbf{P}_{\ell-1}) \otimes \mathbf{A}_\ell
    \label{eq:heat_multilevel_matrix}
\end{align}
is spectrally equivalent to the discrete energy norm matrix. Let $\mathbf{M}_{t,\ell}$ be the temporal mass matrix on level $\ell$. The coefficient matrix of $Q_\ell$ is $\mathbf{M}_{t,L}^{-1}\mathbf{P}_\ell$. Since the differences $Q_\ell-Q_{\ell-1}$ project onto mutually orthogonal temporal detail spaces, the inverse of \eqref{eq:heat_multilevel_matrix} is
\begin{align}
    \mathbf{D}_h^{-1}
    = \sum_{\ell=0}^L \bigl[\mathbf{M}_{t,L}^{-1}(\mathbf{P}_\ell-\mathbf{P}_{\ell-1})\mathbf{M}_{t,L}^{-1}\bigr]
      \otimes \mathbf{A}_\ell^{-1}.
\end{align}
Let $\mathbf{R}_\ell$ denote prolongation from $X_\ell$ to $X_L$ and set
$\mathbf{C}_\ell=\mathbf{R}_\ell\mathbf{M}_{t,\ell}^{-1}\mathbf{R}_\ell^\top
=\mathbf{M}_{t,L}^{-1}\mathbf{P}_\ell\mathbf{M}_{t,L}^{-1}$.
Collecting the terms for each temporal level gives
\begin{align}
    \mathbf{D}_h^{-1}
    = \mathbf{M}_{t,L}^{-1}\otimes\mathbf{A}_L^{-1}
      +\sum_{\ell=0}^{L-1}\mathbf{C}_\ell\otimes
       \bigl(\mathbf{A}_\ell^{-1}-\mathbf{A}_{\ell+1}^{-1}\bigr).
    \label{eq:heat_bpx_inverse}
\end{align}
For a spatial eigenvector satisfying $\mathbf{K}_x\mathbf{z}=\mu\mathbf{M}_x\mathbf{z}$, we have $\mathbf{A}_\ell^{-1}\mathbf{M}_x\mathbf{z}=(h_{t,\ell}^{-1}+\mu)^{-1}\mathbf{z}$. Thus the spatial factors in \eqref{eq:heat_bpx_inverse} give precisely the ODE weights in \eqref{eq:ode_bpx_inverse} on each eigenmode.

\begin{remark}
    For large spatial matrices, the application of $\left(h_{t,\ell}^{-1} \mathbf{M}_x + \mathbf{K}_x\right)^{-1}$ can be realized using an iterative solver. As these matrices correspond to spatial elliptic problems with only an additional mass term, there is a wide range of suitable preconditioners available, e.g. multigrid methods or domain decomposition methods \cite{TarekDomainDecompositionMethods2008, TrottenbergOosterleeSchullerMultigrid2001}. 
    If this is done an overall solver with complexity $\mathcal{O}(N_t \log(N_t) N_x)$ can still be achieved, where $N_x = \dim(V_{h_x})$.
\end{remark}

\subsection{Numerical Examples}\label{sec:pde:numerical_examples}

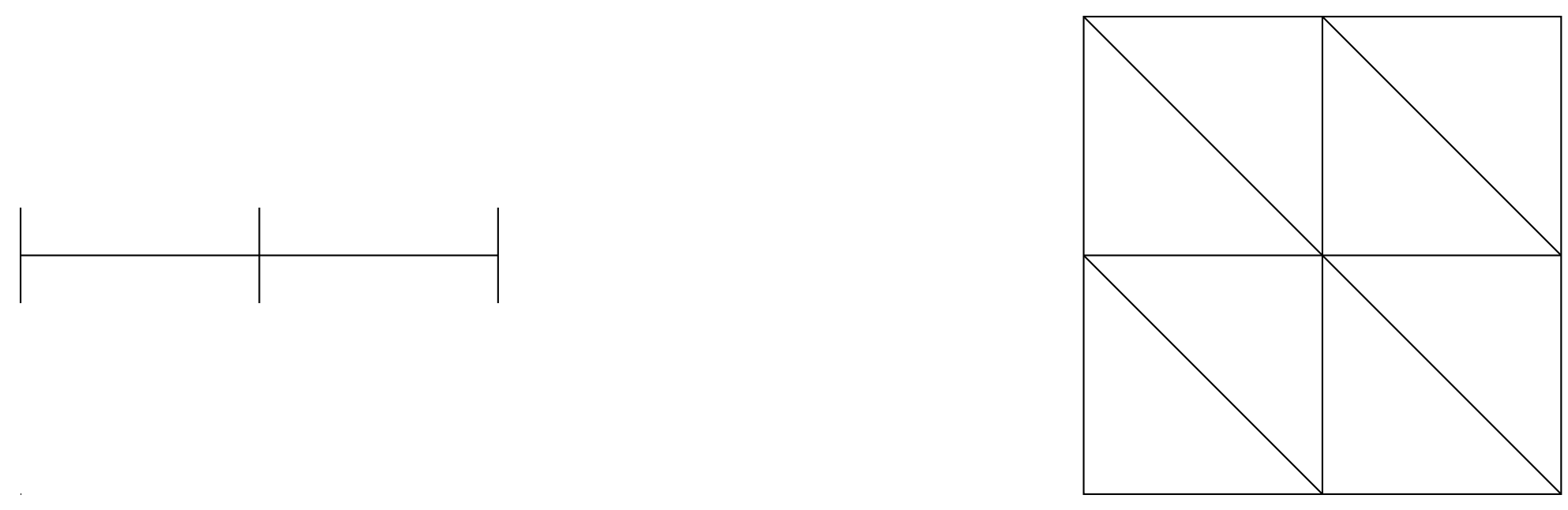
\begin{figure}[h]
    \centering
    \begin{subfigure}[t]{0.45\textwidth}
        \centering
        \begin{tikzpicture}[scale=4]
    \draw (0,0.5) -- (1,0.5); 
    \foreach \x in {0,0.5,1} 
        \draw (\x,0.5-0.1) -- (\x,0.5+0.1); 
    \draw(0.0, 0.0) -- (0.0,0.0); 
\end{tikzpicture}
        \caption{Initial temporal mesh with $N_t = 2$.}
    \end{subfigure}
    \hfill
    \begin{subfigure}[t]{0.45\textwidth}
        \centering
        \begin{tikzpicture}[scale=4]
    \draw (0,0) -- (1,0) -- (1,1) -- (0,1) -- cycle; 
    \draw(0, 0.5) -- (1,0.5); 
    \draw(0.5, 0) -- (0.5,1); 
    \draw(0, 0.5) -- (0.5,0);
    \draw(0.5, 1) -- (1,0.5);
    \draw(0, 1) -- (1, 0); 
\end{tikzpicture}
        \caption{Initial spatial mesh with two intervals in each coordinate direction and eight triangles.}
    \end{subfigure}
    \caption{Initial meshes for the heat equation examples with isotropic scaling.}
    \label{fig:initial_mesh_heat_equation}
\end{figure}

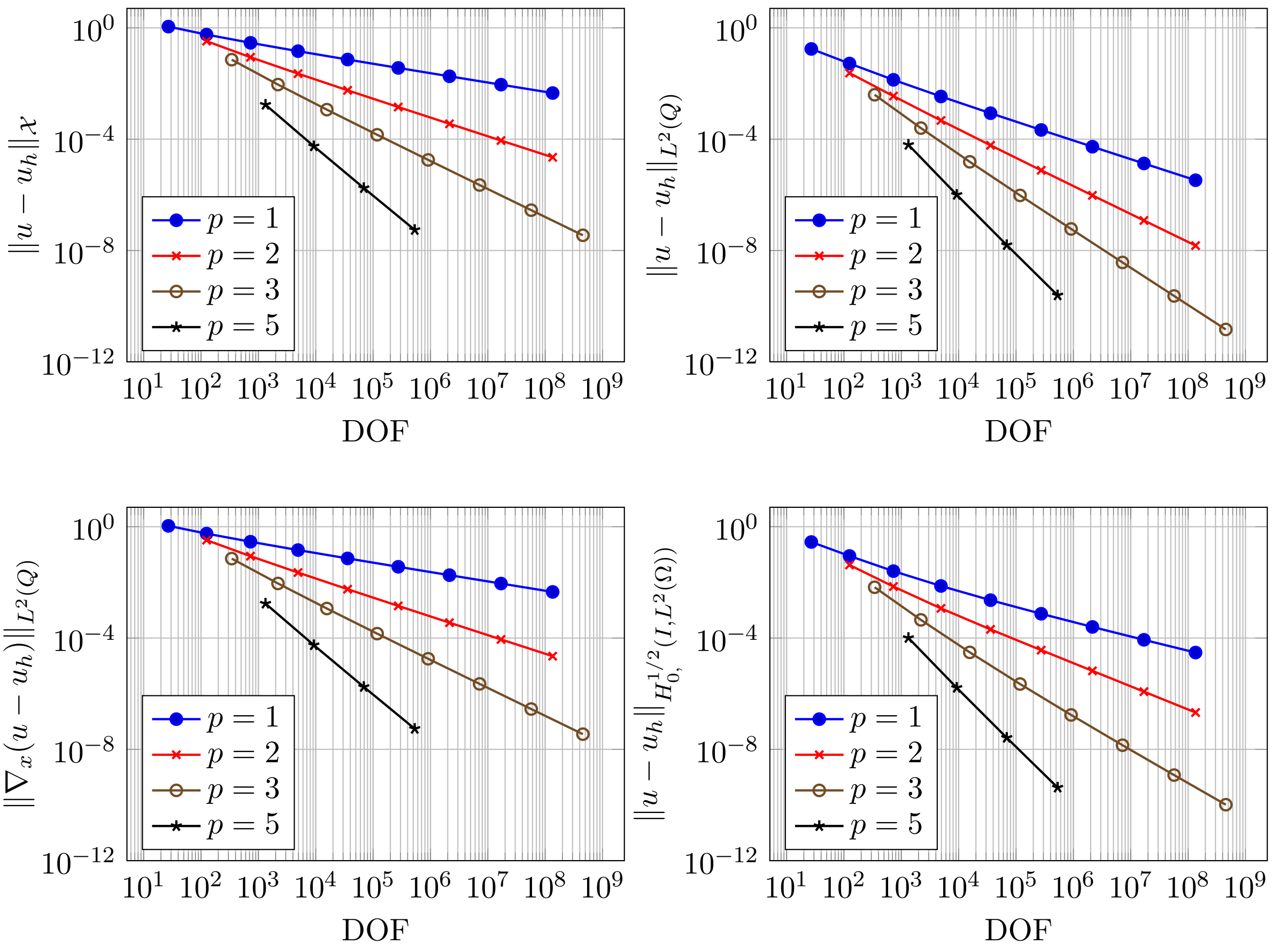
\begin{figure}
    \centering
    \def\csvfolder{simulations/heat/smooth/}
    \begin{tikzpicture}
  \small
  \begin{groupplot}[
    group style={
      group size=2 by 2,       
      horizontal sep=1.8cm,    
      vertical sep=1.8cm,      
    },
    width=0.380\linewidth,      
    height=\HeatSmoothAxisHeight\linewidth,
    scale only axis,
    grid=both,
  ]

    \nextgroupplot[
      xlabel={DOF},
      ylabel={$\|u - u_h\|_{\mathcal{X}}$},
      ylabel near ticks,
      ylabel shift = -7pt,
      ymax = 0.5*1e1,
      ymin = 1e-12,
      ymode=log, xmode=log,
      legend pos=south west,
    ]
    \addplot+[
      thick,
      mark=*,
    ]
    table[
      col sep=comma,
      x expr={0<1e-5 ? \thisrow{dofsQ} : nan},
      y expr={0<1e-5 ? \thisrow{errorEnergy} : nan},
    ]{\csvfolder/results_pol1.csv};
    \addlegendentry{$p=1$}

    \addplot+[
      thick,
      mark=x,
    ]
    table[
      col sep=comma,
      x expr={0<1e-5 ? \thisrow{dofsQ} : nan},
      y expr={0<1e-5 ? \thisrow{errorEnergy} : nan},
    ]{\csvfolder/results_pol2.csv};
    \addlegendentry{$p=2$}

    \addplot+[
      thick,
      mark=o,
    ]
    table[
      col sep=comma,
      x expr={0<1e-5 ? \thisrow{dofsQ} : nan},
      y expr={0<1e-5 ? \thisrow{errorEnergy} : nan},
    ]{\csvfolder/results_pol3.csv};
    \addlegendentry{$p=3$}

    \addplot+[
      thick,
      mark=star,
    ]
    table[
      col sep=comma,
      x expr={0<1e-5 ? \thisrow{dofsQ} : nan},
      y expr={0<1e-5 ? \thisrow{errorEnergy} : nan},
    ]{\csvfolder/results_pol5.csv};
    \addlegendentry{$p=5$}

    \nextgroupplot[
      xlabel={DOF},
      ylabel={$\|u - u_h\|_{L^2(Q)}$},
      ylabel near ticks,
      ylabel shift = -7pt,
      ymax = 0.5*1e1,
      ymin = 1e-12,
      ymode=log, xmode=log,
      legend pos=south west,
    ]
    \addplot+[
      thick,
      mark=*,
    ]
    table[
      col sep=comma,
      x expr={0<1e-5 ? \thisrow{dofsQ} : nan},
      y expr={0<1e-5 ? \thisrow{errorL2} : nan},
    ]{\csvfolder/results_pol1.csv};
    \addlegendentry{$p=1$}

    \addplot+[
      thick,
      mark=x,
    ]
    table[
      col sep=comma,
      x expr={0<1e-5 ? \thisrow{dofsQ} : nan},
      y expr={0<1e-5 ? \thisrow{errorL2} : nan},
    ]{\csvfolder/results_pol2.csv};
    \addlegendentry{$p=2$}

    \addplot+[
      thick,
      mark=o,
    ]
    table[
      col sep=comma,
      x expr={0<1e-5 ? \thisrow{dofsQ} : nan},
      y expr={0<1e-5 ? \thisrow{errorL2} : nan},
    ]{\csvfolder/results_pol3.csv};
    \addlegendentry{$p=3$}

    \addplot+[
      thick,
      mark=star,
    ]
    table[
      col sep=comma,
      x expr={0<1e-5 ? \thisrow{dofsQ} : nan},
      y expr={0<1e-5 ? \thisrow{errorL2} : nan},
    ]{\csvfolder/results_pol5.csv};
    \addlegendentry{$p=5$}

    \nextgroupplot[
      xlabel={DOF},
      ylabel={$\|\nabla_x(u-u_h)\|_{L^2(Q)}$},
      ylabel near ticks,
      ylabel shift = -7pt,
      ymax = 0.5*1e1,
      ymin = 1e-12,
      ymode=log, xmode=log,
      legend pos=south west,
    ]
    \addplot+[
      thick,
      mark=*,
    ]
    table[
      col sep=comma,
      x expr={0<1e-5 ? \thisrow{dofsQ} : nan},
      y expr={0<1e-5 ? \thisrow{errorL2H1} : nan},
    ]{\csvfolder/results_pol1.csv};
    \addlegendentry{$p=1$}

    \addplot+[
      thick,
      mark=x,
    ]
    table[
      col sep=comma,
      x expr={0<1e-5 ? \thisrow{dofsQ} : nan},
      y expr={0<1e-5 ? \thisrow{errorL2H1} : nan},
    ]{\csvfolder/results_pol2.csv};
    \addlegendentry{$p=2$}

    \addplot+[
      thick,
      mark=o,
    ]
    table[
      col sep=comma,
      x expr={0<1e-5 ? \thisrow{dofsQ} : nan},
      y expr={0<1e-5 ? \thisrow{errorL2H1} : nan},
    ]{\csvfolder/results_pol3.csv};
    \addlegendentry{$p=3$}

    \addplot+[
      thick,
      mark=star,
    ]
    table[
      col sep=comma,
      x expr={0<1e-5 ? \thisrow{dofsQ} : nan},
      y expr={0<1e-5 ? \thisrow{errorL2H1} : nan},
    ]{\csvfolder/results_pol5.csv};
    \addlegendentry{$p=5$}

    \nextgroupplot[
      xlabel={DOF},
      ylabel={$\|u - u_h\|_{H^{1/2}_{0,}(I,L^2(\Omega))}$},
      ylabel near ticks,
      ylabel shift = -7pt,
      ymax = 0.5*1e1,
      ymin = 1e-12,
      ymode=log, xmode=log,
      legend pos=south west,
    ]
    \addplot+[
      thick,
      mark=*,
    ]
    table[
      col sep=comma,
      x expr={0<1e-5 ? \thisrow{dofsQ} : nan},
      y expr={0<1e-5 ? \thisrow{errorH12L2} : nan},
    ]{\csvfolder/results_pol1.csv};
    \addlegendentry{$p=1$}

    \addplot+[
      thick,
      mark=x,
    ]
    table[
      col sep=comma,
      x expr={0<1e-5 ? \thisrow{dofsQ} : nan},
      y expr={0<1e-5 ? \thisrow{errorH12L2} : nan},
    ]{\csvfolder/results_pol2.csv};
    \addlegendentry{$p=2$}

    \addplot+[
      thick,
      mark=o,
    ]
    table[
      col sep=comma,
      x expr={0<1e-5 ? \thisrow{dofsQ} : nan},
      y expr={0<1e-5 ? \thisrow{errorH12L2} : nan},
    ]{\csvfolder/results_pol3.csv};
    \addlegendentry{$p=3$}

    \addplot+[
      thick,
      mark=star,
    ]
    table[
      col sep=comma,
      x expr={0<1e-5 ? \thisrow{dofsQ} : nan},
      y expr={0<1e-5 ? \thisrow{errorH12L2} : nan},
    ]{\csvfolder/results_pol5.csv};
    \addlegendentry{$p=5$}

  \end{groupplot}

\end{tikzpicture}
    \caption{Convergence results for the smooth solution $u_r$ for different polynomial degrees $p = p_t = p_x$.}
    \label{fig:pde:regular_convergence}
\end{figure}

To validate the theoretical results, we present numerical examples for the heat equation
in the spatial domain $\Omega = (0,1)^2$ and terminal time $T = 1$. We will consider the regular solution
\begin{align}
    u_r(t,x,y) = \sin(\pi t) \sin(\pi x) \sin(\pi y) \in C^\infty(\overline{Q}),
\end{align}
and the solution with a temporal singularity
\begin{align}
    u_s(t,x,y) = \sqrt{t} \sin(\pi x) \sin(\pi y) \in H^{s}(0,T; C^\infty(\overline{\Omega})),\, s<1.
\end{align}
We discretize the problem in the space-time finite element space $\mathcal{X}_h$ of \eqref{eq:tensorSpace}. The temporal grid has $N_t$ intervals of length $h_t=1/N_t$. On the uniform triangulation of $\Omega$, we denote the spacing along each coordinate axis by $h_x=\sqrt{1/N_x}$. Figure \ref{fig:initial_mesh_heat_equation} shows the initial grids for isotropic refinement.
We again use unrestarted GMRES from a zero initial guess, with the BPX inverse \eqref{eq:heat_bpx_inverse} as a left preconditioner. We initially use an internal tolerance of $10^{-12}$ on the relative preconditioned residual and require $\lVert \underline r_k\rVert_2/\lVert \underline r_0\rVert_2\leq 10^{-10}$, where $\underline r_k$ is the residual of the original linear system. We assemble the spatial finite element matrices using ngsolve \cite{schoberl2014c++} and apply the temporal matrices using the FFT method in \cite{Reichelt_hmod_2025}.
For the regular solution we use $h_t=h_x$ and equal polynomial degrees $p=p_t=p_x$. In Figure \ref{fig:pde:regular_convergence} the energy orders on the finest two grids are $1.00$, $2.00$, $3.00$, and $5.00$ for $p=1,2,3,5$, respectively. The corresponding $L^2(Q)$ orders are $2.00$, $3.00$, $4.00$, and $6.00$. The respective GMRES iterations are found in Figure \ref{fig:pde:regular_gmres}.

\begin{figure}
    \centering
    \def\csvfolder{simulations/heat/singular/}
    \begin{tikzpicture}
  \small
  \begin{groupplot}[
    group style={
      group size=2 by 2,       
      horizontal sep=1.8cm,    
      vertical sep=1.8cm,      
    },
    width=0.385\linewidth,      
    height=\HeatSingularAxisHeight\linewidth,
    scale only axis,
    grid=both,
  ]

    \nextgroupplot[
      xlabel={$N_x$},
      ylabel={$\|u - u_h\|_{\mathcal{X}}$},
      ylabel near ticks,
      ylabel shift = -7pt,
      ymax = 0.5*1e1,
      ymin = 1e-5,
      ymode=log, xmode=log,
      legend pos=south west,
    ]
    \addplot+[
      thick,
      mark=*,
    ]
    table[
      col sep=comma,
      x expr={0<1e-5 ? \thisrow{nx} : nan},
      y expr={0<1e-5 ? \thisrow{errorEnergy} : nan},
    ]{\csvfolder/results_pol1.csv};
    \addlegendentry{$h_t=h_x$}

    \addplot+[
      thick,
      mark=x,
    ]
    table[
      col sep=comma,
      x expr={0<1e-5 ? \thisrow{nx} : nan},
      y expr={0<1e-5 ? \thisrow{errorEnergy} : nan},
    ]{\csvfolder/results_pol1_parabolic.csv};
    \addlegendentry{$h_t=h_x^2$}

    \nextgroupplot[
      xlabel={$N_x$},
      ylabel={$\|u - u_h\|_{L^2(Q)}$},
      ylabel near ticks,
      ylabel shift = -7pt,
      ymax = 0.5*1e1,
      ymin = 1e-5,
      ymode=log, xmode=log,
      legend pos=south west,
    ]
    \addplot+[
      thick,
      mark=*,
    ]
    table[
      col sep=comma,
      x expr={0<1e-5 ? \thisrow{nx} : nan},
      y expr={0<1e-5 ? \thisrow{errorL2} : nan},
    ]{\csvfolder/results_pol1.csv};
    \addlegendentry{$h_t=h_x$}

    \addplot+[
      thick,
      mark=x,
    ]
    table[
      col sep=comma,
      x expr={0<1e-5 ? \thisrow{nx} : nan},
      y expr={0<1e-5 ? \thisrow{errorL2} : nan},
    ]{\csvfolder/results_pol1_parabolic.csv};
    \addlegendentry{$h_t=h_x^2$}

    \nextgroupplot[
      xlabel={$N_x$},
      ylabel={$\|\nabla_x(u-u_h)\|_{L^2(Q)}$},
      ylabel near ticks,
      ylabel shift = -7pt,
      ymax = 0.5*1e1,
      ymin = 1e-5,
      ymode=log, xmode=log,
      legend pos=south west,
    ]
    \addplot+[
      thick,
      mark=*,
    ]
    table[
      col sep=comma,
      x expr={0<1e-5 ? \thisrow{nx} : nan},
      y expr={0<1e-5 ? \thisrow{errorL2H1} : nan},
    ]{\csvfolder/results_pol1.csv};
    \addlegendentry{$h_t=h_x$}

    \addplot+[
      thick,
      mark=x,
    ]
    table[
      col sep=comma,
      x expr={0<1e-5 ? \thisrow{nx} : nan},
      y expr={0<1e-5 ? \thisrow{errorL2H1} : nan},
    ]{\csvfolder/results_pol1_parabolic.csv};
    \addlegendentry{$h_t=h_x^2$}

    \nextgroupplot[
      xlabel={$N_x$},
      ylabel={$\|u - u_h\|_{H^{1/2}_{0,}(I,L^2(\Omega))}$},
      ylabel near ticks,
      ylabel shift = -7pt,
      ymax = 0.5*1e1,
      ymin = 1e-5,
      ymode=log, xmode=log,
      legend pos=south west,
    ]
    \addplot+[
      thick,
      mark=*,
    ]
    table[
      col sep=comma,
      x expr={0<1e-5 ? \thisrow{nx} : nan},
      y expr={0<1e-5 ? \thisrow{errorH12L2} : nan},
    ]{\csvfolder/results_pol1.csv};
    \addlegendentry{$h_t = h_x$}

    \addplot+[
      thick,
      mark=x,
    ]
    table[
      col sep=comma,
      x expr={0<1e-5 ? \thisrow{nx} : nan},
      y expr={0<1e-5 ? \thisrow{errorH12L2} : nan},
    ]{\csvfolder/results_pol1_parabolic.csv};
    \addlegendentry{$h_t = h_x^2$}

  \end{groupplot}

\end{tikzpicture}
    \caption{Convergence results for the singular solution $u_s$ for $p_t =1$ and $h_t = h_x$ or $h_t = h_x^2$.}
    \label{fig:pde_singular_convergence}
\end{figure}
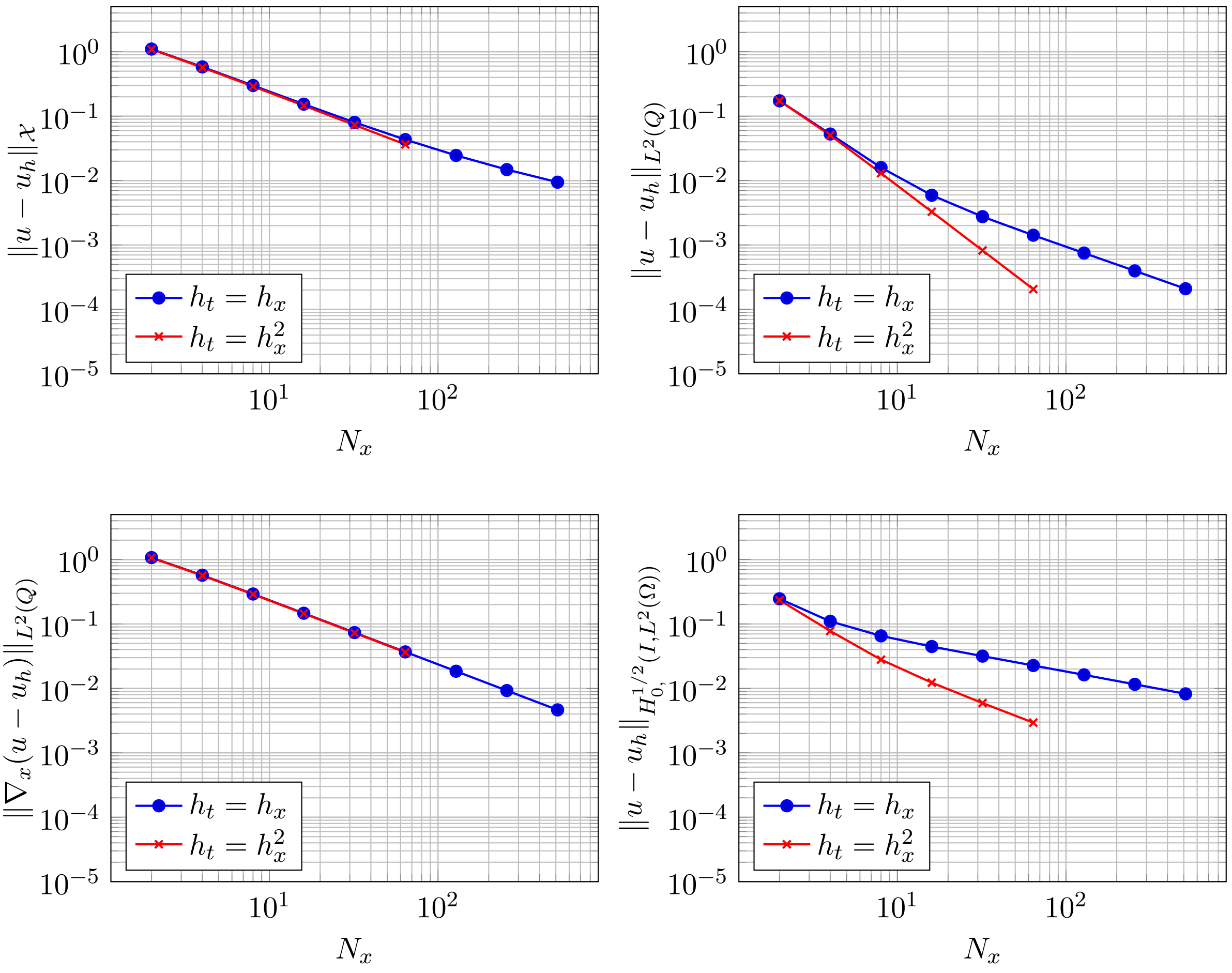

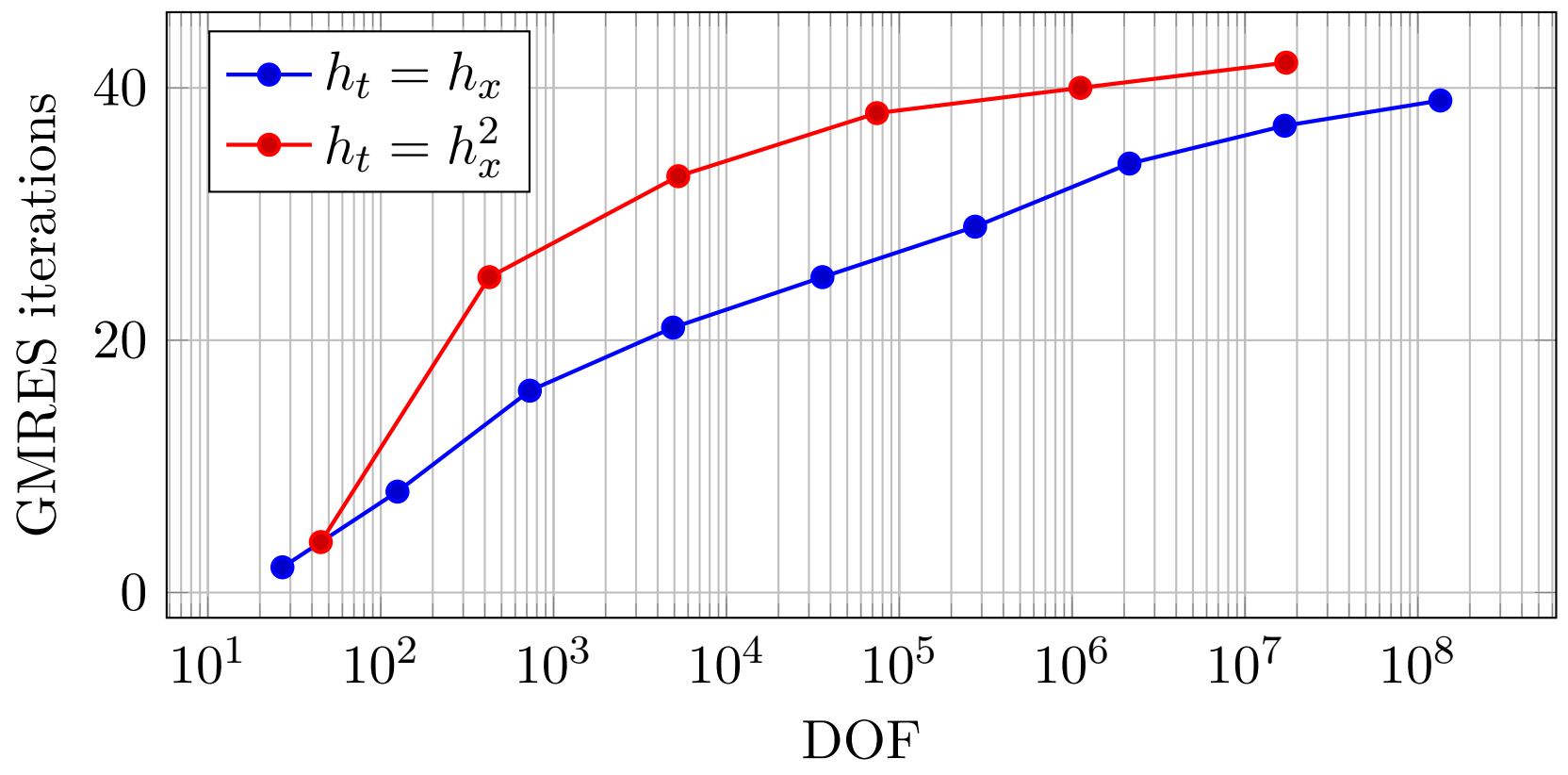
\begin{figure}
    \centering
    \def\csvfolder{simulations/heat/singular/}
    \begin{tikzpicture}
    \small
    \begin{semilogxaxis}[
      xlabel={DOF},
      ylabel={GMRES iterations},
      grid=both,
      width=0.7\linewidth,
      height=\HeatGMRESAxisHeight\linewidth,
      legend pos=north west,
    ]
    \addplot+[
      thick,
      mark=*,
    ]
    table[
      col sep=comma,
      x=dofsQ,
      y=gmresIter,
    ]{\csvfolder/results_pol1.csv};
    \addlegendentry{$h_t=h_x$}
    \addplot+[
      thick,
      mark=*,
    ]
    table[
      col sep=comma,
      x=dofsQ,
      y=gmresIter,
    ]{\csvfolder/results_pol1_parabolic.csv};
    \addlegendentry{$h_t=h_x^2$}
    
    \end{semilogxaxis}
\end{tikzpicture}
    \caption{GMRES iterations for the singular solution $u_s$ for $p_t =1$ and $h_t = h_x$ or $h_t = h_x^2$.}
    \label{fig:pde_singular_gmres}
\end{figure}

For the singular solution we fix $p_t=p_x=1$ and compare isotropic refinement $h_t=h_x$ with parabolic refinement $h_t=h_x^2$. The energy orders on the finest two grids are $0.650$ and $1.000$, respectively, see Figure \ref{fig:pde_singular_convergence}, which are explained by the theoretical estimates in Table \ref{tab:comparisonScaling} . The respective GMRES iterations are depicted in Figure \ref{fig:pde_singular_gmres}.

\begin{remark}
    One can also extend this approach to variable coefficients or non-linear problems. Examples for the ODE can be found in \cite{Reichelt_hmod_2025}.
\end{remark}

\begin{figure}
    \centering
    \def\csvfolder{simulations/heat/smooth/}
    \begin{tikzpicture}
    \small
    \begin{semilogxaxis}[
      xlabel={DOF},
      ylabel={GMRES iterations},
      grid=both,
      width=0.7\linewidth,
      height=\HeatGMRESAxisHeight\linewidth,
      legend pos=north west,
    ]
    \addplot+[
      thick,
      mark=*,
    ]
    table[
      col sep=comma,
      x=dofsQ,
      y=gmresIter,
    ]{\csvfolder/results_pol1.csv};
    \addlegendentry{$p=1$}
    \addplot+[
      thick,
      mark=*,
    ]
    table[
      col sep=comma,
      x=dofsQ,
      y=gmresIter,
    ]{\csvfolder/results_pol2.csv};
    \addlegendentry{$p=2$}
    \addplot+[
      thick,
      mark=*,
    ]
    table[
      col sep=comma,
      x=dofsQ,
      y=gmresIter,
    ]{\csvfolder/results_pol3.csv};
    \addlegendentry{$p=3$}
    \addplot+[
      thick,
      mark=*,
    ]
    table[
      col sep=comma,
      x=dofsQ,
      y=gmresIter,
    ]{\csvfolder/results_pol5.csv};
    \addlegendentry{$p=5$}
    \end{semilogxaxis}
\end{tikzpicture}
    \caption{GMRES iterations for the smooth solution $u_r$ for different polynomial degrees $p = p_t = p_x$.}
    \label{fig:pde:regular_gmres}
\end{figure}
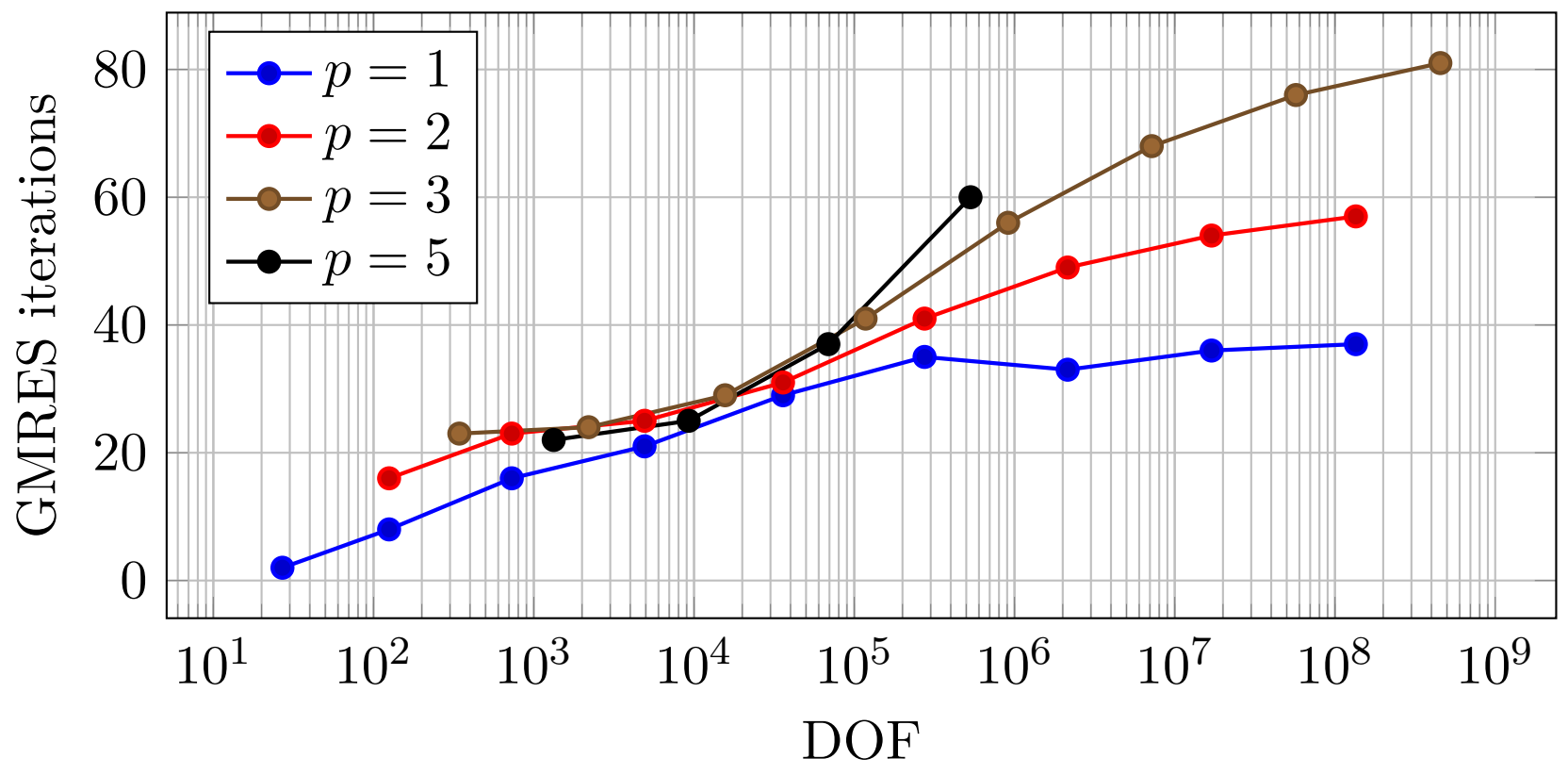

\section{Conclusions and Outlook}
We have presented an operator-theoretic framework based on the polar
decomposition of linear operators for the construction of coercive
variational formulations of non-coercive problems. Its application
to the heat equation reveals connections with several established
approaches. In particular, the unitary factor of the polar decomposition
serves as a transformation operator, providing a systematic construction
of the type of transformations considered in the literature. Its
analytical characterization offers insight into the structure of the
method and facilitates its analysis, while practical implementation
requires an approximation of the transformation. We have analyzed the
resulting perturbation of the bilinear form and established nearly
quasi-optimal error estimates. Moreover, the framework enables the
construction of effective preconditioners for the resulting discrete
systems. Numerical experiments support the theoretical findings.

The scope of this approach extends beyond the heat equation, owing to
the generality of the underlying polar decomposition. Initial
investigations for the wave equation have already been undertaken
\cite{hoonhout2025paving}, while application to the Stokes equations
provides a further direction for research. Related transformation-based
concepts have also been developed for nonlinear operators, notably
the notion of mapped coercivity introduced by Becker and Braack
\cite{BeckerBraak}. Extending the present framework to these settings
and examining its relationship with nonlinear mapped coercivity
remain subjects for future work.

\textbf{Declaration of AI use.} OpenAI’s ChatGPT 5.6-6.2 Sol/Astra was used during the preparation of this manuscript to assist with computations and final polishing (checking mathematical arguments, language revision, and typesetting). The authors independently verified the mathematical content and take full responsibility for the manuscript.

\bibliographystyle{plain}
\bibliography{references_submission}
\end{document}